\documentclass[11pt,oneside,a4paper]{article}
\usepackage{blindtext}
\usepackage{amsmath,amsfonts,amssymb,amsthm,eucal}
\usepackage{microtype}
\usepackage[colorlinks, citecolor=blue, linkcolor=blue, urlcolor=Maroon, filecolor=Maroon, linktocpage=true]{hyperref}
\usepackage[usenames,dvipsnames,svgnames,table]{xcolor}
\usepackage[normalem]{ulem}
\usepackage[T1]{fontenc}
\usepackage{enumerate}
\usepackage{bbm}
\usepackage{tikz}
\usepackage{tikz-cd}
\usepackage[font=footnotesize,labelfont=bf]{caption}
\usepackage{subcaption}
\usepackage{framed}
\usepackage[nottoc,notlot,notlof]{tocbibind}
\usepackage{tocloft}
\usepackage{morefloats}
\usepackage{float}
\usepackage[nottoc,notlot,notlof]{tocbibind}
\usepackage{tocloft}
\usepackage[backgroundcolor=LightBlue,linecolor=LightBlue]{todonotes}

\newcommand\re[1]{(\ref{#1})}
\newcommand{\D}{\mathrm{d}}
\theoremstyle{definition}
\newtheorem{Th}{Theorem}[section]
\newtheorem{Ex}[Th]{Example}
\newtheorem{Prop}[Th]{Proposition}

\newtheorem{Cor}[Th]{Corollary}
\newtheorem{Lem}[Th]{Lemma}
\newtheorem{Def}[Th]{Definition}

\newtheorem{open}[Th]{Open problem}

\theoremstyle{definition}
\newtheorem{rem}[Th]{Remark}
\usepackage[a4paper,
	hcentering,
	text={453pt,686pt},
]{geometry}
\allowdisplaybreaks
\catcode`@=11
\@addtoreset{equation}{section}
\catcode`@=12
\begin{document}
\begin{center}\noindent
\textbf{\Large Limit geometry of complete projective special real manifolds}\\[2em]
{\fontseries{m}\fontfamily{cmss}\selectfont \large David\ Lindemann}\\[1em] 
{\small
Department of Mathematics, Aarhus University\\
Ny Munkegade 118, Bldg 1530, DK-8000 Aarhus C, Denmark\\
\texttt{david.lindemann@math.au.dk}
}
\end{center}
\vspace{1em}
\begin{abstract}
\noindent
We study the limit geometry of complete projective special real manifolds. By limit geometry we mean the limit of the evolution of the defining polynomial and the centro-affine fundamental form along certain curves that leave every compact subset of the initial complete projective special real manifold. We obtain a list of possible limit geometries, which are themselves complete projective special real manifolds, and find a lower bound for the dimension of their respective symmetry groups. We further show that if the initial manifold has regular boundary behaviour, every possible limit geometry is isomorphic to $\mathbb{R}_{>0}\ltimes\mathbb{R}^{n-1}$.
\end{abstract}
\tableofcontents
\section{Introduction and main results}
	In this work we study a certain notion of limit geometry of connected geodesically complete projective special real manifolds. A \textbf{projective special real manifold} (short: \textbf{PSR manifold}) $\mathcal{H}$ of dimension $\dim(\mathcal{H})=n$ is a smooth hypersurface in $\mathbb{R}^{n+1}$ that is contained in the level set $\{h=1\}$ of a hyperbolic cubic homogeneous polynomial $h:\mathbb{R}^{n+1}\to \mathbb{R}$, such that $\mathcal{H}\subset\{h=1\}$ consists only of hyperbolic points of $h$. A point $p\in\{h>0\}$ is called hyperbolic point of $h$ if the negative Hessian $-\partial^2 h_p$ has Lorentz signature, and a homogeneous polynomial of degree at least $2$ admitting a hyperbolic point is called hyperbolic homogeneous polynomial. Every PSR manifold is, equipped with its centro-affine fundamental form $g_\mathcal{H}=-\frac{1}{3}\partial^2 h|_{T\mathcal{H}\times T\mathcal{H}}$, a Riemannian manifold. This follows easily from Euler's Homogeneous Function Theorem. It was shown in \cite{CNS} that a PSR manifold $(\mathcal{H},g_\mathcal{H})$ is geodesically complete if and only if it is closed in its ambient space $\mathbb{R}^{n+1}$. Two PSR manifolds of the same dimension are called \textbf{equivalent} if they are related by a linear transformation of the ambient space, meaning in particular that their defining polynomials are also related by said transformation. In the following, we will abbreviate \textbf{closed connected PSR} manifolds to \textbf{CCPSR} manifolds. A connected PSR manifold $\mathcal{H}\subset\{h=1\}$ is called \textbf{maximal} if it consists of one or more connected components of the set $\{h=1\}\cap\{\text{hyperbolic points of }h\}$. Note that maximality does not imply closedness. In order to give meaning to the term \textbf{limit geometry} we will need the main results of \cite{L2}. In the following, $\langle\cdot,\cdot\rangle$ and $\|\cdot\|$ denote the standard Euclidean scalar product in chosen linear coordinates on $\mathbb{R}^N$ for fitting $N$ and its induced norm, respectively. Furthermore we will throughout this work identify homogeneous polynomials with their associated multi-linear forms so that e.g. for a cubic homogeneous polynomial $h:\mathbb{R}^{n+1}\to\mathbb{R}$ we have $\D h_p(v)=3h(p,p,v)$ for all $p,v\in\mathbb{R}^{n+1}$.

	\begin{Prop}\label{prop_standard_form}
		Let $\mathcal{H}\subset\{h=1\}\subset\mathbb{R}^{n+1}$ be an $n$-dimensional CCPSR manifold and let $\left(\begin{smallmatrix}x\\y\end{smallmatrix}\right)=(x,y_1,\ldots,y_n)^\mathrm{T}$ denote linear coordinates on the ambient space $\mathbb{R}^{n+1}$. Then there exists a smooth map
			\begin{equation*}
				A:\mathcal{H}\to \mathrm{GL}(n+1),
			\end{equation*}
		describing a linear change of coordinates in the ambient space $\mathbb{R}^{n+1}$ depending on $p\in\mathcal{H}$, such that
			\begin{enumerate}[(i)]
				\item\label{prop_standard_form_i} $(A(p)^*h)\left(\left(\begin{smallmatrix}x\\y\end{smallmatrix}\right)\right) = x^3 -x\langle y,y\rangle + P_3(y)$,
				\item $A(p)\left(\begin{smallmatrix}1\\0\end{smallmatrix}\right)=p$,
			\end{enumerate}
		for all $p\in\mathcal{H}$, where $P_3:\mathbb{R}^n\to\mathbb{R}$ is a cubic homogeneous polynomial. The most general form of the map $A$ is given by
			\begin{equation}
				A(p)=\left(
				\begin{array}{c|c}
				p_x &  -\left.\frac{\partial_yh}{\partial_xh}\right|_p\circ \underset{}{E(p)} \\ \hline
				p_y  & \overset{}{E(p)}
				\end{array}\right)\label{eqn_pmoving_trafo_explicit}
			\end{equation}
		for all $p\in\mathcal{H}$, where $p=\left(\begin{smallmatrix}p_x\\p_y\end{smallmatrix}\right)$ with $p_x=x(p)$ and $p_y=(y_1(p),\ldots,y_n(p))^\mathrm{T}$, $\partial_yh=(\partial_{y_1} h,\ldots,\partial_{y_n} h)$, and $E:\mathcal{H}\to \mathrm{GL}(n)$ is required to diagonalise the positive definite symmetric bilinear form
			\begin{equation}
				\mathbb{R}^n\times\mathbb{R}^n\ni(v,w)\mapsto -\frac{1}{2}\partial^2 h_p\left(\left(
				\begin{matrix}
				-\left.\frac{\partial_yh}{\partial_xh}\right|_p(v)\\
				v
				\end{matrix}
				\right),\left(
				\begin{matrix}
				-\left.\frac{\partial_yh}{\partial_xh}\right|_p(w)\\
				w
				\end{matrix}
				\right)\right)\in\mathbb{R}\label{eqn_symbilform_to_diag}
			\end{equation}
		for all $p\in\mathcal{H}$.
		\begin{proof}
			This is a special case of \cite[Prop.\,3.1]{L2} for CCPSR manifolds.
		\end{proof}
	\end{Prop}
	
	In order to proceed we will need and consequently make heavy use of the main result of \cite{L2}.
	
	\begin{Th}\label{thm_Cn}
		Let $(x,y_1,\ldots,y_n)^\mathrm{T}$ denote linear coordinates on $\mathbb{R}^{n+1}$ and let $h:\mathbb{R}^{n+1}\to\mathbb{R}$ be a homogeneous cubic polynomial of the form $h=x^3-x\langle y,y\rangle + P_3(y)$ as in Prop. \ref{prop_standard_form} \eqref{prop_standard_form_i}. Then the connected component $\mathcal{H}\subset\{h=1\}$ containing the point $\left(\begin{smallmatrix}x\\y\end{smallmatrix}\right)=\left(\begin{smallmatrix}1\\0\end{smallmatrix}\right)$ is a CCPSR manifold if and only if
			\begin{equation}
				\max\limits_{\|y\|=1} P_3(y)\leq \frac{2}{3\sqrt{3}}.\label{eqn_P3_max_condition}
			\end{equation}
		Hence, the affine subset
			\begin{equation}
				\mathcal{C}_n:=\left.\left\{x^3-x\langle y,y\rangle + P_3(y)\ \right|\ P_3\text{ fulfils \eqref{eqn_P3_max_condition}}\right\}\subset \mathrm{Sym}^3(\mathbb{R}^{n+1})^*\label{eqn_Cn_gen_set}
			\end{equation}
		is a compact convex generating set of the set of $n$-dimensional CCPSR manifolds.
		\begin{proof}
			This is \cite[Thm\,1.1]{L2} and \cite[Thm\,4.15]{L2}.
		\end{proof}
	\end{Th}
	
	Note that the compactness of $\mathcal{C}_n$ is to be understood with respect to the subspace topology of $\mathrm{Sym}^3(\mathbb{R}^{n+1})^*$ which is induced by the real vector space structure of $\mathrm{Sym}^3(\mathbb{R}^{n+1})^*$. Proposition \ref{prop_standard_form} and Theorem \ref{thm_Cn} imply that we can, in dependence of a reference point and the freedom of choice in the transformation \eqref{eqn_pmoving_trafo_explicit} which essentially boils down to choosing an element in $\mathrm{O}(n)$, assume without loss of generality that the defining polynomial of a given CCPSR manifold $\mathcal{H}\subset\{h=1\}$ is of the form
		\begin{equation}\label{eqn_h_general_form}
			h=h\left(\left(\begin{smallmatrix}x\\y\end{smallmatrix}\right)\right)=x^3-x\langle y,y\rangle + P_3(y),
		\end{equation}
	$P_3$ fulfils the maximality condition \eqref{eqn_P3_max_condition}, and $\mathcal{H}$ is precisely the connected component of $\{h=1\}$ that contains the point $\left(\begin{smallmatrix}x\\ y\end{smallmatrix}\right)=\left(\begin{smallmatrix}1\\ 0\end{smallmatrix}\right)$. We will then say that $\mathcal{H}$ and $h$ are in \textbf{standard form}. If $\mathcal{H}$ is the connected component of $\{h=1\}\cap\{\text{hyperbolic points of }h\}$, $h$ of the form \eqref{eqn_h_general_form}, containing the point $\left(\begin{smallmatrix}x\\ y\end{smallmatrix}\right)=\left(\begin{smallmatrix}1\\ 0\end{smallmatrix}\right)$, independent of whether or not $\mathcal{H}$ is closed in its ambient space, we will also say that $\mathcal{H}$ and $h$ are in standard form. Note that the term $P_3$ is never uniquely determined by $\mathcal{H}$ \cite{L2}.
	
	Another necessary concept in order to define the notion for limit geometry of PSR manifolds in our sense is borrowed from the theory of the K\"ahler-Ricci flow on cohomology classes. Let $X$ be a complex 3-dimensional compact K\"ahler manifold and define on the real $(1,1)$-cohomology of $X$ a homogeneous cubic polynomial via
		\begin{equation}
			h:H^{1,1}(X,\mathbb{R})\to\mathbb{R},\quad [\omega]\mapsto \int_X \omega^3.\label{eqn_h_KR_flow}
		\end{equation}
	Let $\mathcal{K}\subset H^{1,1}(X,\mathbb{R})$ denote the K\"ahler cone of $X$. It then follows from the Hodge-Riemann bilinear relations that $\mathcal{H}:=\{h=1\}\cap \mathcal{K}$ is a PSR manifold of dimension $h^{1,1}(X)$ \cite[Chap.\,5,\; Sec.\,6]{We}. One might also consider the index cone $\mathcal{W}\subset H^{1,1}(X,\mathbb{R})$ as in \cite{Wi} which consists of all hyperbolic points of $h$. Then by definition $\{h=1\}\cap\mathcal{W}$ is also a PSR manifold containing $\mathcal{H}$. Note that in general $\mathcal{K}$ and $\mathcal{W}$ need not coincide. This is a consequence of the main result of \cite{DP} which implies that $[\omega]\in\mathcal{W}$ is contained in $\mathcal{K}$ if and only if $[\omega]$ is numerically positive on all irreducible analytic sets in $X$. Note that in the setting of the geometry of K\"ahler cones, completeness questions have been addressed from that point of view in \cite{M}. Now consider the K\"ahler-Ricci flow equation on the level of cohomology,
		\begin{equation}
			\partial_t[\omega_t]=2\pi c_1(X),\label{eqn_KRflow_classes}
		\end{equation}
	where $c_1(X)$ denotes the first Chern class of $X$. The above equation is obtained by considering the cohomology classes of the K\"ahler-Ricci flow equation $\partial_t \omega_t=\mathrm{Ric}^{\omega_t}$ on $X$, for a reference text see e.g. \cite{T}. For initial datum $[\omega_0]\in\mathcal{K}$, the maximal connected solution of \eqref{eqn_KRflow_classes} is given by
		\begin{equation*}
			[\omega_t]:(T_-,T_+)\to\mathcal{K},\quad t\mapsto [\omega_0] + 2\pi t c_1(X),
		\end{equation*}
	where $T_-=\inf\{t<0\ |\ [\omega_0] + 2\pi s c_1(X)\in\mathcal{K}\ \forall s\in(t,0)\}$ and $T_+=\sup\{t>0\ |\ [\omega_0] + 2\pi s c_1(X)\in\mathcal{K}\ \forall s\in(0,t)\}$. If we view $H^{1,1}(X,\mathbb{R})$ as real vector space and $\mathcal{K}$ as an open cone therein, the above solution $[\omega_t]$ is thus simply an affine line segment that is an integral curve of $2\pi c_1(X)$ viewed as a constant vector field on $H^{1,1}(X,\mathbb{R})$. Save for the trivial case $c_1(X)=0$, the K\"ahler-Ricci flow on classes is never volume preserving. This is however a necessity for our endeavour as we want to study the evolution of the centro-affine fundamental form of a PSR manifold along certain curves contained in said manifold. To solve this issue we first switch back to our language and consider a connected PSR manifold $\mathcal{H}\subset\{h=1\}\subset\mathbb{R}^{n+1}$ and a given constant non-vanishing vector field $V\in\mathfrak{X}(\mathbb{R}^{n+1})$ on the ambient space. The picture analogous to the K\"ahler-Ricci flow on classes in our setting would now be to study integral curves $\gamma:I\to\mathbb{R}^{n+1}$, $0\in I$, of $V$ restricted to the cone $\mathbb{R}_{>0}\cdot \mathcal{H}$ with $\gamma(0)\in\mathcal{H}$. As for the K\"ahler-Ricci flow on classes, it is easy to check that then $h(\gamma(t))=1$ only at at most two values for $t$. In order to fix this, we instead consider integral curves of the \textbf{central projection} of $V$ to $T\mathcal{H}$, that is the unique vector field $\widetilde{V}\in\mathfrak{X}(\mathcal{H})$, such that
		\begin{equation*}
			V=\widetilde{V}+f\xi
		\end{equation*}
	along $\mathcal{H}$, where $\xi\in\mathfrak{X}(\mathbb{R}^{n+1})$ is the position vector field and $f\in C^\infty(\mathcal{H})$. One can check that, in fact, $f=\D h(V)$, which follows from Euler's Homogeneous Function Theorem. Each integral curve of $\widetilde{V}$ is contained in subsets of $\mathcal{H}$ of the form
		\begin{equation}
			\mathcal{H}\cap\mathrm{span}_{\mathbb{R}}\{p,\widetilde{V}_p\}\label{eqn_set_containing_considered_curves}
		\end{equation}
	for some $p\in\mathcal{H}$. Since we want to study the evolution of the centro-affine fundamental form of $\mathcal{H}$ along integral curves that leave every compact subset, we might ignore the speed of the integral curves of $\widetilde{V}$ and, hence, make the following ansatz. For $p\in\mathcal{H}$ and $v\in T_p\mathcal{H}\setminus\{0\}$, let $\gamma:(R_-,R_+)\to\mathcal{H}$ for $R_\pm>0$ be a curve with nowhere vanishing velocity with image precisely the set $\mathcal{H}\cap\mathrm{span}_{\mathbb{R}}\{p,v\}$ and $\gamma(0)=p$. How can we, in a realistic way, calculate some sort of limit of the centro-affine fundamental form along $\gamma$, that is $g_{\mathcal{H}}|_{\gamma(t)}$, as $t\to R_-$, respectively $t\to R_+$? Note that this problem is equivalent to just calculating one limit for all such curves $\gamma$, as $t\to R_+$, by simply reversing $\gamma$ for the other limit and up to a sign-change interchanging $R_-$ and $R_+$. Hence, we can without loss of generality assume that $\gamma$ is defined on $[0,R)$ for some $R>0$. In order to tackle that problem we use Proposition \ref{prop_standard_form}, by which we can without loss of generality assume that our given CCPSR manifold $\mathcal{H}\subset\{h=1\}$ is of standard form and $p=\left(\begin{smallmatrix}1\\0\end{smallmatrix}\right)$. We can also set $v$ to be a choice of a coordinate vector, which we will in fact do during the proof of our main result. We can now, at least in theory, calculate for all $t\in[0,R)$ a standard form of the defining polynomial $h$ for the reference point, which essentially means choosing $A(\gamma(t))$ as in \eqref{eqn_pmoving_trafo_explicit}, calculating the $P_3$-part $P_3(y)=P_3(y)(t)$ of $A(\gamma(t))^* h$ for all $t\in[0,R)$, and then calculate the limit of
		\begin{equation}
			h|_t:=A(\gamma(t))^* h=x^3-x\langle y,y\rangle + P_3(y)|_t\label{eqn_h_t}
		\end{equation}
	as $t\to R$, where we view $h|_t$ as a curve in $\mathcal{C}_n$ as in Theorem \ref{thm_Cn}. The main difficulty is showing that such a limit is well defined, which in practice means making a good choice for the transformations $A(\gamma(t))$. Note that for every allowed $t\in[0,R)$, this defines a CCPSR manifold $\mathcal{H}_t$ together with a choice of linear coordinates on the ambient space which is equivalent to $\mathcal{H}=\mathcal{H}_0$. Only for the limit $\mathcal{H}_R$, assuming it is well defined, can we expect to obtain a CCPSR manifold that is not equivalent to the initial $\mathcal{H}$.
	
	A good way to visualise limit geometries of CCPSR manifolds is as follows. In a slightly more general setting it was shown in \cite[Lem.\,1.14]{CNS} that for every CCPSR manifold $\mathcal{H}\subset\{h=1\}$ of dimension $n$ and $p\in\mathcal{H}$ arbitrary, the intersection of the affinely embedded tangent space $p+T_p\mathcal{H}$ with the cone spanned by $\mathcal{H}$, that is $\mathbb{R}_{>0}\cdot\mathcal{H}$, is relatively precompact in $\mathbb{R}^{n+1}$. In \cite{L1,L2} with the assumption that $\mathcal{H}$ is in standard form and $p=(1,0,\ldots,0)^\mathrm{T}$, this set is called $\mathrm{dom}(\mathcal{H})$ and is considered as an open subset of $\mathbb{R}^n$. This makes sense since we can in this case canonically identify $p+T_p\mathcal{H}$ with $p+(\{0\}\times\mathbb{R}^n)$ and obtain linear coordinates by considering only the $y$-coordinates in \eqref{eqn_h_general_form}. The curves with nonvanishing speed contained in a set of the form \eqref{eqn_set_containing_considered_curves} can be, after a possible change of reference point for the standard form, a rotation in the $y$-coordinates, and a reparametrisation of the curve itself, assumed to be of the form
		\begin{equation}
			\gamma:[0,R)\mapsto\mathcal{H},\quad t\mapsto \beta^{-\frac{1}{3}}\cdot(1,0,\ldots,0,t)^\mathrm{T},\label{eqn_curve_for_limitgeo}
		\end{equation}
	where $\beta=h((1,0,\ldots,0,t)^\mathrm{T})$. This means that $\beta$ is the value of $h$ along a ray starting from the point $p=(1,0,\ldots,0)^\mathrm{T}$ in the affinely embedded set $p+\mathrm{dom}(\mathcal{H})$ in direction of the last $y$-coordinate $y_{n}$. We will use the above notation with the symbol $\beta$ throughout this work. Note that we could have of course chosen another $y$-coordinate to be the direction of $\gamma$. Now, at least in simple cases in dimension $n=2$ and $n=3$, one can use plotting software to see how the boundary of $\mathrm{dom}(\mathcal{H}_t)$, that is $\partial\mathrm{dom}(\mathcal{H}_t)$, changes as $t$ varies, see Figure \ref{fig_mot_example} in Example \ref{ex_motivating}.
	
	We now have all concepts at hand to give a rigorous definition of the terms limit polynomial and limit geometry.
	
	\begin{Def}\label{def_limit_poly_limit_behaviour}
		Assuming that $\overline{h}:=\lim\limits_{t\to R} h|_t$ exists, we will call $\overline{h}$ a \textbf{limit polynomial}. The corresponding CCPSR manifold $\overline{\mathcal{H}}$ given by the connected component of $\{\overline{h}=1\}$ containing the point $p=\left(\begin{smallmatrix}1\\0\end{smallmatrix}\right)$ will be called \textbf{limit geometry} of the initial CCPSR manifold $\mathcal{H}$.
	\end{Def}

	When discussing limit geometries, we will always assume that both the initial CCPSR manifold and the considered limit geometry are in standard form.
	
	\begin{Lem}\label{lem_convergence_limit_in_Cn}
		The limit geometry $\overline{\mathcal{H}}\subset\{\overline{h}=1\}$ in Definition \ref{def_limit_poly_limit_behaviour} is, in fact, a CCPSR manifold as claimed.
			\begin{proof}
				Assuming the existence of the limit polynomial $\overline{h}$, it follows by construction that it is in standard form. The generating set $\mathcal{C}_n$ \eqref{eqn_Cn_gen_set} is compact, $A(\gamma(t))$ is smooth, and
					\begin{equation*}
						\max\limits_{\|y\|=1} P_3(y)|_t
					\end{equation*}
				depends continuously, but not necessarily smoothly, cf. the proof of \cite[Thm.\,1.1]{L2}, on $t\in [0,R)$. Hence, $\overline{h}$ as a limit of $h|_t=x^2-x\langle y,y\rangle + P_3(y)|_t$ is contained in $\mathcal{C}_n$ and our claim follows with Theorem \ref{thm_Cn}.
			\end{proof}
	\end{Lem}
	
	Note that this process in particular also allows us to obtain a limit of the evolution of the centro-affine fundamental form in the following sense.
	
	\begin{Prop}\label{prop_metric_convergence}
		Let $(\overline{\mathcal{H}},g_{\overline{\mathcal{H}}})$ be a limit geometry of a CCPSR manifold $(\mathcal{H},g_\mathcal{H})$ with respect to a curve $\gamma:[0,R)\to\mathcal{H}$ as in \eqref{eqn_curve_for_limitgeo}. Then for every compactly embedded open subset $U\subset\overline{\mathcal{H}}$ and every $\varepsilon>0$ there exists a compactly embedded open subset $U'\subset\mathcal{H}$ and a diffeomorphism $F:\overline{U}\to \overline{U'}$, such that
			\begin{equation*}
				\|g_{\overline{\mathcal{H}}}-F^*g_\mathcal{H}\|_{g_{\overline{\mathcal{H}}}}<\varepsilon
			\end{equation*}
		in $\overline{U}$. If $U$ contains the point $\left(\begin{smallmatrix}1\\0\end{smallmatrix}\right)\in\overline{\mathcal{H}}$, there exists $N\in[0,R)$ such that for all $t\in[N,R)$, $U'$ can be chosen to contain the point $\gamma(t)\in\mathcal{H}$.
			\begin{proof}
				We use the notation as in \eqref{eqn_curve_for_limitgeo} so that $\mathcal{H}\subset\{h=h|_0=1\}$, $\overline{\mathcal{H}}\subset\left\{\lim\limits_{t\to R}h|_t=\overline{h}=1\right\}$ and define $\mathcal{H}_t$ to be the CCPSR manifold in standard form contained in the level set $\{h|_t=1\}$ for all $t\in[0,R)$. In particular $\mathcal{H}=\mathcal{H}_0$ and we might identify $\mathcal{H}_R:=\overline{\mathcal{H}}$. $\overline{\mathcal{H}}$ being a limit geometry of $\mathcal{H}$ means that with
					\begin{equation*}
						\overline{h}=x^3-x\langle y,y\rangle + \overline{P_3}(y),
					\end{equation*}
				$\overline{P_3}(y)$ is the limit $\lim\limits_{t\to R} P_3(y)|_t$ in $\mathrm{Sym}^3(\mathbb{R}^n)^*$. Recall the previously mentioned notation for any CCPSR manifold in standard form $\mathcal{H}$
					\begin{equation*}
						\mathrm{dom}(\mathcal{H})=\mathrm{pr}_{\mathbb{R}^n}\left((\mathbb{R}_{>0}\cdot\mathcal{H})\cap\left(\left(\begin{smallmatrix}1\\0\end{smallmatrix}\right)+T_{\left(\begin{smallmatrix}1\\0\end{smallmatrix}\right)}\mathcal{H}\right)\right)\subset\mathbb{R}^n,
					\end{equation*}
				where $\mathrm{pr}_{\mathbb{R}^n}$ denotes the canonical projection onto the last $n$ coordinates. Note that $\mathrm{dom}(\mathcal{H})$ depends not only on the equivalence class of $\mathcal{H}$, but also on the chosen linear coordinates $(x,y_1,\ldots,y_n)$ of the ambient space $\mathbb{R}^{n+1}$. The two sets $\mathrm{dom}(\mathcal{H})$ and $\mathcal{H}$ are diffeomorphic via the central parametrisation
					\begin{equation*}
						f:\mathrm{dom}(\mathcal{H})\to\mathcal{H},\quad y\mapsto h\left(\left(\begin{smallmatrix}1\\y\end{smallmatrix}\right)\right)^{-\frac{1}{3}}\left(\begin{smallmatrix}1\\y\end{smallmatrix}\right).
					\end{equation*}
				Furthermore, $\mathrm{dom}(\mathcal{H})$ is convex and precompact \cite[Sect.\,1.3]{CNS}, and one can show that for all CCPSR manifolds in standard form $\mathcal{H}$, $B^{\mathrm{Eucl.}}_{\sqrt{3}/2}(0)\subset U\subset B^{\mathrm{Eucl.}}_{\sqrt{3}}(0)$, where $B^{\mathrm{Eucl.}}_{r}(0)\subset\mathbb{R}^n$ denotes the open Euclidean ball of radius $r>0$ in the coordinates $(y_1,\ldots,y_n)$ \cite[Cor.\,4.4]{L2}. Instead of working with the CCPSR manifolds $(\mathcal{H}_t,g_{\mathcal{H}_t})$ in this proof directly, we will consider $(\mathrm{dom}(\mathcal{H}_t),g_t)$ whith $g_t:=f^*g_{\mathcal{H}_t}$ for all $t\in [0,R]$. First observe that for any compactly embedded open subset $U\subset\mathrm{dom}(\overline{\mathcal{H}})$ there exists $N\in[0,R)$ such that for all $t\in[N,R)$, $U$ is a compactly embedded open subset of $\mathrm{dom}(\mathcal{H}_t)$. This follows from the convergence of $P_3(y)|_t\to \overline{P}_3(y)$ as $t\to R$, where we note that for $q\in\partial B^{\mathrm{Eucl.}}_1(0)$ arbitrary but fixed, $\partial \mathrm{dom}(\mathcal{H}_t)\cap(\mathbb{R}_{>0}\cdot q)$ is precisely the smallest positive zero in $s$ of
					\begin{equation*}
						h|_t\left(\left(\begin{smallmatrix}1\\ sq\end{smallmatrix}\right)\right)=1-s^2+s^3P_3(q)|_t,
					\end{equation*}
				and we can thus use the compactness of $q\in\partial B^{\mathrm{Eucl.}}_1(0)$ to find $N$ as required. Now, for $t\in[N,R)$ in the coordinates $(y_1,\ldots,y_n)$, the coefficients of the Riemannian metric $g_t$ converge on $\overline{U}\subset\mathrm{dom}(\mathcal{H}_t)$ uniformly to the coefficients of $g_R=f^*g_{\overline{\mathcal{H}}}$ as $t\to R$. To see this one uses the convergence of $P_3(y)|_t\to \overline{P}_3(y)$ as $t\to R$, $U$ being compactly embedded in $\mathrm{dom}(\mathcal{H}_t)$ for all $t\in[N,R)$ and hence bounded away in Euclidean distance from the boundary of $\mathrm{dom}(\mathcal{H}_t)\cap\mathrm{dom}(\overline{\mathcal{H}})$, and the local formula for $g_t$, $t\in[N,R)$, on $\overline{U}$
					\begin{equation*}
						g_t= -\frac{\partial^2 h|_t\left(\left(\begin{smallmatrix}1\\y\end{smallmatrix}\right)\right)}{3 h|_t\left(\left(\begin{smallmatrix}1\\y\end{smallmatrix}\right)\right)} + \frac{2 \,\D\left(h|_t\left(\left(\begin{smallmatrix}1\\y\end{smallmatrix}\right)\right)\right)^2}{9 h|_t\left(\left(\begin{smallmatrix}1\\y\end{smallmatrix}\right)\right)^2},
					\end{equation*}
				cf. \cite[Lem.\,1.12]{CNS} or \cite[Lem.\,2.8]{L2} for our conventions, where the operator $\partial^2$ is to be understood in the coordinates $(y_1,\ldots,y_n)$. This completes the proof of the first claim of this proposition. The second claim, that is $\gamma(t)\in U'$, follows automatically since in the coordinates $(y_1,\ldots,y_n)$ we are working in, $U=U'$ and $y(\gamma(t))=0$ for all $t\in[0,R]$ by our definition of limit geometries of CCPSR manifolds.
			\end{proof}
	\end{Prop}
	
	Note that our definition of limit geometry is strictly speaking defined up to equivalence of the obtained CCPSR manifold $\overline{\mathcal{H}}$, as we still have the freedom to multiply $A(\gamma(t))$ with an element in $\mathrm{Id}_{x}\times\mathrm{O}(n)$, both from the left and from the right. Aside from the problems of actually calculating such limits, it is a priori not clear which CCPSR manifolds can even occur as limits. It will turn out that we need one more definition and an associated result for CCPSR manifolds.
		
	\begin{Def}\label{def_singular_at_infinity}
		A CCPSR manifold $\mathcal{H}\subset\{h=1\}$ is called \textbf{singular at infinity} if there exists a non-zero point in the boundary of the cone $\mathbb{R}_{>0}\cdot\mathcal{H}\subset\mathbb{R}^{n+1}$ spanned by $\mathcal{H}$, $p\in\partial(\mathbb{R}_{>0}\cdot\mathcal{H})\setminus\{0\}$, such that $\D h_p=0$.
	\end{Def}

	In order to, realistically, check whether or not a CCPSR manifold is singular at infinity we have the following tool at hand.
	
	\begin{Lem}\label{lem_sing_at_infty_iff_max_cond}
		A CCPSR manifold $\mathcal{H}\subset\{h=x^3-x\langle y,y\rangle + P_3(y)=1\}$ in standard form is singular at infinity if and only if $\max\limits_{\|y\|=1} P_3(y)= \frac{2}{3\sqrt{3}}$.
			\begin{proof}
				\cite[Lem.\,4.5\,]{L2}.
			\end{proof}
	\end{Lem}

	For CCPSR manifolds $\mathcal{H}$ the concept of singular at infinity is closely related to the behaviour of the Hessian of the defining polynomial along the boundary of their cone $\mathbb{R}_{>0}\cdot\mathcal{H}$.
	
	\begin{Def}
		A CCPSR manifold $\mathcal{H}\subset\{h=1\}$ has \textbf{regular boundary behaviour} if it is not singular at infinity and
			\begin{equation*}
				-\partial^2h|_{T (\partial(\mathbb{R}_{>0}\cdot\mathcal{H})\setminus\{0\})\times T (\partial(\mathbb{R}_{>0}\cdot\mathcal{H})\setminus\{0\})}\geq 0,\quad \dim\left(-\partial^2h|_{T (\partial(\mathbb{R}_{>0}\cdot\mathcal{H})\setminus\{0\})\times T (\partial(\mathbb{R}_{>0}\cdot\mathcal{H})\setminus\{0\})}\right)=1,
			\end{equation*}
		for all $p\in\partial(\mathbb{R}_{>0}\cdot\mathcal{H})\setminus\{0\}$.
	\end{Def}

	One can show that for CCPSR manifolds $\mathcal{H}\subset\{h=1\}$, the behaviour of the then possibly positive semi-definite bilinear form \eqref{eqn_symbilform_to_diag} along the set $\partial(\mathbb{R}_{>0}\cdot\mathcal{H})\setminus\{0\}$ can only have more than one zero eigenvalue if $h$ is already singular at infinity at that specific point.
	
	\begin{Th}
		A CCPSR manifold has regular boundary behaviour if and only if it is not singular at infinity.
			\begin{proof}
				\cite[Thm.\,4.12]{L2}.
			\end{proof}
	\end{Th}

	Before stating our main result which answers all of the above questions for every dimension we will study an explicit example in dimension $2$.
	
	\begin{Ex}\label{ex_motivating}
		Let $h:\mathbb{R}^3\to\mathbb{R}$, $h=x^3-x(y^2+z^2) + \frac{2}{3\sqrt{3}}y^3$, and let $\mathcal{H}\subset\{h=1\}$ be the corresponding CCPSR surface in standard form. $\mathcal{H}$ is equivalent to \cite[Thm.\,1\,e)]{CDL}, cf. \cite[Ex.\,3.2]{L2}. With Lemma \ref{lem_sing_at_infty_iff_max_cond} it follows that $\mathcal{H}$ is singular at infinity. For the constant vector field on the ambient space $\mathbb{R}^3$ we consider $V=\partial_y$. Up to reparametrisation, the maximal integral curve of the corresponding central projection of $V$ to $\mathcal{H}$ is given by
			\begin{equation*}
				\gamma:\left(-\frac{\sqrt{3}}{2},\sqrt{3}\right)\to\mathcal{H},\quad \gamma(t)=\left(\begin{matrix} \beta^{-\frac{1}{3}}\\ t\beta^{-\frac{1}{3}}\\0\end{matrix}\right),
			\end{equation*}
		where $\beta=\beta(t):=h\left(\left(\begin{smallmatrix}1\\ t\\ 0\end{smallmatrix}\right)\right)=1-t^2 + \frac{2}{3\sqrt{3}}t^3$. For the transformation $A$ as in \eqref{eqn_pmoving_trafo_explicit} along $\gamma$ we choose
			\begin{equation*}
				A(\gamma(t))=\left(\begin{array}{c|c|c}
					\beta^{-\frac{1}{3}} & \frac{2t}{3}\underset{}{\beta^{-\frac{1}{3}}} & 0\\ \hline
					\overset{}{t} \overset{}{\beta^{-\frac{1}{3}}} & \left(1+\frac{t}{\sqrt{3}}\right) \underset{}{\beta^{-\frac{1}{3}}} & 0\\ \hline
					0 & 0 & \overset{}{\beta^{\frac{1}{6}}}
				\end{array}\right)
			\end{equation*}
		and obtain
			\begin{equation*}
				h|_t=A(\gamma(t))^* h = x^3 - x(y^2 + z^2) -\frac{2t}{3}yz^2 + \frac{2}{3\sqrt{3}} y^3
			\end{equation*}
		for all $t\in \left(-\frac{\sqrt{3}}{2},\sqrt{3}\right)$. In this specific case, the only difficulty lies in determining the above choice for $A(\gamma(t))$ as the two possible limit polynomials are very easy to calculate. We obtain
			\begin{equation*}
				\overline{h}_+:=\lim\limits_{t\to \sqrt{3}} h|_t = x^3 - x(y^2 + z^2) +\frac{2}{3\sqrt{3}} y^3 -\frac{2}{\sqrt{3}}yz^2,
			\end{equation*}
		which is equivalent to \cite[Thm.\,1\,a)]{CDL}, cf. \cite[Ex.\,3.2]{L2}, and
			\begin{equation*}
				\overline{h}_-:=\lim\limits_{t\to -\frac{\sqrt{3}}{2}} h|_t = x^3 - x(y^2 + z^2) +\frac{2}{3\sqrt{3}} y^3+\frac{1}{\sqrt{3}}yz^2,
			\end{equation*}
		which is equivalent to \cite[Thm.\,1\,b)]{CDL}, cf. \cite[Ex.\,3.2]{L2}. In particular, the corresponding limit CCPSR curves are not equivalent, neither to each other nor to the initial CCPSR surface $\mathcal{H}$, as the maximal hyperbolic component of $\{\overline{h}_+=1\}$ containing $(x,y,z)^\mathrm{T}=(1,0,0)^\mathrm{T}$ is isometric to the flat $\mathbb{R}^2$, the maximal hyperbolic component of $\{\overline{h}_-=1\}$ containing $(x,y,z)=(1,0,0)$ is isometric to the hyperbolic plane, and $\mathcal{H}$ is not a homogeneous space. In order to better understand what is happening geometrically for different values of $t$ we plot $\partial\mathrm{dom}(\mathcal{H}_t)$, where $\mathcal{H}_t\subset\{h|_t=1\}$, for $t\in\left\{-\frac{\sqrt{3}}{2},0,\frac{\sqrt{3}}{2},\sqrt{3}\right\}$ and obtain Figure \ref{fig_mot_example}.
			\begin{figure}[H]%
			\centering%
			\includegraphics[scale=0.19]{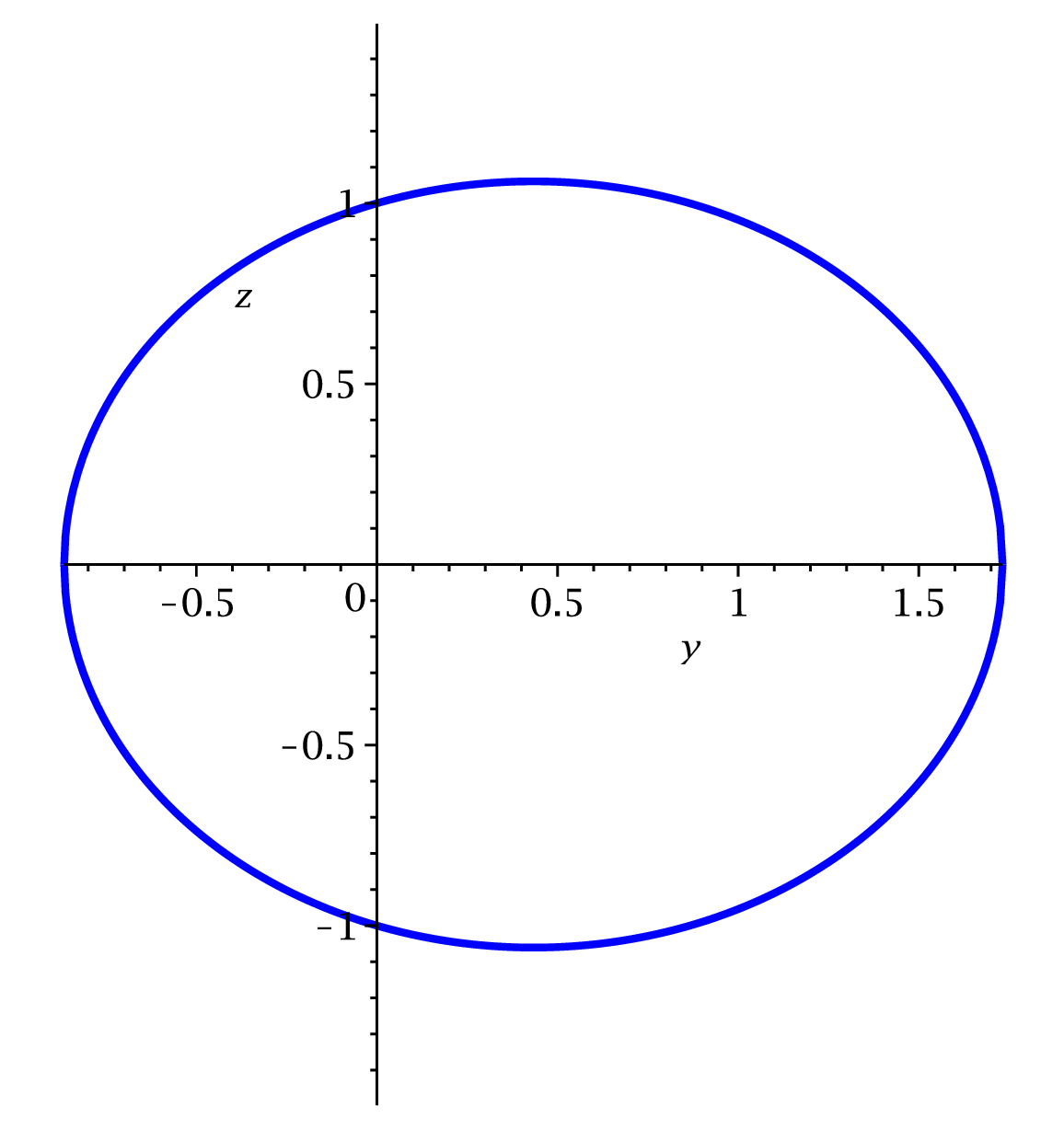}%
			\includegraphics[scale=0.19]{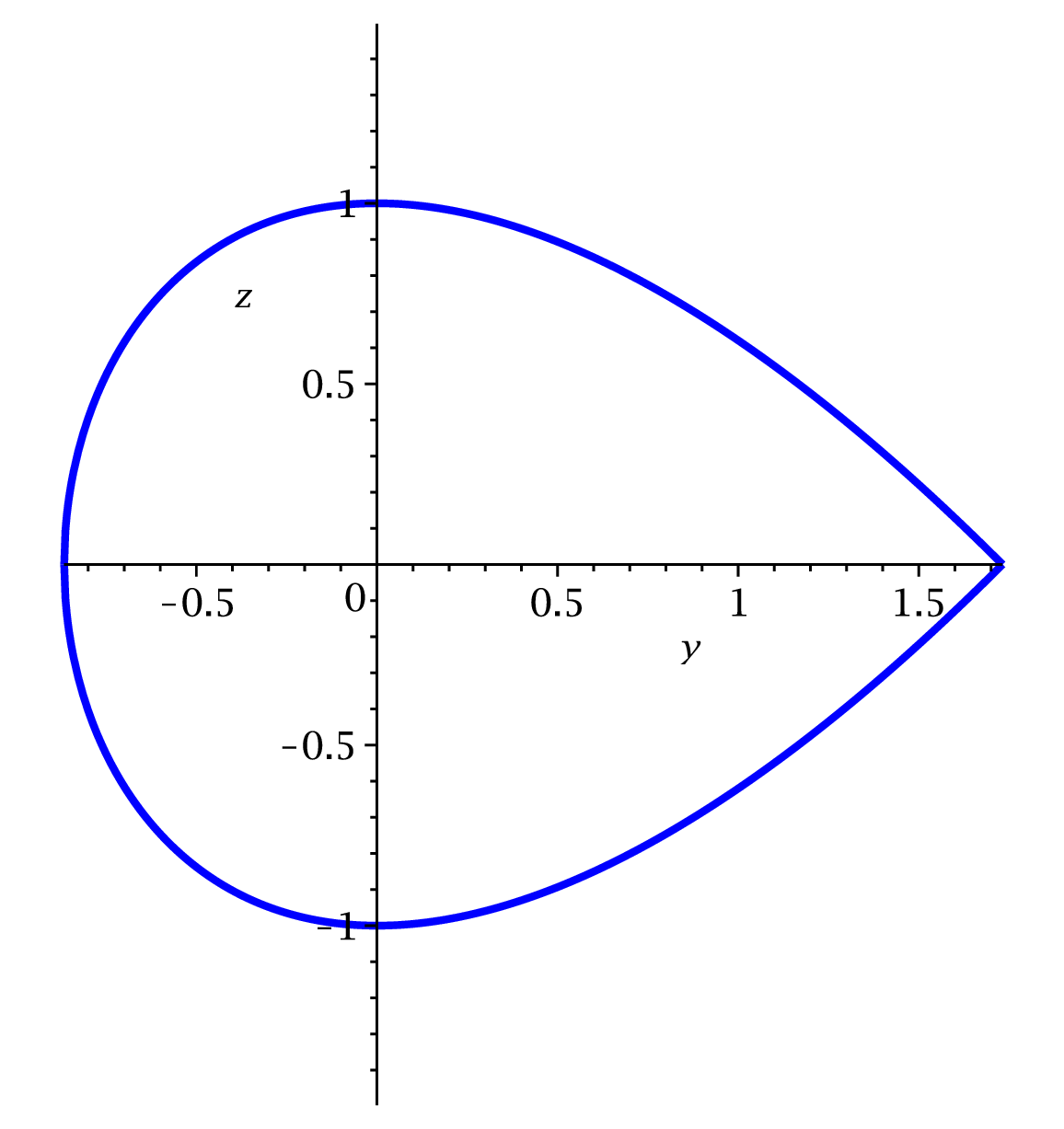}%
			\includegraphics[scale=0.19]{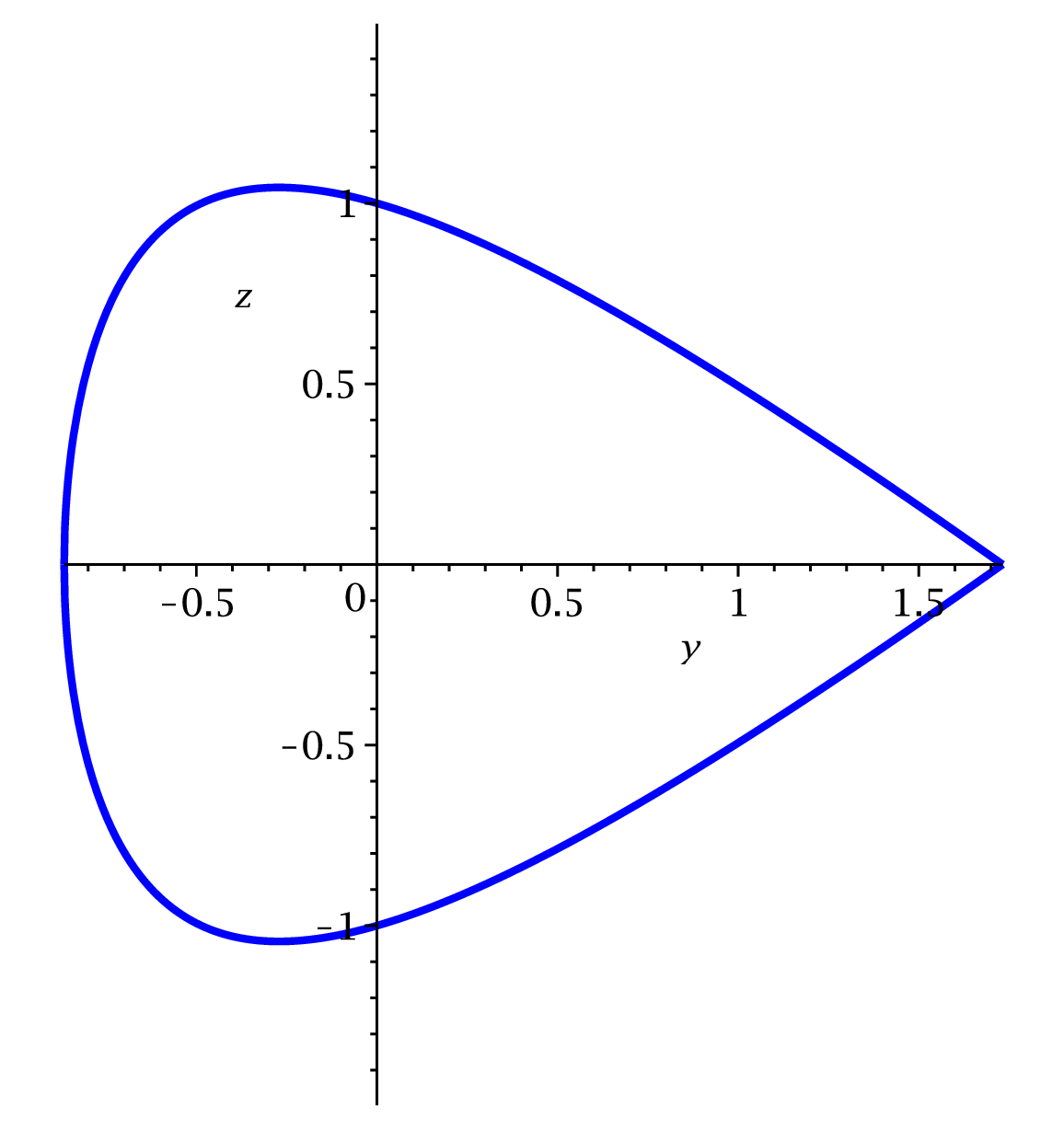}%
			\includegraphics[scale=0.19]{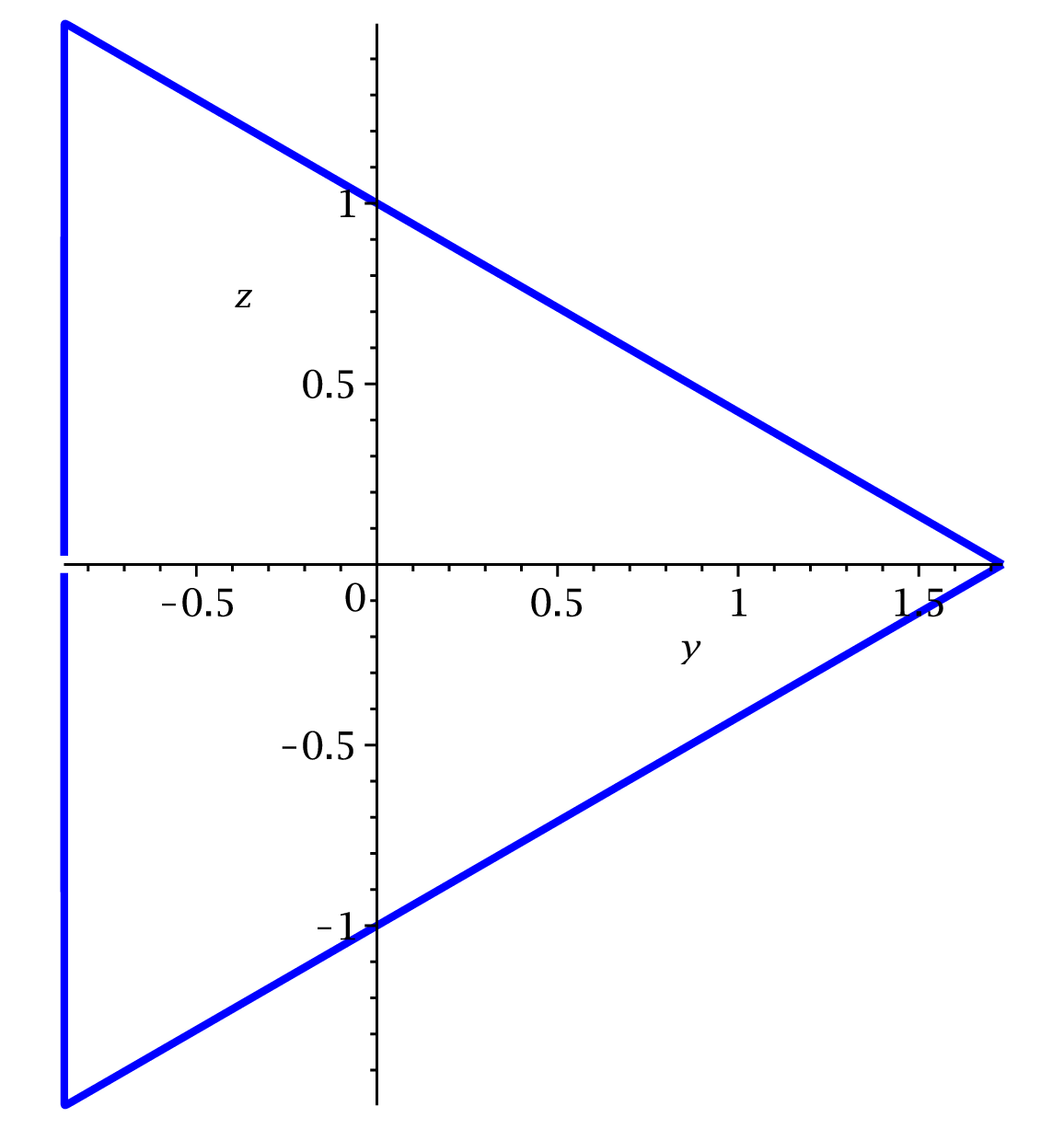}%
			\caption{Plots of $\partial\mathrm{dom}(\mathcal{H}_t)$ for, from left to right, $t=-\frac{\sqrt{3}}{2}$, $t=0$, $t=\frac{\sqrt{3}}{2}$, and $t=\sqrt{3}$.}\label{fig_mot_example}%
			\end{figure}
		\noindent
		Observe that for all $t\in\left[-\frac{\sqrt{3}}{2},\sqrt{3}\right]$, including in particular the limits, the CCPSR surfaces in standard form corresponding to $h|_t$ are singular at infinity. For a related discussion see also \cite[Ex.\,4.18]{L2}.
	\end{Ex}
	
	We now have all tools at hand to state the main result of this work.
	
	\begin{Th}\label{thm_limit_behaviour_ccpsr}
		Let $\mathcal{H}\subset\{h=1\}$ be a CCPSR manifold in standard form. Then every possible limit polynomial $\overline{h}=\lim\limits_{t\to R}h|_t$ is well defined and the corresponding limit geometry $\overline{\mathcal{H}}$ is a CCPSR manifold. The limit polynomial $\overline{h}$ is equivalent to one of the following:
			\begin{enumerate}[(i)]
				\item\label{thm_limit_behaviour_ccpsr_i} \underline{$\dim(\mathcal{H})=1$:}
					\begin{equation*}
						\quad\overline{h}=x^3-xy^2 -\frac{2}{3\sqrt{3}}y^3,
					\end{equation*}
				$\overline{\mathcal{H}}$ is a homogeneous space and equivalent to \cite[Thm.\,8\,a)]{CHM}.
				\item\label{thm_limit_behaviour_ccpsr_ii} \underline{$\dim(\mathcal{H})=2$:}
					\begin{enumerate}[(a)]
						\item\label{thm_limit_behaviour_ccpsr_ii_a}
							\begin{equation*}
								\overline{h}=x^3 - x(y^2 + z^2) -\frac{1}{\sqrt{3}}yz^2 - \frac{2}{3\sqrt{3}} y^3,
							\end{equation*}
							$\overline{\mathcal{H}}$ is a homogeneous space isomorphic to the hyperbolic plane.
						\item\label{thm_limit_behaviour_ccpsr_ii_b}
							\begin{equation*}
								\overline{h}=x^3 - x(y^2 + z^2) +\frac{2}{\sqrt{3}}yz^2 - \frac{2}{3\sqrt{3}} y^3,
							\end{equation*}
						$\overline{\mathcal{H}}$ is a homogeneous space isomorphic to the flat $\mathbb{R}^2$.
					\end{enumerate}
				\item\label{thm_limit_behaviour_ccpsr_iii} \underline{$\dim(\mathcal{H})\geq 3$:}
						\begin{align}
							\overline{h}&=x^3-x(\langle s,s\rangle + \langle u,u\rangle + w^2)\notag\\
							&\quad +\sum\limits_{i=1}^{m}  s_i\langle u,F_i u\rangle + \left(\frac{2}{\sqrt{3}}\langle s,s\rangle - \frac{1}{\sqrt{3}}\langle u,u\rangle\right)w - \frac{2}{3\sqrt{3}} w^3,\label{eqn_dim_geq_3_limits}
						\end{align}
					where $s=(s_1,\ldots,s_m)^\mathrm{T}$, $u=(u_1,\ldots,u_{n-1-m})^\mathrm{T}$ for some $0\leq m\leq n-1$, and each $F_i$, $1\leq i\leq m$, is a symmetric $(n-1-m)\times (n-1-m)$-matrix such that the eigenvalues of all matrices of the form
						\begin{equation*}
							\sum\limits_{i=1}^{m}c_i F_i
						\end{equation*}
					for all $c=(c_1,\ldots,c_{m})^\mathrm{T}\in\mathbb{R}^{m}$ with $\|c\|=1$ are contained in $[-1,1]$. The symmetry group of the corresponding CCPSR manifold $\overline{\mathcal{H}}$ is of dimension at least $1$.
			\end{enumerate}
		All of the above CCPSR manifolds in standard form in \eqref{thm_limit_behaviour_ccpsr_i}--\eqref{thm_limit_behaviour_ccpsr_iii} are singular at infinity and can be realised as a limit geometry of a CCPSR manifold, in particular including themselves in each case.
	\end{Th}

	In the proof of Theorem \ref{thm_limit_behaviour_ccpsr} we will make use of the following convenient tool which has been introduces in \cite[Def.\,3.7]{L2} in a more general setting. It measures the infinitesimal change of the term $P_3(y)$ in the standard form of the defining polynomial of a CCPSR manifold \eqref{eqn_h_general_form} when changing the reference point for said standard form.
	
	\begin{Def}\label{def_first_var_P3}
		Let $\mathcal{H}\subset\{h=x^3-x\langle y,y\rangle + P_3(y)=1\}$ be a maximal connected PSR manifold in standard form. Let $U\subset \mathcal{H}$, $\left(\begin{smallmatrix}1\\0\end{smallmatrix}\right)$, be an open set and $A:U\to\mathrm{GL}(n+1)$ as in \eqref{eqn_pmoving_trafo_explicit}. Denote
			\begin{equation*}
				A(p)^*h=x^3-x\langle y,y\rangle + P_3(y)(p),
			\end{equation*}
		that is view $P_3(y)(\cdot):U\to\mathrm{Sym}^3(\mathbb{R}^n)^*$ with $P_3(y)\left(\left(\begin{smallmatrix}1\\0\end{smallmatrix}\right)\right)=P_3(y)$. Interpreting $A(\cdot)^*h:U\to\mathrm{Sym}^{3}(\mathbb{R}^{n+1})^*$ as a smooth map, we call the term
			\begin{equation*}
				\delta P_3(y):= \left.\D\left(A(\cdot)^*h\right)\right|_{p=\left(\begin{smallmatrix}1\\0\end{smallmatrix}\right)}
			\end{equation*}
		the \textbf{first variation of $P_3$}. We view $\delta P_3(y)$ as a linear map $\delta P_3(y):\mathbb{R}^n\to\mathrm{Sym}^3(\mathbb{R}^n)^*$.
	\end{Def}
	
	One can show that $\delta P_3(y)$ is of the following form.
	
	\begin{Lem}\label{lem_delta_P3_form}
		There exists a linear map $L:\mathbb{R}^n\to\mathfrak{so}(n)$ of the form $L=\sum\limits_{i=1}^n L_i\otimes \D y_i$, $L_i\in\mathfrak{so}(n)$ for all $1\leq i\leq n$, such that with $Ly=\sum\limits_{i=1}^n L_iy\otimes \D y_i$, $\delta P_3(y)$ is of the form
			\begin{equation*}
				\delta P_3(y)=-\frac{2}{3}\langle y,y\rangle \langle y,\D y\rangle + dP_3|_y\left(Ly+\frac{1}{4}\partial^2 P_3|_y \D y\right).
			\end{equation*}
		\begin{proof}
			\cite[Prop.\,3.6,\;Def.\,3.7]{L2}.
		\end{proof}
	\end{Lem}
	
	An immediate consequence of the construction of the first variation of $P_3$ is the following relation to symmetry groups of the corresponding PSR manifold.
	
	\begin{Lem}\label{lem_first_var_P3_sym_grps}
		Let $\mathcal{H}\subset\{h=x^3-x\langle y,y\rangle + P_3(y)=1\}$ be a CCPSR manifold in standard form. Suppose there exists $L:\mathbb{R}^n\to\mathfrak{so}(n)$ as in Lemma \ref{lem_delta_P3_form}, such that $\dim\ker\delta P_3(y)= k$, $1\leq k\leq n$, where we view $\delta P_3(y):\mathbb{R}^n\to\mathrm{Sym}^3(\mathbb{R}^n)^*$ as a linear map. Then the symmetry group of $\mathcal{H}$ is of dimension at least $k$.
			\begin{proof}
				This follows by writing out the differential of $h\left(A(p)\left(\begin{smallmatrix}x\\y\end{smallmatrix}\right)\right)$, which yields $k$ linearly independent linear maps $B_i\in\mathrm{GL}(n+1)$, $1\leq i\leq k$, such that $\D h_{\left(\begin{smallmatrix}x\\y\end{smallmatrix}\right)}\left(B_i\left(\begin{smallmatrix}x\\y\end{smallmatrix}\right)\right)\equiv 0$.
			\end{proof}
	\end{Lem}

	Lastly, we will prove the following and, at least to the author, surprising result about possible limit geometries of CCPSR manifolds that have regular boundary behaviour. 
	
	\begin{Th}\label{thm_generic_limit_behaviour}
		Let $\mathcal{H}$ be an $n\geq 2$-dimensional CCPSR manifold with regular boundary behaviour and let $(x,s_1,\ldots,s_{n-1},w)^\mathrm{T}=\left(\begin{smallmatrix}x\\s\\w\end{smallmatrix}\right)$ denote linear coordinates on $\mathbb{R}^{n+1}$. Then every limit geometry of $\mathcal{H}$ is, up to equivalence, given by the CCPSR manifold in standard form $\overline{\mathcal{H}}\subset\{\overline{h}=1\}$ corresponding to the cubic polynomial
			\begin{equation}
				\overline{h}=x^3-x(\langle s,s\rangle + w^2) -\frac{1}{\sqrt{3}}\langle s,s\rangle w -\frac{2}{3\sqrt{3}} w^3.\label{eqn_h_lie_grp_generic_limit}
			\end{equation}
		The CCPSR manifold $\overline{\mathcal{H}}$ is isomorphic to the Lie group $\mathbb{R}_{>0}\ltimes\mathbb{R}^{n-1}$. The action of $\mathbb{R}_{>0}\ltimes\mathbb{R}^{n-1}$ on the coordinates $\left(\begin{smallmatrix}x\\s\\w\end{smallmatrix}\right)$ is given by
			\begin{align*}
				&(\lambda,v)\cdot \left(\begin{matrix}x\\s\\w\end{matrix}\right)\\
				&=\left(\begin{array}{l}
						\lambda \left(2\cdot 3^{-3/2}w+2\cdot 3^{-1}x\right)\\ \underset{}{+} \lambda^2\left(-2\cdot3^{-3/2}w+3^{-1}x+\langle v,v\rangle\left(2^{1/2}\cdot 3^{-3/2}w + 2^{1/2}\cdot3^{-1}x\right)-2^{5/4}\cdot3^{-1}\langle s,v\rangle\right)\\\hline
						\overset{}{ \lambda^{-1/2}\left(s+\left(-2^{1/4}\cdot3^{-1/2}w-2^{1/4}x\right)v\right)}\underset{}{\vphantom{+} }\\\hline
						\overset{}{\lambda\left(3^{-1}w+3^{-1/2}x\right)}\\ + \lambda^{-2}\left(2\cdot 3^{-1}w-3^{-1/2}x+\langle v,v\rangle\left(-2^{1/2}\cdot3^{-1}w-2^{1/2}\cdot3^{-1/2}x\right)+2^{5/4}\cdot3^{-1/2}\langle s,v\rangle\right)\\
					\end{array}\right)
			\end{align*}
		for all $(\lambda,v)\in\mathbb{R}_{>0}\ltimes\mathbb{R}^{n-1}$, and the group multiplication is given by $(\lambda_1,v_1)\cdot (\lambda_2,v_2)=(\lambda_1\lambda_2,v_1+\lambda_1^{3/2}v_2)$ for all $(\lambda_1,v_1),(\lambda_2,v_2)\in\mathbb{R}_{>0}\ltimes\mathbb{R}^{n-1}$. For $n=2$, that is CCPSR surfaces, we have $\overline{\mathcal{H}}\cong\mathbb{R}_{>0}\ltimes\mathbb{R}\cong H$, where $H$ is the hyperbolic plane.
	\end{Th}
	
	Note that for dimension $\dim(\mathcal{H})=1$ a similar result holds, as there is in fact only one possible limit geometry for CCPSR curves, namely Thm. \ref{thm_limit_behaviour_ccpsr} \eqref{thm_limit_behaviour_ccpsr_i}.
	
	A consequence of Theorem \ref{thm_generic_limit_behaviour} about the moduli space of CCPSR manifolds when equipped with the quotient topology is as follows.
	
	\begin{Cor}\label{cor_non_hd}
		For $n\geq 2$ let $\mathrm{Sym}^3(\mathbb{R}^{n+1})^*_{\mathrm{CCPSR}}\subset \mathrm{Sym}^3(\mathbb{R}^{n+1})^*$ denote the set of cubic homogeneous polynomials in $n+1$ real variables containing a CCPSR manifold as a level set. Let further $\mathrm{Sym}^3(\mathbb{R}^{n+1})^*_{\mathrm{CCPSR,\ reg.}}\subset \mathrm{Sym}^3(\mathbb{R}^{n+1})^*_{\mathrm{CCPSR}}$ denote the subset thereof consisting of cubic homogeneous polynomials in $n+1$ real variables containing a CCPSR manifold as a level with regular boundary behaviour, and let $\overline{h}\in\mathrm{Sym}^3(\mathbb{R}^{n+1})^*_{\mathrm{CCPSR}}$ be as in equation \eqref{eqn_h_lie_grp_generic_limit}. In the quotient topology of $\mathrm{Sym}^3(\mathbb{R}^{n+1})^*_{\mathrm{CCPSR}}/\mathrm{GL}(n+1)$, the point $[\overline{h}]$ cannot be separated from any point in $\mathrm{Sym}^3(\mathbb{R}^{n+1})^*_{\mathrm{CCPSR,\ reg.}}/\mathrm{GL}(n+1)$ with disjoint open sets.
	\end{Cor}
	
	This in particular shows that $\mathrm{Sym}^3(\mathbb{R}^{n+1})^*_{\mathrm{CCPSR}}/\mathrm{GL}(n+1)$ equipped with the quotient topology is not Hausdorff for $n\geq 2$. This in fact also holds for $n=1$, that is CCPSR curves, which directly follows from Theorem \ref{thm_limit_behaviour_ccpsr}.
	
\section{Proof of Theorem \ref{thm_limit_behaviour_ccpsr}}
	In this section we will prove Theorem \ref{thm_limit_behaviour_ccpsr}. We will proceed by dimensions of the CCPSR manifold $\mathcal{H}$ in consideration. We start with $\dim(\mathcal{H})=1$ in Section \ref{sect_main_thm_pf_dim1} which is the easiest and most straightforward case. In Sections \ref{subsect_dim2_b_ne_0} and \ref{sect_b_0} we will study the limit geometry of CCPSR surfaces in detail. These cases are already a lot more involved than the $1$-dimensional cases. The employed techniques and results will also be used in Sections \ref{section_dimH_geq3_bnon0} and \ref{subsection_dimHgeq3_b0} where we deal with CCPSR manifolds $\mathcal{H}$ of dimension $\dim(\mathcal{H})\geq 3$. The main difficulty in this part of the proof lies in the fact that, in comparison with the $2$-dimensional case, the obtained limit geometry CCPSR manifolds are not necessarily homogeneous, and it is at the beginning unclear what one should expect to obtain as the most general form of limit geometries instead. After having solved this problem we will study the details of limit geometries in the $\dim(\mathcal{H})\geq 3$-cases, which amounts to studying the $F_i$-terms in \eqref{eqn_dim_geq_3_limits} and, depending on the integer $m$ in \eqref{eqn_dim_geq_3_limits}, prove the claimed lower bound of the dimension of the symmetry group of a given limit geometry CCPSR manifold. For the latter we will use the concept of first variation of the $P_3$-term in \eqref{eqn_h_general_form} as introduced in \cite[Def.\,3.7]{L2}, cf. Definition \ref{def_first_var_P3}.

\subsection{$\dim(\mathcal{H})=1$}\label{sect_main_thm_pf_dim1}
	For $\dim(\mathcal{H})=1$ we can write down the set $\mathcal{C}_1$ \eqref{eqn_Cn_gen_set} explicitly. Since $P_3(y)=a y^3$, $a\in\mathbb{R}$, for all defining polynomials $h=x^3-xy^2+P_3(y)$ of a connected PSR curve in standard form we have
		\begin{equation*}
			\max\limits_{y^2=1} P_3(y)=|a|
		\end{equation*}
	and, hence,
		\begin{equation*}
			\mathcal{C}_1=\left\{x^3-xy^2 + ay^3\ \left|\ |a|\leq\frac{2}{3\sqrt{3}}\right.\right\}.
		\end{equation*}
	Let
		\begin{equation*}
			\beta:=\beta(r)=h\left(\left(\begin{smallmatrix}1\\r\end{smallmatrix}\right)\right).
		\end{equation*}
	We will use a similar notation for the higher-dimensional cases of the proof as well. We will utilize the following property of $\beta$, and this property will also be utilized for all higher dimensional cases.
	
	\begin{Lem}\label{lem_beta_zeros}
		Let $\mathcal{H}\subset\{h=1\}$ be a CCPSR manifold in standard form and define
			\begin{equation}
				\beta:=h\left((1,0,\ldots,0,r)^\mathrm{T}\right).\label{eqn_beta_def_general}
			\end{equation}
		Then $\beta=1-r^2+ar^3$ has precisely three simple zeros if $|a|<\frac{2}{3\sqrt{3}}$, and one simple zero and one double zero if $|a|=\frac{2}{3\sqrt{3}}$. Furthermore, $|a|=\frac{2}{3\sqrt{3}}$ implies that $\mathcal{H}$ is singular at infinity. If $R$ denotes the smallest positive zero of $\beta$, $a=\frac{R^2-1}{R^3}$ and, hence, $R=\sqrt{3}$ or $R=\frac{\sqrt{3}}{2}$ implies that $\mathcal{H}$ is singular at infinity.
	\end{Lem}
	
	The proof of the above lemma is not difficult and left to the reader. Since by assumption $|a|\leq \frac{2}{3\sqrt{3}}$, $\beta$ has either $3$ distinct simple zeros if $|a|<\frac{2}{3\sqrt{3}}$, or $1$ simple zero and $1$ double zero if $|a|=\frac{2}{3\sqrt{3}}$. In all of the cases the zeros of $\beta$ are real. If $R>0$ is the smallest positive zero of $\beta$, $a$ is given by
		\begin{equation*}
			a=\frac{R^2-1}{R^3}.
		\end{equation*}
	Note that this means that $\beta$ has a simple zero in $R$ if $R\in\left[\frac{\sqrt{3}}{2},\sqrt{3}\right)$ and a double zero if $R=\sqrt{3}$. The only curve as in the introduction we need to consider in for $\dim(\mathcal{H})=1$ is
		\begin{equation*}
			\gamma:[0,R)\to\mathcal{H},\quad r\mapsto \beta^{-\frac{1}{3}}\left(\begin{smallmatrix}1\\r\end{smallmatrix}\right).
		\end{equation*}
	We further split up the $\dim(\mathcal{H})$-case in two cases.
	
	First, assume that the smallest positive zero $R>0$ of $\beta$ fulfils $R<\sqrt{3}$. This, together with our assumption that $\mathcal{H}$ is closed and thus $|a|\leq\frac{2}{3\sqrt{3}}$ by Theorem \ref{thm_Cn}, is equivalent to $a\in\left[-\frac{2}{3\sqrt{3}},\frac{2}{3\sqrt{3}}\right)$. We need to determine $A(\gamma(r))$ as in \eqref{eqn_pmoving_trafo_explicit}. In this case, the necessary calculations leave no freedom of choice except for a sign choice since $\mathrm{O}(1)=\mathbb{Z}_2$, and we obtain
		\begin{equation*}
			A(\gamma(r))=
				\left(
					\begin{array}{c|c}
						\underset{}{\beta^{-\frac{1}{3}}} & \frac{r(2-3ar)}{\sqrt{3}\sqrt{3-9ar+r^2}}\beta^{-\frac{1}{3}}\\ \hline
						r\overset{}{\beta^{-\frac{1}{3}}} & \frac{3-r^2}{\sqrt{3}\sqrt{3-9ar+r^2}}\beta^{-\frac{1}{3}}
					\end{array}
				\right).
		\end{equation*}
	We further calculate
		\begin{equation}
			A(\gamma(r))^* h = x^3 - xy^2 + \frac{27a-18r+27ar^2+(2-27a^2)r^3}{3\sqrt{3}\sqrt{3-9ar+r^2}^3}y^3.\label{eqn_pullback_1dim}
		\end{equation}
	Note that for $a=-\frac{2}{3\sqrt{3}}$, the $y^3$-coefficient in \eqref{eqn_pullback_1dim} is constant and equal to $-\frac{2}{3\sqrt{3}}$. This is the expected behaviour, as for $a=-\frac{2}{3\sqrt{3}}$ the corresponding CCPSR curve in standard form $\mathcal{H}$ is homogeneous as noted in Theorem \ref{thm_limit_behaviour_ccpsr}, which in this case in particular holds. For $a\in \left(-\frac{2}{3\sqrt{3}},\frac{2}{3\sqrt{3}}\right)$ the denominator in the $y^3$-coefficient in \eqref{eqn_pullback_1dim} has no zeros in $[0,R]$, so calculating the limit reduces to inserting $a=\frac{R^2-1}{R^3}$ and $r=R$ to obtain
		\begin{equation*}
			\lim\limits_{r\nearrow R} \frac{27a-18r+27ar^2+(2-27a^2)r^3}{3\sqrt{3}\sqrt{3-9ar+r^2}^3} = -\frac{2}{3\sqrt{3}}.
		\end{equation*}
	Hence, Theorem \ref{thm_limit_behaviour_ccpsr} holds in these cases.
	
	Lastly consider the case $R=\sqrt{3}$ or, equivalently, $a=\frac{2}{3\sqrt{3}}$. In this case, the the denominator of the $y^3$-coefficient in \eqref{eqn_pullback_1dim} has a zero in $r=R$, it however turns out that
		\begin{equation*}
			\left.\frac{27a-18r+27ar^2+(2-27a^2)r^3}{3\sqrt{3}\sqrt{3-9ar+r^2}^3}\right|_{a=\frac{2}{3\sqrt{3}}}\equiv \frac{2}{3\sqrt{3}},
		\end{equation*}
	thus holding in particular for the limit $r\to\sqrt{3}$. Note that this is, again, expected behaviour as the corresponding CCPSR curve in standard form $\mathcal{H}$ is homogeneous and equivalent to the aforementioned homogeneous case by a sign flip in $y$. Hence, Theorem \ref{thm_limit_behaviour_ccpsr} holds in this case as well and we have proven Thm. \ref{thm_limit_behaviour_ccpsr} \eqref{thm_limit_behaviour_ccpsr_i}.

\subsection{$\dim(\mathcal{H})=2$, $b\ne 0$}\label{subsect_dim2_b_ne_0}
Next, we will show that Theorem \ref{thm_limit_behaviour_ccpsr} holds in dimension $2$, that is for CCPSR surfaces. We start with $P_3\left(\left(\begin{smallmatrix}v\\ w\end{smallmatrix}\right)\right):=cv^3+qv^2w+bvw^2+aw^3$, $a,b,c,q\in\mathbb{R}$, so that
	\begin{equation}\label{eqn_starting_P3_dim2}
		\max\limits_{\{v^2+w^2=1\}}P_3\left(\left(\begin{smallmatrix}v\\ w\end{smallmatrix}\right)\right)\leq \frac{2}{3\sqrt{3}},
	\end{equation}
and set
	\begin{equation*}
		h\left(\left(\begin{smallmatrix}x\\ v\\ w\end{smallmatrix}\right)\right):=x^3-x(v^2+w^2)+P_3\left(\left(\begin{smallmatrix}v\\ w\end{smallmatrix}\right)\right).
	\end{equation*}
Now define $\mathcal{H}\subset\{h=1\}$ to be the connected component of $\{h=1\}$ that contains the point $(x,v,w)^\mathrm{T}=(1,0,0)^\mathrm{T}$. Recall that the condition \eqref{eqn_starting_P3_dim2} is equivalent to $\mathcal{H}\subset\mathbb{R}^3$ being closed in the subspace topology. We will study the standard form of $h$ along
	\begin{equation*}
		r\mapsto\left(\begin{smallmatrix}x\\ v\\ w\end{smallmatrix}\right)=\beta^{-\frac{1}{3}}\cdot\left(\begin{smallmatrix} 1\\ 0\\ r\end{smallmatrix}\right)=:p(r)=\left(\begin{smallmatrix}p_x(r)\\ p_v(r)\\ p_w(r)\end{smallmatrix}\right)\in\mathcal{H},\quad \beta:=	h\left(\left(\begin{smallmatrix}1\\ 0\\ r\end{smallmatrix}\right)\right)=1-r^2+ar^3.
	\end{equation*}
In the above equation, the domain for $r$ is given by $[0,R)$, where
	\begin{equation}\label{eqn_dim2_R_def}
		R:=\sup\limits_{\left\{r>0,\ p(r)\in\mathcal{H}\right\}}r,
	\end{equation}
that is $r=R$ is the smallest positive real solution of $\beta=0$.
Equivalently, $R>0$ is the smallest, and by $\mathcal{H}$ being connected also the unique, positive real number, such that $(1,0,R)^\mathrm{T}\in\partial(\mathbb{R}_{>0}\cdot\mathcal{H})$. Following Proposition \ref{prop_standard_form}, we want to determine a smooth map $[0,R)\ni r\mapsto E(r)\in\mathrm{GL}(2)$, such that
	\begin{equation}\label{eqn_dim2_hE_def}
		h\left(\left(\begin{array}{c|cc}
			p_x(r)	&	\left.-\frac{\partial_v h}{\partial_x h}\right|_{p(r)}	&	\left.-\frac{\partial_w h}{\partial_x h}\right|_{p(r)}\\ \hline
			p_v(r)	&	1	&	0\\
			p_w(r)	&	0	&	1
		\end{array}\right)
		\cdot
		\left(\begin{array}{c|c}
			1	&	\\ \hline
				&	E(r)
		\end{array}\right)
		\cdot\left(\begin{matrix}x\\ v\\ w\end{matrix}\right)\right)=:h_E\left(\left(\begin{smallmatrix} x\\ v\\ w\end{smallmatrix}\right)\right)
	\end{equation}
is of standard form \eqref{eqn_h_general_form}. Then we will need to study the existence of the limit of the above polynomial as $r\to R$ and check if the corresponding CCPSR surface is a homogeneous space. The process of determining $E(r)$ and studying limits will be split into several steps. We will assume from here on that $b\ne 0$ as the case $b=0$ will be dealt with in Section \ref{sect_b_0} for all dimensions $\dim(\mathcal{H})\geq 2$. Note that this in particular means $a\notin\left\{-\frac{2}{3\sqrt{3}},\ \frac{2}{3\sqrt{3}}\right\}$ or, equivalently, $R\notin\left\{\frac{\sqrt{3}}{2},\sqrt{3}\right\}$, since otherwise it is easy to check that $P_3$ attains the value $\frac{2}{3\sqrt{3}}$ at one of the two points $\left(\begin{smallmatrix}x\\v\\ w\end{smallmatrix}\right)=\left(\begin{smallmatrix}0\\0\\ \pm1\end{smallmatrix}\right)$, but said point is not a critical value of $P_3$ restricted to $\{\langle v,v\rangle + w^2=1\}$, thereby violating the property of the initial CCPSR manifold being closed, cf. Theorem \ref{thm_Cn}. This implies that the allowed values for $a$ and, equivalently, $R$ are precisely
	\begin{equation}\label{eqn_dim2_a_R_ranges}
		a\in \left(-\frac{2}{3\sqrt{3}},\ \frac{2}{3\sqrt{3}}\right),\quad\text{respectively}\quad R\in \left(\frac{\sqrt{3}}{2},\sqrt{3}\right).
	\end{equation}
We make the ansatz $E=E_1$ in \eqref{eqn_dim2_hE_def} with
	\begin{equation}\label{eqn_dim2_E1_def}
		E_1(r):=\beta^{\frac{1}{6}}\cdot(3-r^2)\cdot\mathbbm{1}
	\end{equation}
and, using
	\begin{equation*}
		\D h=(3x^2-v^2-w^2)\D x + (-2xv+3cv^2+2qvw+bw^2)\D v + (-2xw+qv^2+2bvw+3aw^2)\D w,
	\end{equation*}
obtain
	\begin{equation*}
		h_{E_1}\left(\left(\begin{smallmatrix} x\\ v\\ w\end{smallmatrix}\right)\right)=x^3-x (\lambda v^2 + \chi vw + \mu w^2)+\Theta_{v^3}v^3+\Theta_{v^2w}v^2w+\Theta_{vw^2}vw^2+\Theta_{w^3}w^3
	\end{equation*}
with
	\begin{align}
		\lambda&=-3b^2r^4+(1-qr)(3-r^2)^2,\label{eqn_dim2_lambda_def}\\
		\chi&=-18br\beta,\\
		\mu&=3(3-9ar+r^2)\beta,\label{eqn_dim2_mu_def}\\
		\Theta_{v^3}&=\left(-b^3r^6+br^2(3-r^2)^2+c(3-r^2)^3\right)\sqrt{\beta},\label{eqn_dim2_Theta_vvv}\\
		\Theta_{v^2w}&=\left(3b^2r^5(2-3ar)-r(2-3ar)(3-r^2)^2+q(3-r^2)^3\right)\sqrt{\beta},\label{eqn_dim2_Theta_vvw}\\
		\Theta_{vw^2}&=9b\left(3+ r^2-3ar^3\right)\beta\sqrt{\beta},\label{eqn_dim2_Theta_vww}\\
		\Theta_{w^3}&=\left((2-27a^2)r^3+27ar^2-18r+27a\right)\beta\sqrt{\beta}.\label{eqn_dim2_Theta_www}
	\end{align}
Thus, the task of finding the standard form of $h$ along $p(r)$ has reduced to diagonalizing the bilinear form $\lambda v^2 + \chi vw + \mu w^2$ in $v,w$, which is by assumption positive definite for all $r\in[0,R)$. We will, however, instead of first diagonalizing and then studying the limit as $r\to R$ of the corresponding polynomial, only rescale the diagonal parts of the the latter bilinear form first, then study the limit $r\to R$, and finally diagonalize. These steps will be discussed in detail below. The reason for this approach is that the calculations turned out to be shorter and less prone to errors. In order to follow this plan, we first need to check that the bilinear form,
	\begin{equation*}
		v^2+\frac{\chi}{\sqrt{\lambda \mu}}vw+w^2.
	\end{equation*}
converges to a positive definite bilinear as $r\to R$. This is equivalent to showing that
	\begin{equation}\label{eqn_dim2_chi_lambda_mu_limit}
		\lim\limits_{r\nearrow R}\frac{\chi}{\sqrt{\lambda \mu}}\in (-2,2),
	\end{equation}
which in particular requires the existence of the above limit for all allowed values of $a,b,c,q$, cf. equations \eqref{eqn_starting_P3_dim2} and \eqref{eqn_dim2_a_R_ranges}. From now on we will further assume that $b>0$, which can always be achieved by a sign-flip $v\to -v$ if necessary and, hence, is no restriction of generality under our previous assumptions.
The limit in \eqref{eqn_dim2_chi_lambda_mu_limit} motivates to study the cases
	\begin{equation}\label{eqn_dim2_lambda_0_limit}
		\lim\limits_{r\nearrow R}\lambda=0
	\end{equation}
and
	\begin{equation}\label{eqn_dim2_lambda_NONZERO_limit}
		\lim\limits_{r\nearrow R}\lambda>0
	\end{equation}
separately.

We will start with the case \eqref{eqn_dim2_lambda_0_limit}. For our calculations it is convenient to write $a$ in $R$-dependence, and we obtain from the condition $(1,0,R)^\mathrm{T}\in \partial(\mathbb{R}_{>0}\cdot\mathcal{H})$ that
	\begin{equation}\label{eqn_dim2_a_R_dependence}
		a=\frac{R^2-1}{R^3}.
	\end{equation}
The above assignment is a diffeomorphism of the respective ranges of $a$ and $R$, see equation \eqref{eqn_dim2_a_R_ranges}. We find that
	\begin{equation}\label{eqn_dim2_lamda_0_R_formula}
		\lim\limits_{r\nearrow R}\lambda=0\quad \Leftrightarrow \quad -3b^2R^4+(1-qR)(3-R^2)^2=0.
	\end{equation}
Now, we will now show that in the case \eqref{eqn_dim2_lambda_0_limit}
	\begin{equation}\label{eqn_dim2_lambda_0_derivative_nonzero}
		\left.\partial_r\lambda\right|_{r=R}<0.
	\end{equation}
We calculate 
	\begin{equation*}
		\left.\partial_r\lambda\right|_{r=R}=-12b^2R^3-q(3-R^2)^2-4R(1-qR)(3-R^2).
	\end{equation*}
Equation \eqref{eqn_dim2_lamda_0_R_formula} and $\sqrt{3}/2<R<\sqrt{3}$ imply $(1-qR)(3-R^2)=3b^2R^4/(3-R^2)$ and, hence, we obtain
	\begin{align}
		\left.\partial_r\lambda\right|_{r=R}&=-12b^2R^3-q(3-R^2)^2-\frac{12b^2R^5}{3-R^2}\notag\\
			&=\frac{-36b^2R^3-q(3-R^2)^3}{3-R^2}.\label{eqn_dim2_partial_lamda_rR}
	\end{align}
We can immediately exclude $\left.\partial_r\lambda\right|_{r=R}>0$ since otherwise $\lambda|_{r=R}=0$ would imply that there exist some $0<\widetilde{r}<R$, such that $\lambda|_{r=\widetilde{r}}<0$. Suppose that $\left.\partial_r\lambda\right|_{r=R}=0$. By equation \eqref{eqn_dim2_partial_lamda_rR} and $\sqrt{3}/2<R<\sqrt{3}$ this is fulfilled if and only if $-36b^2R^3-q(3-R^2)^3=0$. This, together with the assumption $b>0$, implies $q<0$. Next, we check that
	\begin{align*}
		0&=12\lambda|_{r=R}-R\left(-36b^2R^3-q(3-R^2)^3\right)\\
			&=(3-R^2)^2\left(qR(-9-R^2)+12\right).
	\end{align*}
But $q<0$ implies $qR(-9-R^2)+12>0$, hence the above equation cannot be satisfied. This shows that \eqref{eqn_dim2_lambda_0_derivative_nonzero} holds in the case \eqref{eqn_dim2_lambda_0_limit}.

Now consider the formula for $\mu$ \eqref{eqn_dim2_mu_def} and observe with \eqref{eqn_dim2_a_R_dependence} that
	\begin{equation*}
		\left.(3-9aR+R^2)\right|_{a=\frac{R^2-1}{R^3}}=0\quad \Leftrightarrow \quad R=\pm\sqrt{3},
	\end{equation*}
and that $\left.(3-9aR+R^2)\right|_{a=\frac{R^2-1}{R^3}}>0$ for all $R\in\left(\sqrt{3}/2,\sqrt{3}\right)$. Also recall that $\beta=1-r^2+ar^3$ has a simple zero in $r=R$ for all $R\in\left(\sqrt{3}/2,\sqrt{3}\right)$. This shows the existence of the limit on the left hand side of \eqref{eqn_dim2_chi_lambda_mu_limit}, since as we have seen both $\mu$ \eqref{eqn_dim2_mu_def} and $\lambda$ \eqref{eqn_dim2_lambda_def} have a simple zero in $r=R$. In order to show that it is an element of $(-2,2)$, it suffices to study the equation
	\begin{equation}\label{eqn_dim2_chilimit_not2}
		\lim\limits_{r\nearrow R}\frac{\chi}{\sqrt{\lambda \mu}}=-2
	\end{equation}
and show that this cannot be satisfied under the assumption that $\mathcal{H}\subset\mathbb{R}^3$ is closed. Note that the sign on the right hand side of equation \eqref{eqn_dim2_chilimit_not2} comes from our assumption $b>0$. Observe that
	\begin{equation}\label{eqn_dim2_q_lamda_zero}
		\lambda|_{r=R}=0\quad \Leftrightarrow\quad q=\frac{-3b^2R^4+(3-R^2)^2}{R(3-R^2)^2}.
	\end{equation}
We calculate
	\begin{equation}\label{eqn_dim2_chilambdamu_sqarelim}
		\lim\limits_{r\nearrow R}\frac{\chi^2}{\lambda \mu}=\frac{108b^2R^4}{3b^2R^4(9+R^2)+(3-R^2)^3}.
	\end{equation}
Note that in the above calculation, no usage of L'H\^opital's rule is needed. Solving $\frac{108b^2R^4}{3b^2R^4(9+R^2)+(3-R^2)^3}=4$ for $b$ symbolically, we find
	\begin{equation*}
		\lim\limits_{r\nearrow R}\frac{\chi^2}{\lambda \mu}=4\quad \Leftrightarrow\quad b=\pm\frac{(3-R^2)\sqrt{3R^2-9}}{3R^3}.
	\end{equation*}
But $R\in\left(\frac{\sqrt{3}}{2},\sqrt{3}\right)$, so the term $\sqrt{3R^2-9}$ would be imaginary. This is a contradiction to $b$ being real. We conclude together with $\frac{\chi}{\sqrt{\lambda \mu}}\in(-2,0)$ for all $r\in[0,R)$, $\mathcal{H}\subset\mathbb{R}^3$ being closed, and using $b>0$, that \eqref{eqn_dim2_chi_lambda_mu_limit} is indeed satisfied, that is
	\begin{equation}\label{eqn_dim2_eqn_dim2_chilambdamu_limit_formula}
		\lim\limits_{r\nearrow R}\frac{\chi}{\sqrt{\lambda \mu}}=-\frac{6}{\sqrt{9+R^2+\frac{1}{3b^2R^4}(3-R^2)^3}}\in (-2,0).
	\end{equation}
	
At this point, we have seen that $a$ and $q$ are completely determined by the choice of $b$ and $R$ in the case \eqref{eqn_dim2_lambda_0_limit}. To see that $c$ is also determined by that choice, we show that
	\begin{equation}\label{eqn_dim2_Theta_at_R_lambda_0}
		\left.\left(\frac{\Theta_{v^3}}{\sqrt{\beta}}\right)\right|_{r=R}=\left.\left(\frac{\Theta_{v^2w}}{\sqrt{\beta}}\right)\right|_{r=R}=\left.\left(\frac{\Theta_{vw^2}}{\sqrt{\beta}}\right)\right|_{r=R}=\left.\left(\frac{\Theta_{w^3}}{\sqrt{\beta}}\right)\right|_{r=R}=0,
	\end{equation}
i.e. the polynomial factors of the $\Theta_I$, $I\in\left\{v^3,\ v^2w,\ vw^2,\ w^3\right\}$ \eqref{eqn_dim2_Theta_vvv}--\eqref{eqn_dim2_Theta_www}, must vanish at $r=R$. Note that it would be enough to prove \eqref{eqn_dim2_Theta_at_R_lambda_0} in order to show that $c$ is determined by $b,R$, but checking the other three equations is a good sanity test for the consistency of our calculations. Verifying $\left.\left(\frac{\Theta_{v^2w}}{\sqrt{\beta}}\right)\right|_{r=R}=\left.\left(\frac{\Theta_{vw^2}}{\sqrt{\beta}}\right)\right|_{r=R}=\left.\left(\frac{\Theta_{w^3}}{\sqrt{\beta}}\right)\right|_{r=R}=0$ with equations \eqref{eqn_dim2_a_R_dependence} and \eqref{eqn_dim2_q_lamda_zero} is thus just a direct computation. Next, suppose that $\left.\left(\frac{\Theta_{v^3}}{\sqrt{\beta}}\right)\right|_{r=R}\ne 0$ for some choice of $b$ and $R$. We will use Theorem \ref{thm_Cn} to show that this would be a contradiction to $\mathcal{H}\subset\mathbb{R}^3$ being closed. By rescaling $h_{E_1}$, cf. \eqref{eqn_dim2_hE_def} and \eqref{eqn_dim2_E1_def}, in $v$ with $1/{\sqrt{\lambda}}$ and in $w$ with $1/\sqrt{\mu}$, we obtain
	\begin{equation}\label{eqn_dim2_hE2}
		h_{E_2}\left(\left(\begin{smallmatrix}x\\ v\\ w\end{smallmatrix}\right)\right)=x^3-x\left(v^2+\frac{\chi}{\sqrt{\lambda\mu}}vw+w^2\right)+ \frac{\Theta_{v^3}}{\lambda\sqrt{\lambda}}v^3 + \frac{\Theta_{v^2w}}{\lambda\sqrt{\mu}}v^2w + \frac{\Theta_{vw^2}}{\mu\sqrt{\lambda}}vw^2 + \frac{\Theta_{w^3}}{\mu\sqrt{\mu}}w^3,
	\end{equation}
where $E_2(r)=E_1(r)\cdot\left(\begin{smallmatrix}1/{\sqrt{\lambda}} & 0\\ 0 & 1/{\sqrt{\mu}}\end{smallmatrix}\right)$. We have shown above that the bilinear form $v^2+\frac{\chi}{\sqrt{\lambda\mu}}vw+w^2$ converges to a positive definite bilinear form as $r\to R$. This of course implies that the eigenvalues of the corresponding symmetric matrix
	\begin{equation}\label{eqn_dim2_postrescale_bilform}
		\left(\begin{matrix}1	&	\frac{\chi}{2\sqrt{\lambda\mu}}\\
			\frac{\chi}{2\sqrt{\lambda\mu}}	&	1
			\end{matrix}\right)
	\end{equation}
are positive and bounded from above by $2-\varepsilon$ for some $\varepsilon>0$ for all $r\in[0,R]$. Hence,
	\begin{equation}\label{eqn_dim2_case_lambda_zero_bilinear_estimate_postrescale}
		v^2+\frac{\chi}{\sqrt{\lambda\mu}}vw+w^2\leq (2-\varepsilon)(v^2+w^2)\quad \forall r\in[0,R]\ \forall v,w\in\mathbb{R}.
	\end{equation}
Recall that by proving \eqref{eqn_dim2_lambda_0_derivative_nonzero} we have shown that $\lambda$ has a simple zero in $r=R$ in the case \eqref{eqn_dim2_lambda_0_limit}, and since $a<\frac{2}{3\sqrt{3}}$ \eqref{eqn_dim2_a_R_ranges}, $\beta$ also has a simple zero in $r=R$, cf. Lemma \ref{lem_beta_zeros}. Hence,
	\begin{equation*}
		\left.\left(\frac{\Theta_{v^3}}{\sqrt{\beta}}\right)\right|_{r=R}\ne 0\quad \Leftrightarrow\quad 	\left.\left(\frac{\Theta_{v^3}}{\sqrt{\lambda}}\right)\right|_{r=R}\ne 0.
	\end{equation*}
But this means that the $v^3$-term in $h_{E_2}$ \eqref{eqn_dim2_hE2} will not converge as $r\to R$. In order to bring \eqref{eqn_dim2_hE2} to standard form, we need to diagonalize \eqref{eqn_dim2_postrescale_bilform}. Using \eqref{eqn_dim2_case_lambda_zero_bilinear_estimate_postrescale} together with the fact that, as we have just seen, $\frac{\Theta_{v^3}}{\lambda\sqrt{\lambda}}$ does not converge as $r\to R$, we conclude that after the diagonalization of \eqref{eqn_dim2_postrescale_bilform} via some matrix $U\in\mathrm{GL}(2)$ which brings $h_{E_2}$ to standard form, cf. Proposition \ref{prop_standard_form}, we can still find $\widetilde{r}\in(0,R)$ and a vector $\gamma\in\mathbb{R}^2$ of Euclidean unit length, such that
	\begin{equation*}
		\left.\left(\frac{\Theta_{v^3}}{\lambda\sqrt{\lambda}}v^3 + \frac{\Theta_{v^2w}}{\lambda\sqrt{\mu}}v^2w + \frac{\Theta_{vw^2}}{\mu\sqrt{\lambda}}vw^2 + \frac{\Theta_{w^3}}{\mu\sqrt{\mu}}w^3\right)\right|_{\left(\begin{smallmatrix}v\\ w\end{smallmatrix}\right)=U\gamma,\ r=\widetilde{r}}>\frac{2}{3\sqrt{3}}.
	\end{equation*}
We now use Theorem \ref{thm_Cn} and find that this a contradiction to $\mathcal{H}\subset\mathbb{R}^3$ being closed. Hence, $\left.\left(\frac{\Theta_{v^3}}{\sqrt{\beta}}\right)\right|_{r=R}=0$ must hold, and we obtain the formula
	\begin{equation}\label{eqn_dim2_c_bR_dep}
		c=\frac{bR^2 (b^2R^4-(3-R^2)^2)}{(3-R^2)^3}.
	\end{equation}
Summarizing up to this point, we have obtained in the case \eqref{eqn_dim2_lambda_0_limit} that $a$ \eqref{eqn_dim2_a_R_dependence}, $q$ \eqref{eqn_dim2_q_lamda_zero}, and $c$ \eqref{eqn_dim2_c_bR_dep} are completely determined by the choices of $b$ and $R$. Now we are at the point where we can calculate the limit of the polynomial $h_{E_2}$ \eqref{eqn_dim2_hE2} as $r\to R$. In order to avoid errors and make our calculations easier to verify, we will split this calculation into a few separate steps. Using the determined values for $a$, $q$, and $c$, we first calculate
	\begin{align}\label{eqn_dim2_chi_by_sqrtlambdamu_limit}
		\lim\limits_{r\nearrow R}\frac{\chi}{\sqrt{\lambda\mu}}&=\frac{-6\sqrt{3}bR^2}{\sqrt{3b^2R^4(9+R^2)+(3-R^2)^3}},\\
		\label{eqn_dim2_sqrtlambda_by_sqrtmu_limit}
		\lim\limits_{r\nearrow R}\frac{\sqrt{\lambda}}{\sqrt{\mu}}&=\frac{R\sqrt{3b^2R^4(9+R^2)+(3-R^2)^3}}{\sqrt{3}(3-R^2)^2},\\
		\label{eqn_dim2_sqrtmu_by_sqrtbeta_limit}
		\lim\limits_{r\nearrow R}\frac{\sqrt{\mu}}{\sqrt{\beta}}&=\frac{\sqrt{3}(3-R^2)}{R},
	\end{align}
where, as for equation \eqref{eqn_dim2_chilambdamu_sqarelim}, no L'H\^opital's rule is needed. Care with the signs of the formulas is however necessary. We further find 
	\begin{align}
		\lim\limits_{r\nearrow R}\frac{\Theta_{v^3}}{\beta\sqrt{\beta}}&=\frac{6bR^2(3b^2R^4-(3-R^2)^2)}{(3-R^2)^2},\label{eqn_dim2_theta_vvv_by_beta_limit}\\
		\lim\limits_{r\nearrow R}\frac{\Theta_{v^2w}}{\beta\sqrt{\beta}}&=\frac{2(3b^2R^4(-9+R^2)+(3-R^2)^3)}{R(3-R^2)},\label{eqn_dim2_theta_vvw_by_beta_limit}\\
		\lim\limits_{r\nearrow R}\frac{\Theta_{vw^2}}{\beta\sqrt{\beta}}&=18b(3-R^2),\label{eqn_dim2_theta_vww_by_beta_limit}\\
		\lim\limits_{r\nearrow R}\frac{\Theta_{w^3}}{\beta\sqrt{\beta}}&=\frac{-2(3-R^2)^3}{R^3}.\label{eqn_dim2_theta_www_by_beta_limit}
	\end{align}
Hence,
	\begin{align}
		\widehat{h}\left(\left(\begin{smallmatrix}x\\ v\\ w\end{smallmatrix}\right)\right):=\lim\limits_{r\nearrow R}	h_{E_2}\left(\left(\begin{smallmatrix}x\\ v\\ w\end{smallmatrix}\right)\right)
		&=x^3-x (v^2+\zeta vw+w^2)\notag\\
		&\quad +\Xi_{v^3}v^3+\Xi_{v^2w}v^2w+\Xi_{vw^2}vw^2+\Xi_{w^3}w^3\label{eqn_dim2_hE2_limit}
	\end{align}
with
	\begin{align}
		\zeta&= \frac{-6\sqrt{3}bR^2}{\sqrt{3b^2R^4(9+R^2)+(3-R^2)^3}}\\
		\Xi_{v^3}&= \frac{-6bR^2(3b^2R^4(-3+R^2)+(3-R^2)^3)}{\left(\sqrt{3b^2R^4(9+R^2)+(3-R^2)^3}\right)^3}\label{eqn_dim2_Xi_vvv},\\
		\Xi_{v^2w}&= \frac{2(3b^2R^4(-9+R^2)+(3-R^2)^3)}{\sqrt{3}\left(3b^2R^4(9+R^2)+(3-R^2)^3\right)}\label{eqn_dim2_Xi_vvw},\\
		\Xi_{vw^2}&= \frac{6bR^2}{\sqrt{3b^2R^4(9+R^2)+(3-R^2)^3}}\label{eqn_dim2_Xi_vww},\\
		\Xi_{w^3}&= \frac{-2}{3\sqrt{3}}\label{eqn_dim2_Xi_www}.
	\end{align}
Next, we need to diagonalize and normalize the bilinear form $v^2+\frac{-6\sqrt{3}bR^2}{\sqrt{3b^2R^4(9+R^2)+(3-R^2)^3}}vw+w^2$. This will be done in two steps. First we check that $\left\{\left(\begin{smallmatrix}1\\ 1\end{smallmatrix}\right),\ \left(\begin{smallmatrix}1\\ -1\end{smallmatrix}\right)\right\}$ is an orthogonal basis of $\mathbb{R}^2$ with respect to that bilinear form. For the following calculations, let
	\begin{equation*}
		\rho:=\sqrt{3b^2R^4(9+R^2)+(3-R^2)^3},\quad \kappa:=bR^2,
	\end{equation*}
so that
	\begin{equation*}
		\zeta=\frac{-6\sqrt{3}\kappa}{\rho},\quad\Xi_{v^3}=\frac{-6\kappa(\rho^2-36\kappa^2)}{\rho^3},\quad\Xi_{v^2w}=\frac{2(\rho^2-54\kappa^2)}{\sqrt{3}\rho^2},\quad\Xi_{vw^2}=\frac{6\kappa}{\rho}.
	\end{equation*}
We set
	\begin{equation*}
		E_3:=\left(\begin{array}{rr}1 & 1\\ 1 & -1\end{array}\right)
	\end{equation*}
and find
	\begin{align*}
		\left(\widehat{h}\circ\left(\begin{smallmatrix}\mathbbm{1} & \\ & E_3\end{smallmatrix}\right)\right)\left(\left(\begin{smallmatrix}x\\ v\\ w\end{smallmatrix}\right)\right)
			&= x^3 - x\left( \frac{2\rho-6\sqrt{3}\kappa}{\rho} v^2 +\frac{2\rho+6\sqrt{3}\kappa}{\rho}w^2\right)\\
			&\quad + \left(\Xi_{v^3}+\Xi_{v^2w}+\Xi_{vw^2}+\Xi_{w^3}\right)v^3\\
			&\quad + \left(3\Xi_{v^3}+\Xi_{v^2w}-\Xi_{vw^2}-3\Xi_{w^3}\right)v^2w\\
			&\quad + \left(3\Xi_{v^3}-\Xi_{v^2w}-\Xi_{vw^2}+3\Xi_{w^3}\right)vw^2\\
			&\quad + \left(\Xi_{v^3}-\Xi_{v^2w}+\Xi_{vw^2}-\Xi_{w^3}\right)w^3.
	\end{align*}
We further calculate
	\begin{align*}
		\Xi_{v^3}+\Xi_{v^2w}+\Xi_{vw^2}+\Xi_{w^3}&=\frac{4\left(\sqrt{3}\rho^3-81\sqrt{3}\kappa^2\rho+486\kappa^3\right)}{9\rho^3},\\
		3\Xi_{v^3}+\Xi_{v^2w}-\Xi_{vw^2}-3\Xi_{w^3}&=\frac{4\left(\rho^2-27\kappa^2\right)\left(\sqrt{3}\rho-18\kappa\right)}{3\rho^3},\\
		3\Xi_{v^3}-\Xi_{v^2w}-\Xi_{vw^2}+3\Xi_{w^3}&=\frac{4\left(\rho^2-27\kappa^2\right)\left(-\sqrt{3}\rho-18\kappa\right)}{3\rho^3},\\
		\Xi_{v^3}-\Xi_{v^2w}+\Xi_{vw^2}-\Xi_{w^3}&=\frac{4\left(-\sqrt{3}\rho^3+81\sqrt{3}\kappa^2\rho+486\kappa^3\right)}{9\rho^3}.
	\end{align*}
In order to bring $\left(\widehat{h}\circ\left(\begin{smallmatrix}\mathbbm{1} & \\ & E_3\end{smallmatrix}\right)\right)\left(\left(\begin{smallmatrix}x\\ v\\ w\end{smallmatrix}\right)\right)$ to standard form \eqref{eqn_h_general_form}, we only need to rescale $v$ and $w$. Before doing that, we transform the coefficients $\rho$ and $\kappa$ via
	\begin{equation*}
		s=\rho+3\sqrt{3}\kappa,\quad t=\rho-3\sqrt{3}\kappa.
	\end{equation*}
Note that $\rho=(s+t)/2$ and $\kappa=(s-t)/(6\sqrt{3})$, and also that $s$ and $t$ are necessarily positive in case \eqref{eqn_dim2_lambda_0_limit}. We finally obtain
	\begin{align}
		h_{(s,t)}\left(\left(\begin{smallmatrix}x\\ v\\ w\end{smallmatrix}\right)\right):=&\!\; \left(\widehat{h}\circ\left(\begin{smallmatrix}\mathbbm{1} & \\ & E_3\end{smallmatrix}\right)\right)\left(\left(\begin{matrix}x\\ \frac{\sqrt{s+t}}{2\sqrt{t}}v\\ \frac{\sqrt{s+t}}{2\sqrt{s}}w\end{matrix}\right)\right)\\
		=&\ x^3-x\left(v^2+w^2\right)\\
		&\quad + \frac{2}{3\sqrt{3}}\cdot\frac{(3s-t)\sqrt{t}}{\sqrt{s+t}^3}v^3\\
		&\quad + \frac{2}{\sqrt{3}}\cdot\frac{(-s+3t)\sqrt{s}}{\sqrt{s+t}^3}v^2w\\
		&\quad + \frac{2}{\sqrt{3}}\cdot\frac{(-3s+t)\sqrt{t}}{\sqrt{s+t}^3}vw^2\\
		&\quad + \frac{2}{3\sqrt{3}}\cdot\frac{(s-3t)\sqrt{s}}{\sqrt{s+t}^3}w^3.
	\end{align}
We see that the above polynomial $h_{(s,t)}$ is in standard form \eqref{eqn_h_general_form}, independent of the values for $s,t$ for which the corresponding maximal connected PSR surface containing the point $(x,v,w)^\mathrm{T}=(1,0,0)^\mathrm{T}$ is closed. Now we observe that all the coefficients in the $P_3$-part of $h_{(s,t)}$ are homogeneous functions of degree $0$ in $s,t$. Since, as mentioned above, $s>0$ and $t>0$ for all allowed values of $b$ and $R$, we deduce that the possible outcomes of $h_{(s,t)}$ in dependence of the starting data $b,R$ can be described by fixing $s=1$. This means that
	\begin{align*}
		&\left.\left\{h_{(s,t)}\ \right|\ s,t\text{ obtained from allowed starting data }b,R\right\}\\
		=\ &\left.\left\{h_{(1,t)}\ \right|\ t\text{ obtained from allowed starting data }b,R\text{ so that }s=1\right\}.
	\end{align*}
Hence, in order to prove Theorem \ref{thm_limit_behaviour_ccpsr} for $\dim(\mathcal{H})=2$ in the case \eqref{eqn_dim2_lambda_0_limit}, it suffices to prove that for all $t>0$ such that the connected component $\mathcal{H}_t$ of $\left\{h_{(1,t)}=1\right\}\subset\mathbb{R}^3$ that contains the point $(x,v,w)^\mathrm{T}=(1,0,0)^\mathrm{T}$ is a CCPSR surface, $\mathcal{H}_t$ is also a homogeneous space. We show that $\mathcal{H}_t$ is indeed a CCPSR surface for all $t>0$ and that it is a homogeneous space in one step. Consider for $t>0$ the linear transformation
	\begin{equation*}
		E_4(t):=\frac{1}{\sqrt{1+t}}\left(\begin{matrix}\sqrt{t} & 1 \\ -1 & \sqrt{t}\end{matrix}\right)\in\mathrm{O}(2).
	\end{equation*}
For any fixed $t>0$, we transform $h_{(1,t)}$ and obtain
	\begin{equation}
			\left(h_{(1,t)}\circ\left(\begin{smallmatrix}\mathbbm{1} & \\ & E_4(t)\end{smallmatrix}\right)\right)\left(\left(\begin{smallmatrix}x\\ v\\ w\end{smallmatrix}\right)\right)
				=x^3-x(v^2+w^2)-\frac{2}{3\sqrt{3}}v^3+\frac{2}{\sqrt{3}}vw^2.\label{eqn_dim2_bne0_lambda0}
	\end{equation}
We conclude that for all $t>0$, $\mathcal{H}_t$ is equivalent to \cite[Thm.\,1\,a)]{CDL}, cf. \cite[Ex.\,3.2]{L2}, which is a CCPSR surface and furthermore a homogeneous space. This finishes our treatment of case \eqref{eqn_dim2_lambda_0_limit}. Note that this case actually occurs. We leave it to the reader to check that one example corresponds to moving to infinity in the $z$-direction in Example \ref{ex_motivating}. A hint how to show this can be found by studying the high-dimensional analogue in Example \ref{ex_ker01_n_geq_3}.

Next, we will deal with the case \eqref{eqn_dim2_lambda_NONZERO_limit}. Recall that we assume $a<\frac{2}{3\sqrt{3}}$ and $b\ne0$. For $b=0$ see the next Section \ref{sect_b_0}. We proceed similarly to the case \eqref{eqn_dim2_lambda_0_limit} and arrive at equations \eqref{eqn_dim2_lambda_def}--\eqref{eqn_dim2_Theta_www}. The difference now is that we assume $\lambda|_{r=R}=-3b^2R^4+(1-qR)(3-R^2)^2>0$. Hence, we obtain
	\begin{equation*}
		\lim\limits_{r\nearrow R}\frac{\chi}{\sqrt{\lambda\mu}}=0
	\end{equation*}
instead of the non-trivial formula \eqref{eqn_dim2_chi_by_sqrtlambdamu_limit} in the case \eqref{eqn_dim2_lambda_0_limit}. The formulas \eqref{eqn_dim2_sqrtmu_by_sqrtbeta_limit}, \eqref{eqn_dim2_theta_vww_by_beta_limit}, and \eqref{eqn_dim2_theta_www_by_beta_limit} however also hold in the case \eqref{eqn_dim2_lambda_NONZERO_limit}, which implies that
	\begin{align*}
		\lim\limits_{r\nearrow R}\frac{\Theta_{w^3}}{\mu\sqrt{\mu}}&=\frac{-2}{3\sqrt{3}},\\
		\lim\limits_{r\nearrow R}\frac{\Theta_{vw^2}}{\mu\sqrt{\lambda}}&=0.
	\end{align*}
Furthermore, we immediately see from equation \eqref{eqn_dim2_Theta_vvv} for $\Theta_{v^3}$, which has the overall factor $\sqrt{\beta}$, that
	\begin{equation*}
		\lim\limits_{r\nearrow R}\frac{\Theta_{v^3}}{\lambda\sqrt{\lambda}}=0.
	\end{equation*}
We now calculate that, independent of the possible starting data $c,q,b,R$ that satisfy \eqref{eqn_dim2_lambda_NONZERO_limit},
	\begin{equation*}
		\lim\limits_{r\nearrow R}\frac{\Theta_{v^2w}}{\lambda\sqrt{\mu}}=-\frac{1}{\sqrt{3}}.
	\end{equation*}
Summarizing, we have shown that in the case \eqref{eqn_dim2_lambda_NONZERO_limit} the limit of $h_{E_2}$ \eqref{eqn_dim2_hE2} as $r\to R$ is given by
	\begin{equation}\label{eqn_dim2_limitpoly_lambda_nonzero_limit}
			\lim\limits_{r\nearrow R}h_{E_2}\left(\left(\begin{smallmatrix}x\\ v\\ w\end{smallmatrix}\right)\right)=x^3-x(v^2+w^2) -\frac{1}{\sqrt{3}}v^2w - \frac{2}{3\sqrt{3}} w^3.
	\end{equation}
By comparing with \cite[Thm.\,1\,b)]{CDL} and \cite[Ex.\,3.2]{L2}, we see that the connected component $\widetilde{\mathcal{H}}$ of $\{\widetilde{h}=1\}\subset\mathbb{R}^3$ that contains the point $(x,v,w)^\mathrm{T}=(1,0,0)^\mathrm{T}$ is a homogeneous CCPSR surface as claimed. Hence, Theorem \ref{thm_limit_behaviour_ccpsr}, respectively Thm. \eqref{thm_limit_behaviour_ccpsr_ii_a}, holds in this case. Summarising this section of the proof of Theorem \ref{thm_limit_behaviour_ccpsr}, we have obtained the following characterisation for limit geometries:
		\begin{equation*}
			\begin{array}{|c|c|c|c|}
				\hline \dim(\mathcal{H}) & b & \lambda|_{r=R} & \text{eqn.}\\
				\hline \hline 2 & \overset{\vphantom{\frac{1}{R}}}{}\underset{}{\vphantom{0}}\ne0  & 0 & \eqref{eqn_dim2_bne0_lambda0}\\
				\hline 2 & \overset{\vphantom{\frac{1}{R}}}{}\underset{}{\vphantom{0}}\ne0 & >0 & \eqref{eqn_dim2_limitpoly_lambda_nonzero_limit}\\
				\hline
			\end{array}
		\end{equation*}
\subsection{$\dim(\mathcal{H})= 2$, $b=0$}\label{sect_b_0}
	Note that $\mathcal{H}$ being closed and $b=0$ implies that $R\in\left[\frac{\sqrt{3}}{2},\sqrt{3}\right]$ or, equivalently, $a\in\left[-\frac{2}{3\sqrt{3}},\frac{2}{3\sqrt{3}}\right]$ are the allowed values for $R$ and $a$, respectively. Since $b=0$, equations \eqref{eqn_dim2_lambda_def}--\eqref{eqn_dim2_Theta_www} are of the simpler form
		\begin{align}
			\lambda &=(1-qr)(3-r^2)^2,\label{eqn_dim2_lambda_def_b0}\\
			\chi &=0,\\
			\mu &=3(3-9ar+r^2)\beta,\label{eqn_dim2_mu_def_b0}\\
			\Theta_{v^3} &=c(3-r^2)^3\sqrt{\beta},\label{eqn_dim2_Theta_vvv_b0}\\
			\Theta_{v^2w} &=\left(-r(2-3ar)(3-r^2)^2+q(3-r^2)^3\right)\sqrt{\beta},\label{eqn_dim2_Theta_vvw_b0}\\
			\Theta_{vw^2} &=0,\label{eqn_dim2_Theta_vww_b0}\\
			\Theta_{w^3} &=\left((2-27a^2)r^3+27ar^2-18r+27a\right)\beta\sqrt{\beta}.\label{eqn_dim2_Theta_www_b0}
		\end{align}
	As for $b\ne 0$, we will separately consider the cases $\lim\limits_{r\nearrow R}\lambda=0$ \eqref{eqn_dim2_lambda_0_limit} and $\lim\limits_{r\nearrow R}\lambda>0$ \eqref{eqn_dim2_lambda_NONZERO_limit}.

	First assume that $\lim\limits_{r\nearrow R}\lambda=0$. We start with the case $R<\sqrt{3}$. In this case $q=\frac{1}{R}$ and, hence, $\lambda\big|_{b=0}=\left(1-\frac{r}{R}\right)(3-r^2)^2$ has a simple zero in $r=R$. Since $R<\sqrt{3}$, $\beta$ also has a simple zero in $r=R$, where we recall that $\beta$ has a simple zero in $r=R$ if $a\in\left[-\frac{2}{3\sqrt{3}},\frac{2}{3\sqrt{3}}\right)$. Hence, $c=0$ must hold, since otherwise
		\begin{equation*}
			\lim\limits_{r\nearrow R}\frac{\Theta_{v^3} }{\lambda\sqrt{\lambda}}\rightarrow \pm\infty,
		\end{equation*}
	which would violate $\mathcal{H}$ being closed by Theorem \ref{thm_Cn}. We use $R<\sqrt{3}$ again and find that equations \eqref{eqn_dim2_sqrtmu_by_sqrtbeta_limit} and \eqref{eqn_dim2_theta_www_by_beta_limit} hold independently of the choice for $b$, so we obtain
		\begin{equation*}
			\lim\limits_{r\nearrow R}\frac{\Theta_{w^3}}{\mu\sqrt{\mu}}=-\frac{2}{3\sqrt{3}}.
		\end{equation*}
	It remains to calculate the limit $\lim\limits_{r\nearrow R}\frac{\Theta_{v^2w}}{\mu\sqrt{\lambda}}$. As above, we find that equations \eqref{eqn_dim2_sqrtlambda_by_sqrtmu_limit} and \eqref{eqn_dim2_theta_vvw_by_beta_limit} hold for $b=0$, and together with equation \eqref{eqn_dim2_sqrtmu_by_sqrtbeta_limit} we obtain
		\begin{equation*}
			\lim\limits_{r\nearrow R}\frac{\Theta_{v^2w}}{\mu\sqrt{\lambda}}=\frac{2}{\sqrt{3}}.
		\end{equation*}
	Hence, the limit polynomial $\lim\limits_{r\nearrow R}h_{E_2}$ \eqref{eqn_dim2_hE2} is given by
		\begin{equation}
			\lim\limits_{r\nearrow R}h_{E_2}\left(\left(\begin{smallmatrix}x\\ v\\ w\end{smallmatrix}\right)\right)=x^3-x(v^2+w^2) +\frac{2}{\sqrt{3}}v^2w - \frac{2}{3\sqrt{3}} w^3,\label{eqn_2dim_lambda0_Rlesssqrt3_limitpoly}
		\end{equation}
	which coincides with Thm. \ref{thm_limit_behaviour_ccpsr} \eqref{thm_limit_behaviour_ccpsr_ii_b}, showing that Theorem \ref{thm_limit_behaviour_ccpsr} holds in this case.
	
	Next consider the case $R=\sqrt{3}$. In this case we cannot use equations \eqref{eqn_dim2_chi_by_sqrtlambdamu_limit}--\eqref{eqn_dim2_theta_www_by_beta_limit} for $b=0$, as these hold only for $R<\sqrt{3}$, and we recall that $\mathcal{H}$ being closed and $b\ne 0$ imply $R<\sqrt{3}$. Note that in the case $R=\sqrt{3}$, $\mu=3(\sqrt{3}-r)^2\beta$. We find using equations \eqref{eqn_dim2_lambda_def_b0}--\eqref{eqn_dim2_Theta_www_b0}
		\begin{equation}\label{eqn_Theta_v3_w3_2dim_b0}
			\frac{\Theta_{v^3}}{\lambda\sqrt{\lambda}}=\frac{c\sqrt{\beta}}{(1-qr)\sqrt{1-qr}},\quad\frac{\Theta_{w^3}}{\mu\sqrt{\mu}}=\frac{2}{3\sqrt{3}}.
		\end{equation}
	This implies that if $q=\frac{1}{\sqrt{3}}$, $c=0$ must hold. Otherwise $\lim\limits_{r\nearrow R}\frac{\Theta_{v^3}}{\lambda\sqrt{\lambda}}=\pm\infty$ since then both $\sqrt{\beta}$ and $1-qr$ have a simple zero in $r=\sqrt{3}=R$. Furthermore, for $q=\frac{1}{\sqrt{3}}$ one finds that the term
		\begin{equation}\label{eqn_Theta_v2w_2dim_b0_q1sqrt3}
			\frac{\Theta_{v^2w}}{\lambda\sqrt{\mu}}=\frac{1}{\sqrt{3}}
		\end{equation}
	is constant. This is expected, as for $c=0$, $q=\frac{1}{\sqrt{3}}$ the initial CCPSR surface is already homogeneous, cf. \cite[Thm.\,1\,b)]{CDL} and \cite[Ex.\,3.2]{L2}. Hence, the limit polynomial $\lim\limits_{r\nearrow R}h_{E_2}$ \eqref{eqn_dim2_hE2} for $R=\sqrt{3}$, $q=\frac{1}{\sqrt{3}}$, is given by
		\begin{equation}
			\lim\limits_{r\nearrow R}h_{E_2}\left(\left(\begin{smallmatrix}x\\ v\\ w\end{smallmatrix}\right)\right)=x^3-x(v^2+w^2) +\frac{1}{\sqrt{3}}v^2w + \frac{2}{3\sqrt{3}} w^3,\label{eqn_dim2_b0_Rsqrt3_q1bysqrt3_limitpoly}
		\end{equation}
	and thus is after a sign-flip in $w$ equivalent to Thm. \ref{thm_limit_behaviour_ccpsr} \eqref{thm_limit_behaviour_ccpsr_ii_a}. Now assume $q<\frac{1}{\sqrt{3}}$. The two equations in \eqref{eqn_Theta_v3_w3_2dim_b0} are also valid in this case, which in particular implies
		\begin{equation*}
			\lim\limits_{r\nearrow R}\frac{\Theta_{v^3}}{\lambda\sqrt{\lambda}}=0
		\end{equation*}
	independent of $c$. The difference to the above case is now that $\frac{\Theta_{v^2w}}{\lambda\sqrt{\mu}}$ is not constant, that is equation \eqref{eqn_Theta_v2w_2dim_b0_q1sqrt3} does not hold. Instead we obtain for the limit
		\begin{equation*}
			\lim\limits_{r\nearrow R}	\frac{\Theta_{v^2w}}{\lambda\sqrt{\mu}}=-\frac{2}{\sqrt{3}}
		\end{equation*}
	independent of $q<\frac{1}{\sqrt{3}}$ and, hence, the limit polynomial $\lim\limits_{r\nearrow R}h_{E_2}$ \eqref{eqn_dim2_hE2} for $R=\sqrt{3}$, $q<\frac{1}{\sqrt{3}}$ is given by
		\begin{equation}\label{eqn_dim2_b0_Rsqrt3_qlessmax}
			\lim\limits_{r\nearrow R}h_{E_2}\left(\left(\begin{smallmatrix}x\\ v\\ w\end{smallmatrix}\right)\right)=x^3-x(v^2+w^2) -\frac{2}{\sqrt{3}}v^2w + \frac{2}{3\sqrt{3}} w^3.
		\end{equation}
	This shows that the Theorem \ref{thm_limit_behaviour_ccpsr} holds in this case by comparing with \cite[Thm.\,1\,a)]{CDL} and \cite[Ex.\,3.2]{L2}.
	
	Now assume that $\lim\limits_{r\nearrow R}\lambda>0$. This implies that $R<\sqrt{3}$. Thus we can proceed exactly as for $b\ne 0$ and obtain that the limit polynomial $	\lim\limits_{r\nearrow R}h_{E_2}$ \eqref{eqn_dim2_hE2} is given by
		\begin{equation}
			\lim\limits_{r\nearrow R}h_{E_2}\left(\left(\begin{smallmatrix}x\\ v\\ w\end{smallmatrix}\right)\right)=x^3-x(v^2+w^2) -\frac{1}{\sqrt{3}}v^2w - \frac{2}{3\sqrt{3}} w^3,\label{eqn_dim2_b0_Rlesssqrt3_lambdapositive_limitpoly}
		\end{equation}
	cf. equation \eqref{eqn_dim2_limitpoly_lambda_nonzero_limit}. By comparing with \cite[Thm.\,1\,b)]{CDL} and \cite[Ex.\,3.2]{L2}, it follows that Theorem \ref{thm_limit_behaviour_ccpsr} holds in this case as well.
	
	This finishes our treatment of CCPSR surfaces and shows that Theorem \ref{thm_limit_behaviour_ccpsr} holds in dimension $2$, i.e. Thm. \ref{thm_limit_behaviour_ccpsr} \eqref{thm_limit_behaviour_ccpsr_ii}. Summarising this part of the proof, we have obtained the following characterisation of limit geometries:
		\begin{equation*}
			\begin{array}{|c|c|c|c|c|}
				\hline \dim(\mathcal{H}) & b & R & q & \text{eqn.}\\
				\hline \hline 2 & \overset{\vphantom{\frac{1}{R}}}{}\underset{}{\vphantom{0}}0 & <\sqrt{3} & \frac{1}{R} & \eqref{eqn_2dim_lambda0_Rlesssqrt3_limitpoly}\\
				\hline 2 & \overset{\vphantom{\frac{1}{R}}}{}\underset{}{\vphantom{0}}0 & \sqrt{3} & \frac{1}{R}=\frac{1}{\sqrt{3}} & \eqref{eqn_dim2_b0_Rsqrt3_q1bysqrt3_limitpoly}\\
				\hline 2 & \overset{\vphantom{\frac{1}{R}}}{}\underset{}{\vphantom{0}}0 & \sqrt{3} & <\frac{1}{R}=\frac{1}{\sqrt{3}} & \eqref{eqn_dim2_b0_Rsqrt3_qlessmax}\\
				\hline 2 & \overset{\vphantom{\frac{1}{R}}}{}\underset{}{\vphantom{0}}0 & <\sqrt{3} & <\frac{1}{R} & \eqref{eqn_dim2_b0_Rlesssqrt3_lambdapositive_limitpoly}\\
				\hline
			\end{array}
		\end{equation*}

\subsection{$\dim(\mathcal{H})\geq 3$, $b\ne 0$}\label{section_dimH_geq3_bnon0}
Let now $n=\dim(\mathcal{H})\geq 3$ and denote $(y_1,\ldots,y_{n-1},y_n)=(v_1,\ldots,v_{n-1},w)$, so that $v=(v_1,\ldots,v_{n-1})$ denotes linear coordinates on $\mathbb{R}^{n-1}$ and we have
\begin{equation}\label{eqn_xvw_coords}
	\left(\begin{matrix}x\\ y\end{matrix}\right)=\left(\begin{matrix} x\\ v\\ w\end{matrix}\right).
\end{equation}
Let further
\begin{equation}\label{eqn_h_n_geq_3_first_form}
	h\left(\left(\begin{smallmatrix}x\\ v\\ w\end{smallmatrix}\right)\right):=x^3-x(\langle v,v\rangle +w^2) + C(v) + Q(v)w+b\langle \gamma,v\rangle w^2 + aw^3.
\end{equation}
In the above equation \re{eqn_h_n_geq_3_first_form}, $a\in\mathbb{R}$, $b\in\mathbb{R}\setminus\{0\}$, the vector $\gamma\in\mathbb{R}^{n-1}$ is a unit vector, i.e. fulfils $\langle \gamma,\gamma\rangle=1$, $Q:\mathbb{R}^{n-1}\to \mathbb{R}$ is a quadratic homogeneous polynomial, and $C:\mathbb{R}^{n-1}\to\mathbb{R}$ is a cubic homogeneous polynomial. For the case $b=0$ see Section \ref{subsection_dimHgeq3_b0}. The $P_3$-term in the above equation is thus given by
	\begin{equation}
		P_3\left(\left(\begin{smallmatrix}v\\w\end{smallmatrix}\right)\right)=C(v) + Q(v)w+b\langle \gamma,v\rangle w^2 + aw^3.\label{eqn_P3_term_bnon0}
	\end{equation}	
In the following we will identify $Q(v)$ with the corresponding bilinear form $Q(v,v)$ and $C(v)$ with the corresponding trilinear form $C(v,v,v)$, so that we can write e.g. $\D Q_v=2Q(v,dv)$. Furthermore, we will assume that $Q(v)$ is of the form
	\begin{equation*}
		Q(v)=\sum\limits_{i=1}^{n-1}q_iv_i^2,\quad q_i\in\mathbb{R},
	\end{equation*}
which can always be achieved via an orthogonal transformation with respect to the standard Euclidean scalar product in the $v$-coordinates $\langle\cdot,\cdot\rangle$ of the $v$-coordinates. Note that these type of transformations commute with the process of obtaining limit geometries as stated in Theorem \ref{thm_limit_behaviour_ccpsr}. Similar to the case $\dim(\mathcal{H})=2$, we want to apply Proposition \ref{prop_standard_form} to a point of the form
	\begin{equation}\label{eqn_p_xvw}
		p=\left(\begin{smallmatrix} p_x\\ p_v\\ p_w\end{smallmatrix}\right)=\beta^{-\frac{1}{3}}\left(\begin{smallmatrix} 1\\ 0\\ r\end{smallmatrix}\right)\in\mathcal{H},
	\end{equation}
where $p_x=x(p)$, $p_v=v(p)$, and $p_w=w(p)$. This explains the emphasis on the coordinate labelling in \re{eqn_xvw_coords}, and we note that $h$ as in \re{eqn_h_n_geq_3_first_form} is of the most general form, cf. equation \re{eqn_h_general_form}. In the following, let $R$ be defined analogously to the $2$-dimensional case in equation \eqref{eqn_dim2_R_def}, i.e. $r=R$ is the smallest positive solution of $h\left(\left(\begin{smallmatrix}1\\0\\r\end{smallmatrix}\right)\right)=0$. Note that $b\ne 0$ implies $R\in\left(\frac{\sqrt{3}}{2},\sqrt{3}\right)$ or, equivalently, $a\in\left(-\frac{2}{3\sqrt{3}},\frac{2}{3\sqrt{3}}\right)$, since if $R$ attained the values $\frac{\sqrt{3}}{2}$ or $\sqrt{3}$ it is easy to check that $\left|P_3\left(\left(\begin{smallmatrix}0\\1\end{smallmatrix}\right)\right)\right|=\frac{2}{3\sqrt{3}}$ but $\left(\left(\begin{smallmatrix}0\\1\end{smallmatrix}\right)\right)$ would not be a critical point of $P_3|_{\{\langle v,v\rangle + w^2=1\}}$ by $b\ne 0 $.
We calculate
	\begin{align}\label{eqn_dh_xvw}
		\D h&=\left(3x^2-\langle v,v\rangle -w^2\right)\D x\notag\\
		&\quad-2x\langle v,\D v\rangle + bw^2\langle \gamma,\D v\rangle + 2Q(v,\D v)+3C(v,v,\D v)\notag\\
		&\quad+\left(-2xw+3aw^2+2bw\langle \gamma,v\rangle + Q(v)\right)\D w
	\end{align}
and
	\begin{align}\label{eqn_hessh_alternative_xvw}
		\partial^2 h &= 6x\, \D x^2 -2x\langle \D v,\D v\rangle + 2Q(\D v,\D v) + 6C(v,\D v,\D v) + (-2x+6aw+2b\langle \gamma,v\rangle)\D w^2\notag\\
		&\quad-4\,\D x\langle v,\D v\rangle -4w\,\D x\D w + 4bw\langle \gamma,\D v\rangle \D w + 4Q(v,\D v)\D w.
	\end{align}
Evaluating $\D h$ \re{eqn_dh_xvw} and $\partial^2h$ \re{eqn_hessh_alternative_xvw} at the point $p$ \re{eqn_p_xvw}, we obtain
\begin{equation}
	\D h_p=\beta^{-\frac{2}{3}}\cdot\left((3-r^2)\D x + br^2\langle \gamma,\D v\rangle + (-2r+3ar^2)\D w\right)
\end{equation}
and
\begin{align}
	\partial^2h_p&=\beta^{-\frac{1}{3}}\cdot\left(6r\,\D x^2-2\langle \D v,\D v\rangle + 2Q(\D v,\D v)\right.\notag\\
		&\quad\left. + (-2+6ar)\D w^2-4r\,\D x\D w+4br\langle \gamma,\D v\rangle \D w\right),
\end{align}
respectively. Next, similar as in equations \eqref{eqn_dim2_hE_def} and \eqref{eqn_dim2_E1_def} for the case $\dim(\mathcal{H})=2$, we define
$h_E$ as in \eqref{eqn_dim2_hE_def} with the difference that now $\left.-\frac{\partial_v h}{\partial_x h}\right|_{p(r)}=\left.\left(-\partial_{v_1}h/\partial_x h,\ldots,-\partial_{v_{n-1}}h/\partial_x h\right)\right|_{p(r)}$ has $(n-1)$ entries instead of just one, and begin with the ansatz in \eqref{eqn_dim2_hE_def}
	\begin{equation}\label{eqn_dimgen_E1_def}
		E(r):=\beta^{\frac{1}{6}}\cdot(3-r^2)\cdot\mathbbm{1}.
	\end{equation}
We obtain
	\begin{equation}\label{eqn_dimgen_hE1_prelimit}
		h_{E}\left(\left(\begin{smallmatrix} x\\ v\\ w\end{smallmatrix}\right)\right)=x^3-x\left(\boldsymbol{\lambda}+\boldsymbol{\chi}+\boldsymbol{\mu}\right)+\boldsymbol{\Theta}_{v^3} +\boldsymbol{\Theta}_{v^2w} +\boldsymbol{\Theta}_{vw^2} +\boldsymbol{\Theta}_{w^3}
	\end{equation}
with
	\begin{align}
		\boldsymbol\lambda&=-3b^2r^4\langle\gamma,v\rangle^2+\sum\limits_{i=1}^{n-1}(1-q_ir)(3-r^2)^2v_i^2,\label{eqn_dimgen_lambda_def}\\
		\boldsymbol\chi&=-18br\beta\langle\gamma,v\rangle w,\label{eqn_dimgen_chi_def}\\
		\boldsymbol\mu&=3(3-9ar+r^2)\beta w^2,\label{eqn_dimgen_mu_def}\\
		\boldsymbol\Theta_{v^3}&=\left(-b^3r^6\langle\gamma,v\rangle^3+br^2(3-r^2)^2\langle v,v\rangle\langle\gamma,v\rangle+(3-r^2)^3C(v)\right)\sqrt{\beta},\label{eqn_dimgen_Theta_vvv}\\
		\boldsymbol\Theta_{v^2w}&= \left(3b^2r^5(2-3ar)\langle\gamma,v\rangle^2-r(2-3ar)(3-r^2)^2\langle v,v\rangle+(3-r^2)^3Q(v)\right)w \sqrt{\beta},\label{eqn_dimgen_Theta_vvw}\\
		\boldsymbol\Theta_{vw^2}&=9b\left(3+ r^2-3ar^3\right)\langle\gamma,v\rangle w^2\beta\sqrt{\beta},\label{eqn_dimgen_Theta_vww}\\
		\boldsymbol\Theta_{w^3}&=\left((2-27a^2)r^3+27ar^2-18r+27a\right)w^3\beta\sqrt{\beta}.\label{eqn_dimgen_Theta_www}
	\end{align}
By comparing equations \eqref{eqn_dimgen_lambda_def}--\eqref{eqn_dimgen_Theta_www} with equations \eqref{eqn_dim2_lambda_def}--\eqref{eqn_dim2_Theta_www} from the $2$-dimensional case, we see that $\boldsymbol\chi=\chi\langle\gamma,v\rangle w$, $\boldsymbol\mu=\mu w^2$, $\boldsymbol\Theta_{vw^2}= \Theta_{vw^2}\langle\gamma,v\rangle w^2$, and $\boldsymbol\Theta_{w^3}=\Theta_{w^3}w^3$. In particular, equation \eqref{eqn_dim2_sqrtmu_by_sqrtbeta_limit}, which is equivalent to
	\begin{equation}
		\lim\limits_{r\nearrow R}\frac{\sqrt{\beta}}{\sqrt{\mu}}=\frac{R}{\sqrt{3}(3-R^2)},\label{eqn_sqrt_mu_by_sqrt_beta_limit}
	\end{equation}
still holds in our present setting, and using equations \eqref{eqn_dim2_theta_vww_by_beta_limit} and \eqref{eqn_dim2_theta_www_by_beta_limit}, where we recall that these two results depend \textit{only} on the values of $a=(R^2-1)/R^3$ and $b$, it is also clear that
	\begin{align}
		\lim\limits_{r\nearrow R}\frac{\boldsymbol\Theta_{vw^2}}{\beta\sqrt{\beta}}&=\lim\limits_{r\nearrow R}\frac{\Theta_{vw^2}}{\beta\sqrt{\beta}}\langle\gamma,v\rangle w^2 =18b(3-R^2)\langle\gamma,v\rangle w^2,\label{eqn_dimgen_theta_vww_by_beta_limit}\\
		\lim\limits_{r\nearrow R}\frac{\boldsymbol\Theta_{w^3}}{\beta\sqrt{\beta}}&=\lim\limits_{r\nearrow R}\frac{\Theta_{w^3}}{\beta\sqrt{\beta}}w^3 =\frac{-2(3-R^2)^3}{R^3}w^3.\label{eqn_dimgen_theta_www_by_beta_limit}
	\end{align}
Equations \eqref{eqn_dim2_theta_vvv_by_beta_limit} and \eqref{eqn_dim2_theta_vvw_by_beta_limit} do not have easy analogues for $\dim(\mathcal{H})\geq 3$. For $\dim(\mathcal{H})=2$ we differentiated between the cases \eqref{eqn_dim2_lambda_0_limit} and \eqref{eqn_dim2_lambda_NONZERO_limit}. Since $\boldsymbol\lambda$ is a positive definite bilinear form for each $r\in[0,R)$ and not simply a function depending on a certain choice of variables, we will have to make a more delicate differentiation for $\dim(\mathcal{H})\geq 3$. In the following, we will denote by
	\begin{equation}\label{eqn_dimgen_q_def}
		q=\left(\begin{matrix}q_1 & & \\ & \ddots & \\ & & q_{n-1}\end{matrix}\right),
	\end{equation}
and by expressions of the form $F(q)$, where $F$ is a smooth function defined on some subset of $\mathbb{R}_{>0}$ and possibly depending on the other variables $a,b,r,R$, we denote
	\begin{equation}\label{eqn_dimgen_Fq_def}
		F(q):=\left(\begin{matrix}F(q_1) & & \\ & \ddots & \\ & & F(q_{n-1})\end{matrix}\right).
	\end{equation}
With the above notation, $\boldsymbol\lambda$ \eqref{eqn_dimgen_lambda_def} can be written as
	\begin{equation*}
		\boldsymbol\lambda=-3b^2r^4\langle\gamma,v\rangle^2+(3-r^2)^2\langle v,(1-qr)v\rangle.
	\end{equation*}

We start with the case $\boldsymbol\lambda|_{r=R}>0$. Let $B\in\mathrm{GL}(n-1)$, such that $B^\ast \boldsymbol\lambda|_{r=R}=\langle\cdot,\cdot\rangle$, i.e.
	\begin{equation}\label{eqn_dimgen_lambda_nokernel_diag}
		-3b^2R^4\langle\gamma,Bv\rangle^2+(3-R^2)^2\langle Bv,(1-qR)Bv\rangle=\langle v,v\rangle\quad \forall v\in\mathbb{R}^{n-1}.
	\end{equation}
We immediately see that in this case
	\begin{equation*}
		\lim\limits_{r\nearrow R}(B\times(\mathbbm{1}_\mathbb{R}/\sqrt{\mu}))^\ast\boldsymbol\chi=0.
	\end{equation*}
In the above equation the map $\mathbbm{1}_\mathbb{R}$ acts on the $w$-variable, so that formally $(\mathbbm{1}_\mathbb{R}/\sqrt{\mu})^\ast\boldsymbol\mu=w^2$ for all $r\in[0,R)$ and also for the limit $r\to R$. Since the term $\boldsymbol\Theta_{v^3}$ \eqref{eqn_dimgen_Theta_vvv} does not depend on the variable $w$, it follows that
	\begin{equation*}
		\lim\limits_{r\nearrow R}\left(B\times(\mathbbm{1}_\mathbb{R}/\sqrt{\mu})\right)^\ast\boldsymbol\Theta_{v^3}=0.
	\end{equation*}
Now we find using equations \eqref{eqn_dim2_sqrtmu_by_sqrtbeta_limit} and \eqref{eqn_dimgen_lambda_nokernel_diag} that
	\begin{align}
		&\lim\limits_{r\nearrow R}\left(B\times(\mathbbm{1}_\mathbb{R}/\sqrt{\mu})\right)^\ast\boldsymbol\Theta_{v^2w}\notag\\
		&=\frac{1}{\sqrt{3}}\left( 3b^2R^4\langle\gamma,Bv\rangle^2-(3-R^2)^2\langle Bv,Bv\rangle+R(3-R^2)^2\langle Bv,qBv\rangle\right)w\notag\\
		&=-\frac{1}{\sqrt{3}}\langle v,v\rangle w.\label{eqn_Theta_vvw_pullback_limit_lambda_greater_0}
	\end{align}
With equations \eqref{eqn_dim2_sqrtmu_by_sqrtbeta_limit} and \eqref{eqn_dimgen_theta_vww_by_beta_limit} and the assumption $\boldsymbol\lambda|_{r=R}>0$ we deduce that
	\begin{equation*}
		\lim\limits_{r\nearrow R}\left(B\times(\mathbbm{1}_\mathbb{R}/\sqrt{\mu})\right)^\ast\boldsymbol\Theta_{vw^2}=0,
	\end{equation*}
and equations \eqref{eqn_dim2_sqrtmu_by_sqrtbeta_limit} and \eqref{eqn_dimgen_theta_www_by_beta_limit} imply
	\begin{equation}\label{eqn_Theta_w3_limit_R_less_sqrt3}
		\lim\limits_{r\nearrow R}\left(B\times(\mathbbm{1}_\mathbb{R}/\sqrt{\mu})\right)^\ast\boldsymbol\Theta_{w^3}=-\frac{2}{3\sqrt{3}}w^3.
	\end{equation}
Summarizing, we have shown that under the assumption $\boldsymbol\lambda|_{r=R}>0$ the limit of $h_{E_1}$ \eqref{eqn_dimgen_hE1_prelimit} as $r\to R$ is given by
	\begin{align}
		\widetilde{h}\left(\left(\begin{smallmatrix}x\\ v\\ w\end{smallmatrix}\right)\right):=&\lim\limits_{r\nearrow R}\left(\mathbbm{1}_\mathbb{R}\times B\times(\mathbbm{1}_\mathbb{R}/\sqrt{\mu})\right)^\ast 	h_{E}\left(\left(\begin{smallmatrix} x\\ v\\ w\end{smallmatrix}\right)\right)\notag\\
		=&\ x^3-x\left(\langle v,v\rangle +w^2\right) -\frac{1}{\sqrt{3}}\langle v,v\rangle w -\frac{2}{3\sqrt{3}} w^3.\label{eqn_dimgen_lambda_nokernel_hE1_limit}
	\end{align}
By comparing \eqref{eqn_dimgen_lambda_nokernel_hE1_limit} to the corresponding limit polynomial in the $2$-dimensional case \eqref{eqn_dim2_limitpoly_lambda_nonzero_limit} when equation \eqref{eqn_dim2_lambda_NONZERO_limit} holds, we find that the assumption $\boldsymbol\lambda|_{r=R}>0$ generalises this case as one might expect. The connected component $\widetilde{\mathcal{H}}\subset\left\{\widetilde{h}=1\right\}$ that contains the point $\left(\begin{smallmatrix}x\\ v\\ w\end{smallmatrix}\right)=\left(\begin{smallmatrix}1\\ 0\\ 0\end{smallmatrix}\right)$ is a homogeneous CCPSR manifold, for details see the proof of \cite[Prop.\,6.9]{L1}.

The next step is dealing with the cases $\dim\ker\boldsymbol\lambda|_{r=R}>0$. By the hyperbolicity of the points in $\mathcal{H}$ we know that $q_i\leq\frac{1}{R}$ for all $1\leq i\leq n-1$, since otherwise it is easy to see that there would exist $\widetilde{r}\in[0,R)$, such that $\boldsymbol\lambda|_{r=\widetilde{r}}$ has at least one negative eigenvalue. From here on we assume without loss of generality that $q_1\geq\ldots\geq q_{n-1}$. We will proceed as follows.

If $\frac{1}{R}>q_1$, $\dim\ker\boldsymbol\lambda|_{r=R}>0$ can only be true if $\dim\ker\boldsymbol\lambda|_{r=R}=1$. This follows from the assumption $q_i\geq q_{i+1}$ for all $1\leq i\leq n-2$ and $\boldsymbol\lambda$ being positive definite as a bilinear form for all $r\in[0,R)$. These cases are surprisingly complex to deal with and will be treated as the first step.

First we transform our linear coordinates $v$ on $\mathbb{R}^{n-1}$ and $w$ on $\mathbb{R}$ for $r\in[0,R)$ via
	\begin{align}
		v&=\frac{1}{(3-r^2)\sqrt{1-qr}}\eta,\label{eqn_ker_0_1_eta_trafo}\\ \overline{w}&=\frac{1}{\sqrt{\mu}}w,\label{eqn_ker_0_1_overlinew_trafo}
	\end{align}
transforming $\boldsymbol\lambda=:\boldsymbol\lambda_{vv}$ \eqref{eqn_dimgen_lambda_def} and $\boldsymbol\chi=:\boldsymbol\chi_{vw}$ \eqref{eqn_dimgen_chi_def} to
	\begin{align}
		\boldsymbol\lambda_{\eta\eta}&=-\frac{3b^2r^4}{(3-r^2)^2}\left\langle\frac{1}{\sqrt{1-qr}}\gamma,\eta\right\rangle^2 + \langle\eta,\eta\rangle,\label{eqn_dimgen_lambda_ker_0_1_1st_trafo}\\
		\boldsymbol\chi_{\eta\overline{w}}&=-\frac{18br\beta}{(3-r^2)\sqrt{\mu}}\left\langle\frac{1}{\sqrt{1-qr}}\gamma,\eta\right\rangle\overline{w},\label{eqn_dimgen_chi_ker_0_1_1st_trafo}
	\end{align}
respectively. Note that under these transformations, $\boldsymbol\mu=:\boldsymbol\mu_{ww}$ \eqref{eqn_dimgen_mu_def} transforms to $\boldsymbol\mu_{\overline{w}\overline{w}}=\overline{w}^2$. Furthermore note that the linear transformation in $v$ \eqref{eqn_ker_0_1_eta_trafo} extends to $r=R$ since $b\ne 0$ by assumption, which means that $R<\sqrt{3}$ must be fulfilled, and because $q_i<\frac{1}{R}$ for all $1\leq i\leq n-1$ by assumption. It follows that $\dim\ker\boldsymbol\lambda|_{r=R}=1$ if and only if
	\begin{equation}
		\left\|\frac{1}{\sqrt{1-qR}}\gamma\right\|^2=\frac{(3-R^2)^2}{3b^2R^4}.\label{eqn_ker_0_1_condition}
	\end{equation}
Now $\frac{1}{\sqrt{1-qr}}\in\mathrm{GL}(n-1)$ for all $r\in[0,R]$ and $\gamma\ne 0$ implies that we can choose a smooth map $B:[0,R]\to\mathrm{SO}(n-1)$, such that
	\begin{equation}
		B(r)^\mathrm{T}\frac{1}{\sqrt{1-qr}}\gamma=\left\|\frac{1}{\sqrt{1-qr}}\gamma\right\|\partial_{n-1}\label{eqn_ker_0_1_B_gamma_condition}
	\end{equation}
for all $r\in[0,R]$, where we have identified $\partial_{n-1}$ with the $(n-1)$-th Euclidean unit vector in $\mathbb{R}^{n-1}$ in the $v$-coordinates. The above equation \eqref{eqn_ker_0_1_B_gamma_condition} is equivalent to
	\begin{equation}
		\frac{1}{\sqrt{1-qr}}B(r)\partial_{n-1}=\left\|\frac{1}{\sqrt{1-qr}}\gamma\right\|^{-1}\frac{1}{1-qr}\gamma\label{eqn_ker_0_1_B_gamma_condition_alternative}
	\end{equation}
for all $r\in[0,R]$.
Denote $\rho^\mathrm{T}=(\rho_1,\ldots,\rho_{n-2})$ and transform $\eta$ via
	\begin{equation}
		\eta=B(r)\left(\begin{matrix}\rho\\ \tau\end{matrix}\right),\quad \rho\in\mathbb{R}^{n-2},\ \tau\in\mathbb{R}.\label{eqn_ker_0_1_rho_tau_trafo}
	\end{equation}
Note that above transformation \eqref{eqn_ker_0_1_rho_tau_trafo} does depend on $r\in[0,R]$. We now obtain
	\begin{align}
		B(r)^*\boldsymbol\lambda_{\eta\eta}&=:\boldsymbol\lambda_{\rho\rho}+\boldsymbol\lambda_{\tau\tau}\notag\\
			&=\langle\rho,\rho\rangle + \left(1-\frac{3b^2r^4}{(3-r^2)^2}\left\|\frac{1}{\sqrt{1-qr}}\gamma\right\|^2\right)\tau^2,\label{eqn_lambda_rhorho_tautau}\\
		B(r)^*\boldsymbol\chi_{\eta\overline{w}}&=:\boldsymbol\chi_{\tau\overline{w}}\notag\\
			&=-\frac{18br\beta}{(3-r^2)\sqrt{\mu}}\left\|\frac{1}{\sqrt{1-qr}}\gamma\right\|\tau\overline{w}.\label{eqn_chi_tau_overlinew}
	\end{align}
To reduce the symbols needed in the following calculations, we define
	\begin{equation}
		\nu=\nu(r):=\frac{\sqrt{3}br^2}{(3-r^2)\sqrt{1-qr}}\gamma,\quad r\in[0,R],\label{eqn_nu_def}
	\end{equation}
with $\|\nu\|^2=\frac{3b^2r^4}{(3-r^2)^2}\left\|\frac{1}{\sqrt{1-qr}}\gamma\right\|^2$. Note that equation \eqref{eqn_ker_0_1_condition} implies that $\|\nu\|^2\to 1$ as $r\to R$. Before considering any kind of limit we need to introduce one more transformation. For $r\in[0,R]$, let
	\begin{equation}
		\tau=\frac{1}{\sqrt{1-\|\nu\|^2}}\alpha,\quad \alpha\in\mathbb{R}.\label{eqn_ker_0_1_tau_alpha_trafo}
	\end{equation}
Under \eqref{eqn_ker_0_1_tau_alpha_trafo}, $\boldsymbol\lambda_{\rho\rho}+\boldsymbol\lambda_{\tau\tau}$ and $\boldsymbol\chi_{\tau\overline{w}}$ transform to
	\begin{align}
		\boldsymbol\lambda_{\rho\rho}+\boldsymbol\lambda_{\alpha\alpha}&=\langle\rho,\rho\rangle + \alpha^2,\label{eqn_lambda_rhorho_alphaalpha}\\
		\boldsymbol\chi_{\alpha\overline{w}}&=-\mathrm{sgn}(b)\frac{6\sqrt{3}\beta\|\nu\|}{r\sqrt{\mu}\sqrt{1-\|\nu\|^2}}\alpha\overline{w},\label{eqn_chi_alpha_overlinew}
	\end{align}
respectively. Our plan now is to take the limit $r\to R$ of our already transformed polynomial $h_E$ \eqref{eqn_dimgen_hE1_prelimit} in the coordinates $\left(x,\rho^\mathrm{T},\alpha,\overline{w}\right)^\mathrm{T}$ \underline{before} transforming these coordinates any further, and in the step thereafter bring the so obtained limit polynomial to standard form. In order for this to make sense in our setting, we need to check that for all admissible initial data with $\mathcal{H}$ closed together with the assumptions we are currently working with,
	\begin{equation}
		\lim\limits_{r\nearrow R}\left|\boldsymbol\chi_{\alpha\overline{w}}\right|<2|\alpha\overline{w}|.\label{eqn_chi_alpha_overlinew_limit_condition}
	\end{equation}
Note that in comparison with the $2$-dimensional cases, the this step is analogous to proving that the right hand side of equation \eqref{eqn_dim2_chilambdamu_sqarelim} cannot attain the value $4$. In order to calculate the left hand side of equation \eqref{eqn_chi_alpha_overlinew_limit_condition}, we recall that the formula for the limit of $\frac{\sqrt{\beta}}{\sqrt{\mu}}$ as $r\to R$ \eqref{eqn_sqrt_mu_by_sqrt_beta_limit}. From equation \eqref{eqn_chi_alpha_overlinew} we now see that it remains to determine the limit of $\frac{\sqrt{\beta}}{\sqrt{1-\|\nu\|^2}}$ as $r\to R$. To do so we will use L'H\^opital's rule. We calculate
	\begin{equation*}
		\partial_r\|\nu\|^2 = \frac{3b^2r^3(9+r^2)}{(3-r^2)^3}\left\langle\gamma,\frac{1}{1-qr}\gamma\right\rangle + \frac{3b^2r^3}{(3-r^2)^2}\left\langle\gamma,\frac{1}{(1-qr)^2}\gamma\right\rangle,
	\end{equation*}
and obtain using equation \eqref{eqn_ker_0_1_condition}
	\begin{align}
		\left.\partial_r\|\nu\|^2\right|_{r=R}&=\frac{9+R^2}{R(3-R^2)} + \frac{3b^2R^3}{(3-R^2)^2}\left\langle\gamma,\frac{1}{(1-qR)^2}\gamma\right\rangle\notag\\
		&=\frac{9+R^2}{R(3-R^2)} + \frac{3b^2R^3}{(3-R^2)^2}\left\|\frac{1}{1-qR}\gamma\right\|^2.\label{eqn_partial_r_nusquare_limit}
	\end{align}
Note that the second equality in \eqref{eqn_partial_r_nusquare_limit} is justified since $\frac{1}{1-qR}$, viewed as a bilinear form, is positive definite since by assumption $q_i<\frac{1}{R}$ for all $1\leq i\leq n-1$. Hence, equation \eqref{eqn_partial_r_nusquare_limit} and $R<\sqrt{3}$ implies $\left.\partial_r\|\nu\|^2\right|_{r=R}>0$. We quickly check that $\left.\partial_r\beta\right|_{r=R}=\frac{-3+R^2}{R}$ and obtain
	\begin{equation}
		\lim\limits_{r\nearrow R}\frac{\sqrt{\beta}}{\sqrt{1-\|\nu\|^2}} = \frac{3-R^2}{\sqrt{9+R^2+\frac{3b^2R^4}{3-R^2}\left\|\frac{1}{1-qR}\gamma\right\|^2}}.\label{eqn_sqrt_beta_sqrt_1minusnusquare_limit}
	\end{equation}
Hence,
	\begin{equation}
		0<\lim\limits_{r\nearrow R}\frac{\sqrt{\beta}}{\sqrt{1-\|\nu\|^2}}<\frac{3-R^2}{\sqrt{9+R^2}}.\label{eqn_estimate_sqrt_beta_sqrt_1minusnusquare_limit}
	\end{equation}
We obtain using $\|\nu\|^2\to 1$ as $r\to R$ and equation \eqref{eqn_sqrt_mu_by_sqrt_beta_limit} that
	\begin{align}
		\lim\limits_{r\nearrow R} \boldsymbol\chi_{\alpha\overline{w}}  &= -\frac{6\;\mathrm{sgn}(b)}{3-R^2}\lim\limits_{r\nearrow R}\frac{\sqrt{\beta}}{\sqrt{1-\|\nu\|^2}} \alpha \overline{w} \notag\\
		&=-\frac{6\;\mathrm{sgn}(b)}{\sqrt{9+R^2+\frac{3b^2R^4}{3-R^2}\left\|\frac{1}{1-qR}\gamma\right\|^2}} \alpha \overline{w},\label{eqn_chi_alpha_overlinew_limit}
	\end{align}
and can now estimate using \eqref{eqn_estimate_sqrt_beta_sqrt_1minusnusquare_limit}
	\begin{align}
		\lim\limits_{r\nearrow R}\left|\boldsymbol\chi_{\alpha\overline{w}}\right| <\frac{6}{\sqrt{9+R^2}}|\alpha \overline{w}|<2|\alpha \overline{w}|.\label{eqn_chi_alpha_overlinew_limit_estimate}
	\end{align}
The reader might want to compare the above equation \eqref{eqn_chi_alpha_overlinew_limit_estimate} with equation \eqref{eqn_dim2_eqn_dim2_chilambdamu_limit_formula} from the $2$-dimensional case. Thus we have shown that \eqref{eqn_chi_alpha_overlinew_limit_condition} holds for all admissible initial data for $\mathcal{H}$ closed in our current case. Hence, in order to calculate the limit polynomial, we can proceed analogously to the $2$-dimensional case. That means that we take the limit in our current coordinates $\left(x,\rho^\mathrm{T},\alpha,w\right)^\mathrm{T}$ first, corresponding to equation \eqref{eqn_dim2_hE2_limit} in the $2$-dimensional case, and \underline{then} bring the polynomial obtained through that process to standard form, which leads to equation \eqref{eqn_dim2_limitpoly_lambda_nonzero_limit} in the $2$-dimensional case. In the latter step we of course need a transformation that leaves the point $\left(x,\rho^\mathrm{T},\alpha,w\right)^\mathrm{T}=\left(1,0^\mathrm{T},0,0\right)^\mathrm{T}$ invariant.

In order to proceed as described above, we need to transform the $\boldsymbol\Theta$-terms \eqref{eqn_dimgen_Theta_vvv}--\eqref{eqn_dimgen_Theta_vww} in our current coordinates $\left(x,\rho^\mathrm{T},\alpha,w\right)^\mathrm{T}$. In order to make clear what is happening, we will first give a name to the transformation leading to our currently in-use coordinates,
	\begin{equation}
		L:[0,R)\to\mathrm{GL}(n-1),\quad L(r)\left(\begin{matrix}\rho\\ \alpha \end{matrix}\right)=\frac{1}{(3-r^2)\sqrt{1-qr}}B(r)\left(\begin{matrix}\rho\\ \frac{\alpha}{\sqrt{1-\|\nu\|^2}}\end{matrix}\right).\label{eqn_L_trafo}
	\end{equation}
Then, using
	\begin{align*}
		L(r)^*\langle\gamma,v\rangle &= \frac{1}{3-r^2}\cdot\frac{\left\|\frac{1}{\sqrt{1-qr}}\gamma\right\|}{\sqrt{1-\|\nu\|^2}}\alpha,\\
		L(r)^*\langle v,v\rangle &= \frac{1}{(3-r^2)^2}\left(\left\|\frac{1}{\sqrt{1-qr}}B(r)\left(\begin{smallmatrix}\rho\\ 0\end{smallmatrix}\right)\right\|^2\right.\\
		&\quad +\frac{2}{\sqrt{1-\|\nu\|^2}}\left\|\frac{1}{\sqrt{1-qr}}\gamma\right\|^{-1} \left\langle \frac{1}{\sqrt{1-qr}}B(r)\left(\begin{smallmatrix}\rho\\ 0\end{smallmatrix}\right),\frac{1}{1-qr}\gamma\right\rangle\alpha\\
		&\quad\left.+\ \frac{1}{1-\|\nu\|^2}\left\|\frac{1}{\sqrt{1-qr}}\gamma\right\|^{-2}\left\|\frac{1}{1-qr}\gamma\right\|^2\alpha^2\right),
	\end{align*}
we obtain for the pullback of $\boldsymbol\Theta_{v^3}$ \eqref{eqn_dimgen_Theta_vvv}
	\begin{equation*}
		L(r)^*\boldsymbol\Theta_{v^3}=\boldsymbol\Theta_{\rho^3}+\boldsymbol\Theta_{\rho^2\alpha}+\boldsymbol\Theta_{\rho\alpha^2}+\boldsymbol\Theta_{\alpha^3}
	\end{equation*}
with
	\begin{align}
		\boldsymbol\Theta_{\rho^3}&=C\left(\frac{1}{\sqrt{1-qr}}B(r)\left(\begin{smallmatrix}\rho\\ 0\end{smallmatrix}\right)\right)\sqrt{\beta},\label{eqn_dimgen_Theta_rhorhorho}\\
		\boldsymbol\Theta_{\rho^2\alpha}&=\left(\vphantom{\left\|\frac{1}{\sqrt{1-qr}}\gamma\right\|^{-1}} \frac{br^2}{3-r^2}\left\|\frac{1}{\sqrt{1-qr}}\gamma\right\|\left\|\frac{1}{\sqrt{1-qr}}B(r)\left(\begin{smallmatrix}\rho\\ 0\end{smallmatrix}\right)\right\|^2\right.\notag\\
		&\quad\left. + 3\left\|\frac{1}{\sqrt{1-qr}}\gamma\right\|^{-1}C\left(\frac{1}{\sqrt{1-qr}}B(r)\left(\begin{smallmatrix}\rho\\ 0\end{smallmatrix}\right),\frac{1}{\sqrt{1-qr}}B(r)\left(\begin{smallmatrix}\rho\\ 0\end{smallmatrix}\right),\frac{1}{1-qr}\gamma\right)\right)\frac{\sqrt{\beta}}{\sqrt{1-\|\nu\|^2}} \alpha,\label{eqn_dimgen_Theta_rhorhoalpha}\\
		\boldsymbol\Theta_{\rho\alpha^2}&=\left(\vphantom{\left\|\frac{1}{\sqrt{1-qr}}\gamma\right\|^{-2}} \frac{2br^2}{3-r^2}\left\langle \frac{1}{\sqrt{1-qr}}B(r)\left(\begin{smallmatrix}\rho\\ 0\end{smallmatrix}\right),\frac{1}{1-qr}\gamma\right\rangle\right.\notag\\
		&\quad\left. + 3\left\|\frac{1}{\sqrt{1-qr}}\gamma\right\|^{-2}C\left(\frac{1}{\sqrt{1-qr}}B(r)\left(\begin{smallmatrix}\rho\\ 0\end{smallmatrix}\right),\frac{1}{1-qr}\gamma,\frac{1}{1-qr}\gamma\right)\right)\frac{\sqrt{\beta}}{1-\|\nu\|^2}\alpha^2,\label{eqn_dimgen_Theta_rhoalphaalpha}\\
		\boldsymbol\Theta_{\alpha^3}&=\left(-\frac{b^3r^6}{(3-r^2)^3}\left\|\frac{1}{\sqrt{1-qr}}\gamma\right\|^3 + \frac{br^2}{3-r^2}\left\|\frac{1}{\sqrt{1-qr}}\gamma\right\|^{-1}\left\|\frac{1}{1-qr}\gamma\right\|^2\right.\notag\\
		&\quad\left. + \left\|\frac{1}{\sqrt{1-qr}}\gamma\right\|^{-3} C\left(\frac{1}{1-qr}\gamma\right)\right)\frac{\sqrt{\beta}}{\sqrt{1-\|\nu\|^2}^3}\alpha^3\label{eqn_dimgen_Theta_alphaalphaalpha},
	\end{align}
for the pullback of $\boldsymbol\Theta_{v^2w}$ \eqref{eqn_dimgen_Theta_vvw}
	\begin{equation*}
		\left(L(r)\times(\mathbbm{1}_\mathbb{R}/\sqrt{\mu})\right)^*\boldsymbol\Theta_{v^2w}=\boldsymbol\Theta_{\rho^2\overline{w}}+\boldsymbol\Theta_{\rho\alpha\overline{w}}+\boldsymbol\Theta_{\alpha^2\overline{w}}
	\end{equation*}
with
	\begin{align}
		\boldsymbol\Theta_{\rho^2\overline{w}}&=\left(\frac{-r^2(2-3ar)+3-r^2}{r}\left\|\frac{1}{\sqrt{1-qr}}B(r)\left(\begin{smallmatrix}\rho\\ 0\end{smallmatrix}\right)\right\|^2 - \frac{3-r^2}{r}\langle\rho,\rho\rangle\right)\frac{\sqrt{\beta}}{\sqrt{\mu}}\overline{w},\label{eqn_dimgen_Theta_rhorhooverlinew}\\
		\boldsymbol\Theta_{\rho\alpha\overline{w}}&=\frac{2(-r^2(2-3ar)+3-r^2)}{r\sqrt{1-\|\nu\|^2}}\left\|\frac{1}{\sqrt{1-qr}}\gamma\right\|^{-1} \left\langle\frac{1}{\sqrt{1-qr}}B(r)\left(\begin{smallmatrix}\rho\\ 0\end{smallmatrix}\right),\frac{1}{1-qr}\gamma\right\rangle\frac{\sqrt{\beta}}{\sqrt{\mu}}\alpha\overline{w},\label{eqn_dimgen_Theta_rhoalphaoverlinew}\\
		\boldsymbol\Theta_{\alpha^2\overline{w}}&=\left(\frac{3b^2r^5(2-3ar)}{(3-r^2)^2}\left\|\frac{1}{\sqrt{1-qr}}\gamma\right\|^2-\frac{3-r^2}{r} \right.\notag\\
		&\quad\left. + \frac{-r^2(2-3ar)+3-r^2}{r}\left\|\frac{1}{\sqrt{1-qr}}\gamma\right\|^{-2}\left\|\frac{1}{1-qr}\gamma\right\|^2\right)\frac{\sqrt{\beta}}{\sqrt{\mu}(1-\|\nu\|^2)}\alpha^2\overline{w},\label{eqn_dimgen_Theta_alphaalphaoverlinew}
	\end{align}
and for the pullback of $\boldsymbol\Theta_{vw^2}$ \eqref{eqn_dimgen_Theta_vww}
	\begin{equation*}
		\left(L(r)\times(\mathbbm{1}_\mathbb{R}/\sqrt{\mu})\right)^*\boldsymbol\Theta_{vw^2}=\boldsymbol\Theta_{\rho\overline{w}^2}+\boldsymbol\Theta_{\alpha\overline{w}^2}
	\end{equation*}
with
	\begin{align}
		\boldsymbol\Theta_{\rho\overline{w}^2}&=0,\label{eqn_dimgen_Theta_rhooverlinewoverlinew}\\
		\boldsymbol\Theta_{\alpha\overline{w}^2}&=b\frac{-3r^4(2-3ar)^2+r^2(3-r^2)^2+(3-r^2)^3}{(3-r^2)}\left\|\frac{1}{\sqrt{1-qr}}\gamma\right\|\frac{\sqrt{\beta}}{\mu\sqrt{1-\|\nu\|^2}}\alpha\overline{w}^2.\label{eqn_dimgen_Theta_alphaoverlinewoverlinew}
	\end{align}
For $\boldsymbol\Theta_{w^3}$ \eqref{eqn_dimgen_Theta_www} we find that the same calculation as for equation \eqref{eqn_Theta_w3_limit_R_less_sqrt3} yields
	\begin{equation}
		\lim\limits_{r\nearrow R}\left(L(r)\times(\mathbbm{1}_\mathbb{R}/\sqrt{\mu})\right)^*\boldsymbol\Theta_{w^3}=-\frac{2}{3\sqrt{3}}\overline{w}^3.\label{eqn_Theta_www_limit_L_trafo_case}
	\end{equation}
We want to explicitly calculate the limits of equations \eqref{eqn_dimgen_Theta_rhorhorho}--\eqref{eqn_dimgen_Theta_alphaoverlinewoverlinew} as $r\to R$. However, in order to do that we will in some cases first need to analyse the implications of the terms $\boldsymbol\Theta_{\rho^3}$--$\boldsymbol\Theta_{\rho\overline{w}^2}$ in equations \eqref{eqn_dimgen_Theta_rhorhorho}--\eqref{eqn_dimgen_Theta_rhooverlinewoverlinew} being bounded as $r\to R$ in the sense that for all these terms there exists a homogeneous cubic polynomial in their respective variables, e.g. a polynomial $P$ in $\rho,\alpha$ for $\boldsymbol\Theta_{\rho^2\alpha}$, such that $\left|\boldsymbol\Theta_{\rho^2\alpha}\right|<|P|$ for all $r\in[0,R)$. This is a consequence of the assumption that $\mathcal{H}$ is closed in its ambient space, Theorem \ref{thm_Cn}, and the fact that we have shown that the limit of $|\boldsymbol\chi_{\alpha\overline{w}}|$ \eqref{eqn_chi_alpha_overlinew} is always bounded away from $2|\alpha\overline{w}|$ from below, cf. equation \eqref{eqn_chi_alpha_overlinew_limit_estimate}. Also recall that that both $\beta$ and $\sqrt{1-\|\nu\|^2}$ have simple zeros in $r=R$, cf. \eqref{eqn_partial_r_nusquare_limit}. Since $B$ is smooth with domain $[0,R]$ and $q_i<\frac{1}{R}$ for all $1\leq i\leq n-1$, we immediately see that
	\begin{equation}
		\lim\limits_{r\nearrow R}\boldsymbol\Theta_{\rho^3}=0.\label{eqn_Theta_rhorhorho_limit}
	\end{equation}
Before calculating the rest of the limits, we define the following abbreviations
	\begin{align}
		K&:=\left\|\frac{1}{\sqrt{1-qR}}\gamma\right\|^{-2}\left\|\frac{1}{1-qR}\gamma\right\|^2,\label{eqn_K_def}\\
		Y&:=9+R^2+(3-R^2)K.\label{eqn_Y_def}
	\end{align}
We can now write rewrite equations \eqref{eqn_partial_r_nusquare_limit} and \eqref{eqn_sqrt_beta_sqrt_1minusnusquare_limit} as 
	\begin{align}
		\left.\partial_r\|\nu\|^2\right|_{r=R}&=\frac{Y}{R(3-R^2)},\label{eqn_partial_r_nusquare_limit_alternative}\\
		\lim\limits_{r\nearrow R}\frac{\sqrt{\beta}}{\sqrt{1-\|\nu\|^2}}&=\frac{3-R^2}{\sqrt{Y}},\label{eqn_sqrt_beta_sqrt_1minusnusquare_limit_alternative}
	\end{align}
respectively. The limit of $\boldsymbol\Theta_{\rho^2\alpha}$ \eqref{eqn_dimgen_Theta_rhorhoalpha} can be calculated directly by using equation \eqref{eqn_sqrt_beta_sqrt_1minusnusquare_limit_alternative} and
	\begin{align}
		\frac{1}{\sqrt{1-qR}}B(R)\partial_{n-1}&=\left\|\frac{1}{\sqrt{1-qR}}\gamma\right\|^{-1}\frac{1}{1-qR}\gamma\notag\\
			&=\frac{\sqrt{3}|b|R^2}{3-R^2}\frac{1}{1-qR}\gamma,\label{eqn_B_at_R_partial_n_minus_1}
	\end{align}
which follows from equations \eqref{eqn_ker_0_1_condition} and \eqref{eqn_ker_0_1_B_gamma_condition_alternative}. We obtain
	\begin{align}
		\lim\limits_{r\nearrow R}\boldsymbol\Theta_{\rho^2\alpha}& = \;\mathrm{sgn}(b)\left(\frac{3-R^2}{\sqrt{3}\sqrt{Y}}\left\|\frac{1}{\sqrt{1-qR}}B(R)\left(\begin{smallmatrix}\rho\\ 0\end{smallmatrix}\right)\right\|^2\right.\notag\\
		&\quad\left. + \frac{3\sqrt{3}|b|R^2}{\sqrt{Y}}C\left(\frac{1}{\sqrt{1-qR}}B(R)\left(\begin{smallmatrix}\rho\\ 0\end{smallmatrix}\right),\frac{1}{\sqrt{1-qR}}B(R)\left(\begin{smallmatrix}\rho\\ 0\end{smallmatrix}\right),\frac{1}{1-qR}\gamma\right)\right) \alpha.\label{eqn_Theta_rhorhoalpha_limit}
	\end{align}
In later calculations, we will need equation \eqref{eqn_Theta_rhorhoalpha_limit}, and to reduce the amount of long formulas we define a bilinear form $W:\mathbb{R}^{n-2}\times\mathbb{R}^{n-2}\to\mathbb{R}$ fulfilling
	\begin{equation}
		\lim\limits_{r\nearrow R}\boldsymbol\Theta_{\rho^2\alpha} =\frac{W(\rho,\rho)}{\sqrt{Y}}\alpha.\label{eqn_W_def}
	\end{equation}
For $|\boldsymbol\Theta_{\rho\alpha^2}|$ \eqref{eqn_dimgen_Theta_rhoalphaalpha} to be bounded from above by the absolute value of some homogeneous cubic polynomial for all $r\in[0,R)$ we see, using the fact that $\beta$ and $1-\|\nu\|^2$ have simple zeros in $r=R$, that the right hand side of
	\begin{align}
		\frac{(1-\|\nu\|^2)\boldsymbol\Theta_{\rho\alpha^2}}{\sqrt{\beta}\alpha^2} &= \frac{2br^2}{3-r^2}\left\langle \frac{1}{\sqrt{1-qr}}B(r)\left(\begin{smallmatrix}\rho\\ 0\end{smallmatrix}\right),\frac{1}{1-qr}\gamma\right\rangle \notag\\
		&\quad  + 3\left\|\frac{1}{\sqrt{1-qr}}\gamma\right\|^{-2}C\left(\frac{1}{\sqrt{1-qr}}B(r)\left(\begin{smallmatrix}\rho\\ 0\end{smallmatrix}\right),\frac{1}{1-qr}\gamma,\frac{1}{1-qr}\gamma\right)\label{eqn_Theta_rhoalphaalpha_vanishing_condition}
	\end{align}
must vanish at $r=R$, which by equation \eqref{eqn_B_at_R_partial_n_minus_1} means that
	\begin{align}
		&C\left(\frac{1}{\sqrt{1-qR}}B(R)\left(\begin{smallmatrix}\rho\\ 0\end{smallmatrix}\right),\frac{1}{1-qR}\gamma,\frac{1}{1-qR}\gamma\right)\notag\\
		&=-\frac{2\;\mathrm{sgn}(b)}{3\sqrt{3}}\left\|\frac{1}{\sqrt{1-qR}}\gamma\right\|\left\langle \frac{1}{\sqrt{1-qR}}B(R)\left(\begin{smallmatrix}\rho\\ 0\end{smallmatrix}\right),\frac{1}{1-qR}\gamma\right\rangle\label{eqn_dC_vanishing_cond_ker_0_1}
	\end{align}
for all $\rho\in\mathbb{R}^{n-2}$. But since the right hand side of equation \eqref{eqn_Theta_rhoalphaalpha_vanishing_condition} is smooth on $[0,R]$, meaning that its derivative in $r$-direction in $r=R$ is bounded for all fixed $\rho\in\mathbb{R}^{n-2}$, we obtain together with the fact that $1-\|\nu\|^2$ has a simple zero in $r=R$, cf. equation \eqref{eqn_partial_r_nusquare_limit}, that
	\begin{equation}
		\lim\limits_{r\nearrow R}\boldsymbol\Theta_{\rho\alpha^2}=0.\label{eqn_Theta_rhoalphaalpha_limit}
	\end{equation}
Arguing similarly for $\boldsymbol\Theta_{\alpha^3}$ \eqref{eqn_dimgen_Theta_alphaalphaalpha} yields that the right hand side of
	\begin{align}
		\frac{\sqrt{1-\|\nu\|^2}^3\boldsymbol\Theta_{\alpha^3}}{\sqrt{\beta}\alpha^3}&=-\frac{b^3r^6}{(3-r^2)^3}\left\|\frac{1}{\sqrt{1-qr}}\gamma\right\|^3 + \frac{br^2}{3-r^2}\left\|\frac{1}{\sqrt{1-qr}}\gamma\right\|^{-1}\left\|\frac{1}{1-qr}\gamma\right\|^2\notag\\
					&\quad + \left\|\frac{1}{\sqrt{1-qr}}\gamma\right\|^{-3} C\left(\frac{1}{1-qr}\gamma\right)\label{eqn_Theta_aaa_vanishing_condition}
	\end{align}
must vanish at $r=R$. This is equivalent to
	\begin{equation}
		C\left(\frac{1}{1-qR}\gamma\right)=\;\mathrm{sgn}(b)\left(\frac{1}{3\sqrt{3}}-\frac{K}{\sqrt{3}}\right)\left\|\frac{1}{\sqrt{1-qR}}\gamma\right\|^3,\label{eqn_C_1_by_1_minus_qR_gamma}
	\end{equation}
where we have used equation \eqref{eqn_B_at_R_partial_n_minus_1}. In order to find the limit of $\boldsymbol\Theta_{\alpha^3}$ \eqref{eqn_dimgen_Theta_alphaalphaalpha} as $r\to R$, we first calculate
	\begin{align}
		\left.\partial_r\frac{1}{\sqrt{1-qr}}\gamma\right|_{r=R}&=\frac{q}{2\sqrt{1-qR}^3}\gamma,\label{eqn_r_der_sqrt_stuff}\\
		\left.\partial_r\left\|\frac{1}{\sqrt{1-qr}}\gamma\right\|\right|_{r=R}&= \frac{-1+K}{2R}\left\|\frac{1}{\sqrt{1-qR}}\gamma\right\|,\label{eqn_r_der_sqrt_stuff_norm}\\
		\left.\partial_r\frac{1}{1-qr}\gamma\right|_{r=R}&=\frac{q}{(1-qR)^2}\gamma,\label{eqn_r_der_stuff}\\
		\left.\partial_r\left\|\frac{1}{1-qr}\gamma\right\|\right|_{r=R}&= \frac{1}{R}\left(-\left\|\frac{1}{1-qR}\gamma\right\|+\left\|\frac{1}{1-qR}\gamma\right\|^{-1}\left\|\frac{1}{\sqrt{1-qR}^3}\gamma\right\|^2\right).\label{eqn_r_der_stuff_norm}
	\end{align}
Furthermore, observe that
	\begin{equation}
		\frac{1}{(1-qR)^2}\gamma=t\frac{1}{1-qR}\gamma + \frac{1}{\sqrt{1-qR}}B(R)\left(\begin{smallmatrix}Z\\ 0\end{smallmatrix}\right)\label{eqn_t_Z_split}
	\end{equation}
has a unique solution $t\in\mathbb{R}$, $Z\in\mathbb{R}^{n-2}$, which follows with equation \eqref{eqn_B_at_R_partial_n_minus_1}. We find by multiplying both hands of equation \eqref{eqn_t_Z_split} with $\sqrt{1-qR}$ from the left and then taking the Euclidean scalar product with $B(R)\partial_{n-1}$ that $t=K$. In order to calculate the limit of $\boldsymbol\Theta_{\alpha^3}$ \eqref{eqn_dimgen_Theta_alphaalphaalpha}, we will need the $r$-derivative of equation \eqref{eqn_Theta_aaa_vanishing_condition} at $r=R$. We obtain with $t=K$ in equation \eqref{eqn_t_Z_split} and the help of equations \eqref{eqn_dC_vanishing_cond_ker_0_1} and \eqref{eqn_C_1_by_1_minus_qR_gamma}--\eqref{eqn_t_Z_split}
	\begin{align}
		\left.\partial_r \frac{\sqrt{1-\|\nu\|^2}^3\boldsymbol\Theta_{\alpha^3}}{\sqrt{\beta}\alpha^3} \right|_{r=R}&= -\;\mathrm{sgn}(b)\frac{2\sqrt{3}(1-K)}{R(3-R^2)}. \label{eqn_Theta_alphaalphaalpha_pre_limit}
	\end{align}
Hence, with the help of L'H\^opital's rule equations \eqref{eqn_partial_r_nusquare_limit_alternative}, \eqref{eqn_sqrt_beta_sqrt_1minusnusquare_limit_alternative}, and \eqref{eqn_Theta_alphaalphaalpha_pre_limit} finally show that
	\begin{align}
		\lim\limits_{r\nearrow R}\boldsymbol\Theta_{\alpha^3}&=\;\mathrm{sgn}(b)\frac{2\sqrt{3}(1-K)(3-R^2)}{\sqrt{Y}^3}.\label{eqn_Theta_alphaalphaalpha_limit}
	\end{align}

Next we will determine the limits of $\boldsymbol\Theta_{\rho^2\overline{w}}$\;--\;$\boldsymbol\Theta_{\alpha^2\overline{w}}$ \eqref{eqn_dimgen_Theta_rhorhooverlinew}--\eqref{eqn_dimgen_Theta_alphaalphaoverlinew} as $r\to R$. First observe that it follows from $a=\frac{R^2-1}{R^3}$ that
	\begin{equation}
		\left.\left(r^2(2-3ar)-(3-r^2)\right)\right|_{r=R}=0,\label{eqn_r_a_stuff_at_R}
	\end{equation}
and we find
	\begin{equation}
		\left.\partial_r\left(r^2(2-3ar)-(3-r^2)\right)\right|_{r=R}=\frac{3(3-R^2)}{R}>0\label{eqn_r_derivative_r_a_stuff_at_R}
	\end{equation}
since $R<\sqrt{3}$. Using \eqref{eqn_sqrt_mu_by_sqrt_beta_limit} and \eqref{eqn_r_a_stuff_at_R}, we thus obtain for the limit of $\boldsymbol\Theta_{\rho^2\overline{w}}$ \eqref{eqn_dimgen_Theta_rhorhooverlinew}
	\begin{equation}
		\lim\limits_{r\nearrow R}\boldsymbol\Theta_{\rho^2\overline{w}}=-\frac{1}{\sqrt{3}}\langle\rho,\rho\rangle\overline{w}.\label{eqn_Theta_rhorhooverlinew_limit}
	\end{equation}
Taking additionally into account that $1-\|\nu\|^2$ has a simple zero in $r=R$, cf. equation \eqref{eqn_partial_r_nusquare_limit}, we obtain for the limit of $\boldsymbol\Theta_{\rho\alpha\overline{w}}$ \eqref{eqn_dimgen_Theta_rhoalphaoverlinew}
	\begin{equation}
		\lim\limits_{r\nearrow R}\boldsymbol\Theta_{\rho\alpha\overline{w}}=0.\label{eqn_Theta_rhoalphaoverlinew_limit}
	\end{equation}
To determine the limit of $\boldsymbol\Theta_{\alpha^2\overline{w}}$ \eqref{eqn_dimgen_Theta_alphaalphaoverlinew}, we first use equations \eqref{eqn_r_der_sqrt_stuff_norm}, \eqref{eqn_r_der_stuff_norm}, and \eqref{eqn_r_derivative_r_a_stuff_at_R} to obtain that
	\begin{align}
		\left.\frac{\sqrt{\mu}(1-\|\nu\|^2)\boldsymbol\Theta_{\alpha^2\overline{w}}}{\sqrt{\beta}\alpha^2\overline{w}}\right|_{r=R} &= 0,\label{eqn_alphaalphaoverlinew_ker_0_1_limit_vanishing}\\
		\left.\partial_r\frac{\sqrt{\mu}(1-\|\nu\|^2)\boldsymbol\Theta_{\alpha^2\overline{w}}}{\sqrt{\beta}\alpha^2\overline{w}}\right|_{r=R} &= \frac{2(9-R^2-K(3-R^2))}{R^2},\label{eqn_alphaalphaoverlinew_ker_0_1_limit_r_der_vanishing}
	\end{align}
Hence, using L'H\^opital's rule together with equations \eqref{eqn_alphaalphaoverlinew_ker_0_1_limit_vanishing}, \eqref{eqn_alphaalphaoverlinew_ker_0_1_limit_r_der_vanishing}, \eqref{eqn_sqrt_mu_by_sqrt_beta_limit}, and \eqref{eqn_partial_r_nusquare_limit_alternative} implies
	\begin{equation}
		\lim\limits_{r\nearrow R}\boldsymbol\Theta_{\alpha^2\overline{w}} = \left(\frac{2}{\sqrt{3}}-\frac{12\sqrt{3}}{Y}\right)\alpha^2\overline{w}.\label{eqn_Theta_alphaalphaoverlinew_limit}
	\end{equation}
Since $\boldsymbol\Theta_{\rho\overline{w}^2}$ \eqref{eqn_dimgen_Theta_rhooverlinewoverlinew} vanishes identically for all $r\in[0,R)$, it follows that
	\begin{equation}
		\lim\limits_{r\nearrow R}\boldsymbol\Theta_{\rho\overline{w}^2}=0.\label{eqn_Theta_rhooverlinewoverlinew_limit}
	\end{equation}
For the limit of $\boldsymbol\Theta_{\alpha\overline{w}^2}$ \eqref{eqn_dimgen_Theta_alphaoverlinewoverlinew}, observe that equation \eqref{eqn_r_a_stuff_at_R} implies that the term $-r^4(2-3ar)^2+(3-r^2)^2$ in the numerator of one of the factors of $\boldsymbol\Theta_{\alpha\overline{w}^2}$ also vanishes in $r=R$. We further calculate
	\begin{equation*}
		\frac{-r^4(2-3ar)^2+(3-r^2)^2}{\mu}= \frac{(3+r^2)R^3-3r^3(1-R^2)}{(3+r^2)R^3+9r(1-R^2)},
	\end{equation*}
which implies
	\begin{equation}
		\lim\limits_{r\nearrow R}\frac{-r^4(2-3ar)^2+(3-r^2)^2}{\mu}=\frac{2R^2}{3-R^2}.\label{eqn_some_stuff_by_mu_limit}
	\end{equation}
It now follows with equation \eqref{eqn_sqrt_beta_sqrt_1minusnusquare_limit_alternative} that
	\begin{equation}
		\lim\limits_{r\nearrow R}\boldsymbol\Theta_{\alpha\overline{w}^2} = \;\mathrm{sgn}(b)\frac{2\sqrt{3}}{\sqrt{Y}}\alpha\overline{w}^2.\label{eqn_Theta_alphaoverlinewoverlinew_limit}
	\end{equation}
Summarizing up to this point, equations \eqref{eqn_chi_alpha_overlinew_limit}, \eqref{eqn_Theta_rhorhorho_limit}, \eqref{eqn_Theta_rhorhoalpha_limit}, \eqref{eqn_Theta_rhoalphaalpha_limit}, \eqref{eqn_Theta_alphaalphaalpha_limit}, \eqref{eqn_Theta_rhorhooverlinew_limit}, \eqref{eqn_Theta_rhoalphaoverlinew_limit}, \eqref{eqn_Theta_alphaalphaoverlinew_limit}, \eqref{eqn_Theta_rhooverlinewoverlinew_limit}, \eqref{eqn_Theta_alphaoverlinewoverlinew_limit}, and \eqref{eqn_Theta_www_limit_L_trafo_case} imply that the limit of our transformed initial polynomial $h_E$ \eqref{eqn_dimgen_hE1_prelimit} under the family of transformations $\mathbbm{1}_\mathbb{R}\times L(r)\times (\mathbbm{1}_\mathbb{R}/\sqrt{\mu})$ \eqref{eqn_L_trafo} is given by
	\begin{align}
		&\lim\limits_{r\nearrow R} \left(\mathbbm{1}_\mathbb{R}\times L(r)\times (\mathbbm{1}_\mathbb{R}/\sqrt{\mu})\right)^* h_{E}\left(\left(\begin{smallmatrix} x\\ \rho\\ \alpha\\ \overline{w}\end{smallmatrix}\right)\right)\notag\\
		&=x^3-x\left(\langle \rho,\rho\rangle + \alpha^2 + \boldsymbol\zeta  + \overline{w}^2\right) + \boldsymbol\Xi_{\rho^2\alpha} + \boldsymbol\Xi_{\alpha^3} + \boldsymbol\Xi_{\rho^2\overline{w}} + \boldsymbol\Xi_{\alpha^2\overline{w}} + \boldsymbol\Xi_{\alpha\overline{w}^2} + \boldsymbol\Xi_{\overline{w}^3},\label{eqn_limit_ker_1_case_pre_trafo}
	\end{align}
where
	\begin{equation*}
		\begin{array}{|c||c|c|c|c|c|c|c|}
			\hline \text{symbol} & \overset{}{\boldsymbol\zeta} & \boldsymbol\Xi_{\rho^2\alpha} & \boldsymbol\Xi_{\alpha^3} & \boldsymbol\Xi_{\rho^2\overline{w}} & \boldsymbol\Xi_{\alpha^2\overline{w}} & \boldsymbol\Xi_{\alpha\overline{w}^2} & \boldsymbol\Xi_{\overline{w}^3}\\
			\hline \text{meaning} & \lim\limits_{r\nearrow R} \boldsymbol\chi_{\alpha\overline{w}} & \lim\limits_{r\nearrow R}\boldsymbol\Theta_{\rho^2\alpha} & \lim\limits_{r\nearrow R}\boldsymbol\Theta_{\alpha^3} & -\frac{1}{\sqrt{3}}\langle\rho,\rho\rangle\overline{w} & \lim\limits_{r\nearrow R}\boldsymbol\Theta_{\alpha^2\overline{w}} & \lim\limits_{r\nearrow R}\boldsymbol\Theta_{\alpha\overline{w}^2} & \overset{}{-\frac{2}{3\sqrt{3}}\overline{w}^3}\\
			\hline \text{cf. eqn.} & \overset{}{\text{\eqref{eqn_chi_alpha_overlinew_limit}}} & \text{\eqref{eqn_Theta_rhorhoalpha_limit}} & \text{\eqref{eqn_Theta_alphaalphaalpha_limit}} & \text{\eqref{eqn_Theta_rhorhooverlinew_limit}} & \text{\eqref{eqn_Theta_alphaalphaoverlinew_limit}} & \text{\eqref{eqn_Theta_alphaoverlinewoverlinew_limit}} & \text{\eqref{eqn_Theta_www_limit_L_trafo_case}}\\
			\hline
		\end{array}
	\end{equation*}
The next steps are similar to the $2$-dimensional case, cf. Section \eqref{subsect_dim2_b_ne_0}, equation \eqref{eqn_dim2_hE2_limit} onward. We need to transform the $\alpha$- and $\overline{w}$-coordinate, so that the limit of the transformed polynomial $h_E$ in \eqref{eqn_limit_ker_1_case_pre_trafo} is brought to standard form \eqref{eqn_h_general_form}. For the following calculations, define
	\begin{equation}
		\zeta:=-\;\mathrm{sgn}(b)\frac{6}{\sqrt{Y}}.\label{eqn_zeta_ker0_1_def}
	\end{equation}
Using the above equation \eqref{eqn_zeta_ker0_1_def}, equations \eqref{eqn_K_def}, \eqref{eqn_Y_def}, and \eqref{eqn_W_def}, $\boldsymbol\zeta$, $\boldsymbol\Xi_{\rho^2\alpha}$, $\boldsymbol\Xi_{\alpha^3}$, $\boldsymbol\Xi_{\alpha^2\overline{w}}$, and $\boldsymbol\Xi_{\alpha\overline{w}^2}$ can be written as
	\begin{align}
		\boldsymbol\zeta &= -\;\mathrm{sgn}(b)\frac{6}{\sqrt{Y}} \alpha \overline{w} =\zeta\alpha\overline{w},\label{eqn_zeta_ker_0_1}\\
		\boldsymbol\Xi_{\rho^2\alpha} &= \frac{W(\rho,\rho)}{\sqrt{Y}}\alpha,\label{eqn_Xi_rhorhoalpha}\\
		\boldsymbol\Xi_{\alpha^3} &= \;\mathrm{sgn}(b)\frac{2\sqrt{3}(1-K)(3-R^2)}{\sqrt{Y}^3} =\frac{1}{3\sqrt{3}}\zeta(3-\zeta^2)\alpha^3,\label{eqn_Xi_alphaalphaalpha}\\
		\boldsymbol\Xi_{\alpha^2\overline{w}} &= \left(\frac{2}{\sqrt{3}}-\frac{12\sqrt{3}}{Y}\right)\alpha^2\overline{w} =\frac{1}{\sqrt{3}}(2-\zeta^2)\alpha^2\overline{w},\\
		\boldsymbol\Xi_{\alpha\overline{w}^2} &= \;\mathrm{sgn}(b) \frac{2\sqrt{3}}{\sqrt{Y}}\alpha\overline{w}^2 =-\frac{1}{\sqrt{3}}\zeta\alpha\overline{w}^2.\label{eqn_Xi_alphaoverlinewoverlinew}
	\end{align}
Now as in the $2$-dimensional case
we set
	\begin{equation}
		T:=1-\frac{\zeta}{2},\quad t:=1+\frac{\zeta}{2},\label{eqn_T_t_dimgen}
	\end{equation}
and transform $\alpha$ and $\overline{w}$ via $M\in\mathrm{GL}(2)$,
	\begin{equation}
		\left(\begin{matrix}\alpha\\ \overline{w}\end{matrix}\right)=\frac{1}{\sqrt{2}}\left(\begin{matrix} \frac{1}{\sqrt{t}} & \frac{1}{\sqrt{T}}\\ \frac{1}{\sqrt{t}} & -\frac{1}{\sqrt{T}}\end{matrix}\right)\left(\begin{matrix}k\\ \ell\end{matrix}\right)=:M \left(\begin{matrix}k\\ \ell\end{matrix}\right).\label{eqn_M_def}
	\end{equation}
Note that $M$ is well defined due to \eqref{eqn_chi_alpha_overlinew_limit_estimate}. We verify that
	\begin{equation}
		M^*(\alpha^2+\boldsymbol\zeta+\overline{w}^2)=k^2+\ell^2\label{eqn_check_M_bilin_pullback}
	\end{equation}
and further calculate
	\begin{align}
		&\left(\mathbbm{1}_{\mathbb{R}^{n-1}}\times M\right)^*\lim\limits_{r\nearrow R} \left(\mathbbm{1}_\mathbb{R}\times L(r)\times (\mathbbm{1}_\mathbb{R}/\sqrt{\mu})\right)^* h_{E}\left(\left(\begin{smallmatrix} x\\ \rho\\ k\\ \ell\end{smallmatrix}\right)\right)\notag\\
		&=x^3-x\left(\langle \rho,\rho\rangle + k^2 + \ell^2\right) + \boldsymbol\Xi_{\rho^2k} + \boldsymbol\Xi_{\rho^2\ell} + \boldsymbol\Xi_{k^3} + \boldsymbol\Xi_{k^2\ell} + \boldsymbol\Xi_{k\ell^2} + \boldsymbol\Xi_{\ell^3}\label{eqn_limit_ker_1_case_post_trafo1}
	\end{align}
with
	\begin{align}
		\boldsymbol\Xi_{\rho^2k}&=\left(\frac{\boldsymbol\Xi_{\rho^2\alpha}}{\alpha} + \frac{\boldsymbol\Xi_{\rho^2\overline{w}}}{\overline{w}}\right)\frac{k}{\sqrt{2t}}\notag\\
		&=\left(\frac{W(\rho,\rho)}{\sqrt{Y}}-\frac{1}{\sqrt{3}}\langle\rho,\rho\rangle\right)\frac{k}{\sqrt{2t}},\label{eqn_Xi_rhorhok}\\
		\boldsymbol\Xi_{\rho^2\ell}&=\left(\frac{\boldsymbol\Xi_{\rho^2\alpha}}{\alpha} - \frac{\boldsymbol\Xi_{\rho^2\overline{w}}}{\overline{w}}\right)\frac{\ell}{\sqrt{2T}}\notag\\
		&=\left(\frac{W(\rho,\rho)}{\sqrt{Y}}+\frac{1}{\sqrt{3}}\langle\rho,\rho\rangle\right)\frac{\ell}{\sqrt{2T}},\label{eqn_Xi_rhorhoell}\\
		\boldsymbol\Xi_{k^3}&=\left(\frac{\boldsymbol\Xi_{\alpha^3}}{\alpha^3} + \frac{\boldsymbol\Xi_{\alpha^2\overline{w}}}{\alpha^2\overline{w}} + \frac{\boldsymbol\Xi_{\alpha\overline{w}^2}}{\alpha\overline{w}^2} + \frac{\boldsymbol\Xi_{\overline{w}^3}}{\overline{w}^3}\right)\frac{k^3}{\sqrt{2t}^3}\notag\\
		&=\frac{2}{3\sqrt{3}}\cdot\frac{(3T-t)\sqrt{t}}{2\sqrt{2}}k^3,\label{eqn_Xi_kkk}\\
		\boldsymbol\Xi_{k^2\ell}&=\left(3\frac{\boldsymbol\Xi_{\alpha^3}}{\alpha^3} + \frac{\boldsymbol\Xi_{\alpha^2\overline{w}}}{\alpha^2\overline{w}} - \frac{\boldsymbol\Xi_{\alpha\overline{w}^2}}{\alpha\overline{w}^2} -3 \frac{\boldsymbol\Xi_{\overline{w}^3}}{\overline{w}^3}\right)\frac{k^2\ell}{2t\sqrt{2T}}\notag\\
		&=\frac{2}{\sqrt{3}}\cdot\frac{(-T+3t)\sqrt{s}}{2\sqrt{2}}k^2\ell,\label{eqn_Xi_kkell}\\
		\boldsymbol\Xi_{k\ell^2}&=\left(3\frac{\boldsymbol\Xi_{\alpha^3}}{\alpha^3} - \frac{\boldsymbol\Xi_{\alpha^2\overline{w}}}{\alpha^2\overline{w}} - \frac{\boldsymbol\Xi_{\alpha\overline{w}^2}}{\alpha\overline{w}^2} +3 \frac{\boldsymbol\Xi_{\overline{w}^3}}{\overline{w}^3}\right)\frac{k\ell^2}{2T\sqrt{2t}}\notag\\
		&=\frac{2}{\sqrt{3}}\cdot\frac{(-3T+t)\sqrt{t}}{2\sqrt{2}}k\ell^2,\label{eqn_Xi_kellell}\\
		\boldsymbol\Xi_{\ell^3}&=\left(\frac{\boldsymbol\Xi_{\alpha^3}}{\alpha^3} - \frac{\boldsymbol\Xi_{\alpha^2\overline{w}}}{\alpha^2\overline{w}} + \frac{\boldsymbol\Xi_{\alpha\overline{w}^2}}{\alpha\overline{w}^2} - \frac{\boldsymbol\Xi_{\overline{w}^3}}{\overline{w}^3}\right)\frac{\ell^3}{\sqrt{2T}^3}\notag\\
		&=\frac{2}{3\sqrt{3}}\cdot\frac{(T-3t)\sqrt{T}}{2\sqrt{2}}\ell^3.\label{eqn_Xi_ellellell}
	\end{align}
We need one more transformation in the $k$-$\ell$-coordinates, given by $N\in\mathrm(O)(2)$,
	\begin{equation}
		\left(\begin{matrix}
			k\\ \ell\end{matrix}\right)
			=\frac{1}{\sqrt{2}}\left(\begin{matrix} \sqrt{t} & \sqrt{T}\\ -\sqrt{T} & \sqrt{t}\end{matrix}\right)\left(\begin{matrix}\widetilde{k}\\ \widetilde{\ell}\end{matrix}\right)=:N \left(\begin{matrix}\widetilde{k}\\ \widetilde{\ell}
		\end{matrix}\right).\label{eqn_N_def}
	\end{equation}
We obtain with equations \eqref{eqn_limit_ker_1_case_post_trafo1}--\eqref{eqn_N_def}
	\begin{align}
		&\left(\mathbbm{1}_{\mathbb{R}^{n-1}}\times N\right)^*\left(x^3-x\left(\langle \rho,\rho\rangle + k^2 + \ell^2\right) + \boldsymbol\Xi_{\rho^2k} + \boldsymbol\Xi_{\rho^2\ell} + \boldsymbol\Xi_{k^3} + \boldsymbol\Xi_{k^2\ell} + \boldsymbol\Xi_{k\ell^2} + \boldsymbol\Xi_{\ell^3}\right)\notag\\
		&=x^3-x\left(\langle \rho,\rho\rangle + \widetilde{k}^2 + \widetilde{\ell}^2\right)\notag\\
		&\quad + \frac{1}{\sqrt{Tt}}\left(\frac{W(\rho,\rho)}{\sqrt{Y}}+\frac{\zeta}{2\sqrt{3}}\langle \rho,\rho\rangle\right)\widetilde{\ell} -\frac{2}{3\sqrt{3}}\widetilde{k}^3-\frac{1}{\sqrt{3}}\langle\rho,\rho\rangle\widetilde{k} +\frac{2}{\sqrt{3}}\widetilde{k}\widetilde{\ell}^2.\label{eqn_limit_poly_ker_0_1}
	\end{align}
The above formula \eqref{eqn_limit_poly_ker_0_1} means that all possible limit polynomials for $\dim\ker\boldsymbol\lambda|_{r=R}=1$, $q_i<\frac{1}{R}$ for all $1\leq i\leq n-1$, are of the required form as in Theorem \ref{thm_limit_behaviour_ccpsr}. Furthermore we have obtained a condition on the $C$-term in the initial polynomial \eqref{eqn_h_n_geq_3_first_form}, namely that
	\begin{equation}
		-\langle\rho,\rho\rangle\leq \frac{1}{\sqrt{Tt}}\left(\frac{W(\rho,\rho)}{\sqrt{Y}}+\frac{\zeta}{2\sqrt{3}}\langle \rho,\rho\rangle\right) \leq \langle \rho,\rho\rangle.\label{eqn_W_condition_ker_0_1}
	\end{equation}
Note that this is not an open condition on $C\in\mathrm{Sym}^3\left(\mathbb{R}^{n-1}\right)^*$, but rather a condition on the restriction to an $(n-2)$-dimensional subspace, cf. equations \eqref{eqn_Theta_rhoalphaalpha_vanishing_condition} and \eqref{eqn_W_def}. To see that \eqref{eqn_W_condition_ker_0_1} must in fact hold see Section \ref{sect_Fi_calcs}. At this point we do not know if an example of our presently studied case $\dim\ker \boldsymbol\lambda|_{r=R}=1$, $\frac{1}{R}>q_1\geq\ldots\geq q_{n-1}$ exists in all dimensions $n\geq 3$. To give such an example we must find $P_3\left(\left(\begin{smallmatrix}v\\w\end{smallmatrix}\right)\right)$ as in \eqref{eqn_P3_term_bnon0}, such that the latter conditions are fulfilled and the maximality condition \eqref{eqn_P3_max_condition} holds.

\begin{Ex}\label{ex_ker01_n_geq_3}
	Let $P_3\left(\left(\begin{smallmatrix}v\\w\end{smallmatrix}\right)\right)=-\frac{2}{3\sqrt{3}}v_{n-1}^3+\frac{2}{\sqrt{3}}v_{n-1}w^2$, corresponding to
		\begin{equation*}
			R=1,\quad b=\frac{2}{\sqrt{3}},\quad q=0,\quad \gamma=\partial_{n-1},\quad C(v)=-\frac{2}{3\sqrt{3}}v_{n-1}^3.
		\end{equation*}
	Then $\dim\ker \boldsymbol\lambda|_{r=R}=1$ is true since all $q_i$ vanish and \eqref{eqn_ker_0_1_condition} is also easily seen to be true. In order to show that $\max\limits_{\{\langle v,v\rangle + w^2=1\}} P_3\left(\left(\begin{smallmatrix}v\\w\end{smallmatrix}\right)\right)\leq \frac{2}{3\sqrt{3}}$ we restrict $P_3\left(\left(\begin{smallmatrix}v\\w\end{smallmatrix}\right)\right)$ to planes of the form $\mathrm{span}\{V,\partial_w\}$, where $V\ne0$ is of the form $V=\sum\limits_{i=1}^{n-1} V^i \partial_{i}$. It is then immediate that $P_3\left(\left(\begin{smallmatrix}v\\w\end{smallmatrix}\right)\right)$ fulfilling \eqref{eqn_P3_max_condition} is equivalent to showing that
		\begin{equation*}
			\widetilde{P_3}\left(\left(\begin{smallmatrix}v\\w\end{smallmatrix}\right)\right):= -\frac{2}{3\sqrt{3}}c^3v+\frac{2}{\sqrt{3}}cvw^2,
		\end{equation*}
	$(v,w)^\mathrm{T}$ being linear coordinates on $\mathbb{R}^2$, fulfils $\max\limits_{\{v^2+w^2=1\}}\widetilde{P_3}\left(\left(\begin{smallmatrix}v\\w\end{smallmatrix}\right)\right)\leq\frac{2}{3\sqrt{3}}$ for all $c\in [0,1]$. For $c=0$ this is clear, and for $c\in(0,1]$ we find that $\widetilde{P_3}\left(\left(\begin{smallmatrix}v\\w\end{smallmatrix}\right)\right)$ attains the value $\frac{2}{3\sqrt{3}}$ at precisely the points
		\begin{equation*}
			\left(\begin{matrix}v\\w\end{matrix}\right)\in\left\{\left(\begin{matrix}\frac{1}{2c}\\ \pm\frac{\sqrt{c^2+2}}{2c}\end{matrix}\right), \left(\begin{matrix}-\frac{1}{c^3}\\0\end{matrix}\right)\right\}.
		\end{equation*}
	Since all of the above vectors have Euclidean norm at least $1$ for all $c\in(0,1]$, which follows from $\frac{c^2+3}{4c^2}\geq 1$ for all $c\in(0,1]$ which in turn follows easily by monotonicity, we obtain that $\widetilde{P_3}\left(\left(\begin{smallmatrix}v\\w\end{smallmatrix}\right)\right)$ does in fact fulfil the maximality condition \eqref{eqn_P3_max_condition} as required.
\end{Ex}

Next, we will consider the cases with $q_i=\frac{1}{R}$ for at least one $1\leq i\leq n-1$. Assume that $q_1\geq\ldots\geq q_{n-1}$ and that for $m\in\{1,\dots,n-2\}$ fixed, $q_m=\frac{1}{R}>q_{m+1}$. Note at this point that $q_1=\ldots=q_{n-1}=\frac{1}{R}$ and $b\ne 0$ implies that near $r=R$, $\boldsymbol{\lambda}$ would necessarily have a negative eigenvalue, violating the condition that $\mathcal{H}$ consists only of hyperbolic points of its defining polynomial. Hence we can exclude this case. Similarly by restriction to $\mathbb{R}^{m}\cong\left\{\left.(v_1,\ldots,v_m,0,\ldots,0)^\mathrm{T}\ \right|\ (v_1,\ldots,v_m)^\mathrm{T}\in\mathbb{R}^{m}\right\}\subset\mathbb{R}^{n-1}$ we see that $\gamma$ must necessarily be of the form
	\begin{equation}
		\gamma=\left(\begin{matrix} 0\\ \widehat{\gamma}\end{matrix}\right),\quad \widehat{\gamma}\in\mathbb{R}^{n-1-m},\quad \langle\widehat{\gamma},\widehat{\gamma}\rangle=1.\label{eqn_gamma_ker_m_and_m_1_cases}
	\end{equation}
In the following we will use the notation
	\begin{equation}\label{eqn_s_u_def}
		s=\left(\begin{matrix}s_1\\ \vdots\\ s_m\end{matrix}\right)=\left(\begin{matrix}v_1\\ \vdots\\ v_m\end{matrix}\right),\quad u=\left(\begin{matrix}u_1\\ \vdots\\ u_{n-1-m}\end{matrix}\right)=\left(\begin{matrix}v_{m+1}\\ \vdots\\ v_{n-1}\end{matrix}\right),
	\end{equation}
so that $v=\left(\begin{smallmatrix}s\\ u\end{smallmatrix}\right)$. We further denote
	\begin{equation*}
		q=\left(\begin{array}{c|c} \frac{1}{R}\mathbbm{1}_m & \\ \hline & \widehat{q}\end{array}\right),
	\end{equation*}
and obtain that $h$ \eqref{eqn_h_n_geq_3_first_form} is of the form
	\begin{align}
		h&=x^3-x\left(\langle s,s\rangle + \langle u,u\rangle + w^2\right)\notag\\
		&\quad+C(s)+3C(s,s,u)+3C(s,u,u)+C(u)\notag\\
		&\quad+\left(\frac{1}{R}\langle s,s\rangle + \langle u,\widehat{q}u\rangle \right)w + b\langle\widehat{\gamma},u\rangle w^2 + \frac{R^2-1}{R^3}w^3.\label{eqn_h_n_geq_3_first_form_m_case}
	\end{align}
Next, we rewrite $\boldsymbol\Theta_{v^3}$ \eqref{eqn_dimgen_Theta_vvv} as $\boldsymbol\Theta_{v^3}=\boldsymbol\Theta_{s^3}+\boldsymbol\Theta_{s^2u}+\boldsymbol\Theta_{su^2}+\boldsymbol\Theta_{u^3}$ with
	\begin{align}
		\boldsymbol\Theta_{s^3}&=(3-r^2)^3 C(s)\sqrt{\beta},\label{eqn_dimgen_Theta_sss}\\
		\boldsymbol\Theta_{s^2u}&=\left(br^2(3-r^2)^2\langle s,s\rangle\langle \widehat{\gamma},u\rangle + 3(3-r^2)^3 C(s,s,u)\right)\sqrt{\beta},\label{eqn_dimgen_Theta_ssu}\\
		\boldsymbol\Theta_{su^2}&=3(3-r^2)^3 C(s,u,u)\sqrt{\beta},\label{eqn_dimgen_Theta_suu}\\
		\boldsymbol\Theta_{u^3}&=\left(-b^3r^6\langle \widehat{\gamma},u\rangle ^3 + br^2(3-r^2)^2\langle u,u\rangle\langle \widehat{\gamma},u\rangle + (3-r^2)^3 C(u)\right)\sqrt{\beta}.\label{eqn_dimgen_Theta_uuu}
	\end{align}
We further write $\boldsymbol\lambda=\boldsymbol\lambda_{s^2}+\boldsymbol\lambda_{u^2}$ with
	\begin{align}
		\boldsymbol\lambda_{s^2}&:=(3-r^2)^2\left(1-\frac{r}{R}\right)\langle s,s\rangle,\label{eqn_lambda_ss}\\
		\boldsymbol\lambda_{u^2}&:=-3b^2r^4\langle \widehat{\gamma},u\rangle^2 + (3-r^2)^2\langle u,(1-\widehat{q}r)u\rangle,\label{eqn_lambda_uu}
	\end{align}
which makes it easy to see that either $\dim\ker \boldsymbol\lambda|_{r=R}=m$ or $\dim\ker \boldsymbol\lambda|_{r=R}=m+1$ holds. By restriction to the linear subspace $\mathbb{R}^{m}\cong \{u=0\}\subset\mathbb{R}^{n-1}$ and using the assumption $R<\sqrt{3}$ it thus follows that $C(s)\equiv 0$. To see the latter, observe that $\boldsymbol\lambda_{s^2}$ vanishes identically as $r\to R$. Hence, after rescaling the $s$-part of $v$ appropriately with the factor $(3-r^2)^{-1}\left(1-\frac{r}{R}\right)^{-1/2}$ transforming $\boldsymbol\lambda_{s^2}$ to $\langle s,s\rangle$, $\boldsymbol\Theta_{s^3}$ \eqref{eqn_dimgen_Theta_sss} transforms to
	\begin{equation}\label{C_s_pullback_b_ne_0_R_not_sqrt2}
		(3-r^2)\left(1-\frac{r}{R}\right)^{-\frac{3}{2}} C(s)\sqrt{\beta}.
	\end{equation}
The above term would have a maximum on $\{\langle s,s\rangle=1\}$ going to infinity as $r\to R$ if $C(s)\not\equiv 0$, where we recall that $R<\sqrt{3}$ implies that $\beta$ has a simple zero in $r=R$, thereby violating $\mathcal{H}$ being closed in the ambient space by Theorem \ref{thm_Cn}. With a similar argument for $\boldsymbol\Theta_{s^2u}$ \eqref{eqn_dimgen_Theta_ssu}, using again that $\beta$ and $\left(1-\frac{r}{R}\right)$ each have simple zeros in $r=R$, it follows that
	\begin{equation}\label{C_ssu_b_ne0_ker_m_mplus1}
		C(s,s,u)=-\frac{bR^2}{3(3-R^2)}\langle s,s\rangle \langle \widehat{\gamma},u\rangle
	\end{equation}
must hold, since otherwise $\mathcal{H}\subset\mathbb{R}^{n+1}$ being closed would be violated at this point. Note that the above formula must hold in both possible scenarios $\dim\ker \boldsymbol\lambda|_{r=R}=m$ and $\dim\ker \boldsymbol\lambda|_{r=R}=m+1$, the latter being the case if and only if $\dim\ker \boldsymbol\lambda|_{r=R}|_{\{s=0\}}=1$. Summarizing, we have shown that $h$ \eqref{eqn_h_n_geq_3_first_form_m_case} must be of the form
	\begin{align}
		h&=x^3-x\left(\langle s,s\rangle + \langle u,u\rangle + w^2\right) + P_3\left(\left(\begin{smallmatrix}s\\ u\\ w\end{smallmatrix}\right)\right)\label{eqn_h_n_geq_3_at_least_m_ker_lambda}
	\end{align}
with
	\begin{align}
		P_3\left(\left(\begin{smallmatrix}s\\ u\\ w\end{smallmatrix}\right)\right)&=-\frac{bR^2}{3-R^2}\langle s,s\rangle\langle\widehat{\gamma},u\rangle+3C(s,u,u)+C(u)\notag\\
		&\quad+\left(\frac{1}{R}\langle s,s\rangle + \langle u,\widehat{q}u\rangle \right)w + b\langle\widehat{\gamma},u\rangle w^2 + \frac{R^2-1}{R^3}w^3.\label{eqn_P_3_n_geq_3_at_least_m_ker_lambda}
	\end{align}
Observe that any point $p\in\mathbb{R}^{n}$ of the form
	\begin{equation}
		p=\left(\begin{matrix} s(p)\\ u(p) \\ w(p)\end{matrix}\right)=\left(\begin{matrix} s(p)\\ 0 \\ \frac{R}{\sqrt{3}}\end{matrix}\right),\quad \|s(p)\|=\frac{\sqrt{3-R^2}}{\sqrt{3}},\label{eqn_ker_m_m_plus_1_explicit_singatinfpoints}
	\end{equation}
fulfils $\|p\|=1$ and $P_3(p)=\frac{2}{3\sqrt{3}}$. This means that $P_3$ in \eqref{eqn_P_3_n_geq_3_at_least_m_ker_lambda} in particular is singular at infinity, cf. Definition \ref{def_singular_at_infinity}. Recall that at this point that both $\dim\ker \boldsymbol\lambda|_{r=R}=m$ or $\dim\ker \boldsymbol\lambda|_{r=R}=m+1$ are allowed.

We will deal with both cases separately and start with $\dim\ker \boldsymbol\lambda|_{r=R}=m$. Fix $\widehat{L}\in\mathrm{GL}(n-1-m)$, such that $\left(\widehat{L}^*\boldsymbol\lambda_{u^2}|_{r=R}\right)(u,u)=\langle u,u\rangle$ and let
	\begin{equation}
		L:[0,R)\to\mathrm{GL}(n+1),\quad L(r)=
			\left(\begin{array}{c|c} (3-r^2)^{-1}\left(1-\frac{r}{R}\right)^{-\frac{1}{2}}\mathbbm{1}_m &  \\ \hline & \overset{}{\widehat{L}} \end{array}\right).\label{eqn_B_r_ker_m_0}
	\end{equation}
Then
	\begin{equation*}
		\lim\limits_{r\nearrow R}\left(L(r)^*\boldsymbol\lambda\right)\left(\left(\begin{smallmatrix} s\\ u\end{smallmatrix}\right),\left(\begin{smallmatrix} s\\ u\end{smallmatrix}\right)\right)=\langle s,s\rangle + \langle u,u\rangle.
	\end{equation*}
Recall that $\mu$ has a simple zero in $r=R$, which follows from equation \eqref{eqn_dim2_sqrtmu_by_sqrtbeta_limit}. Together with $\widehat{B}\in\mathrm{GL}(n-1-m)$ not depending on $r$ and $\gamma$ being of the form $\gamma=\left(\begin{smallmatrix} 0\\ \widehat{\gamma}\end{smallmatrix}\right)$, $\widehat{\gamma}\in\mathbb{R}^{n-1-m}$, this implies that
	\begin{equation*}
		\lim\limits_{r\nearrow R} L(r)^*\boldsymbol\chi =0.
	\end{equation*}
Now we need to study the limit of the pullback of the different $\boldsymbol\Theta$-terms \eqref{eqn_dimgen_Theta_vvv}--\eqref{eqn_dimgen_Theta_www} with respect to $L(r)\times (\mathbbm{1}_\mathbb{R}/\sqrt{\mu})$ as $r\to R$. We first consider $\boldsymbol\Theta_{v^3}=\boldsymbol\Theta_{s^3}+\boldsymbol\Theta_{s^2u}+\boldsymbol\Theta_{su^2}+\boldsymbol\Theta_{u^3}$, cf. equations \eqref{eqn_dimgen_Theta_sss}--\eqref{eqn_dimgen_Theta_uuu}. As in \eqref{C_s_pullback_b_ne_0_R_not_sqrt2} we get
	\begin{equation*}
		L(r)^*\boldsymbol\Theta_{s^3}=(3-r^2)\left(1-\frac{r}{R}\right)^{-\frac{3}{2}} C(s)\sqrt{\beta}\equiv 0,
	\end{equation*}
by equation \eqref{eqn_P_3_n_geq_3_at_least_m_ker_lambda}. Hence,
	\begin{equation}
		\lim\limits_{r\nearrow R}L(r)^*\boldsymbol\Theta_{s^3}=0.\label{eqn_Theta_sss_vanishing_ker_m_0}
	\end{equation}
Recall equation \eqref{C_ssu_b_ne0_ker_m_mplus1}, namely that $C(s,s,u)=-\frac{bR^2}{3(3-R^2)}\langle s,s\rangle \langle \widehat{\gamma},u\rangle$. This means that the pullback of $\boldsymbol\Theta_{s^2u}$ \eqref{eqn_dimgen_Theta_ssu} under $L(r)$ \eqref{eqn_B_r_ker_m_0} is of the form
	\begin{equation*}
		L(r)^*\boldsymbol\Theta_{s^2u}=\left(1-\frac{r}{R}\right)^{-1}\left(br^2-bR^2\frac{3-r^2}{3-R^2}\right)\langle s,s\rangle \langle \widehat{\gamma},u\rangle\sqrt{\beta},
	\end{equation*}
and by observing that $1-\frac{r}{R}$ has a simple zero in $r=R$ and that $br^2-bR^2\frac{3-r^2}{3-R^2}$ has at least a simple zero in $r=R$, we obtain using $\beta(R)=0$
	\begin{equation}
		\lim\limits_{r\nearrow R}L(r)^*\boldsymbol\Theta_{s^2u}=0.\label{eqn_Theta_ssu_vanishing_ker_m_0}
	\end{equation}	
The term $\boldsymbol\Theta_{su^2}$ will be considered in a later step. Since $\dim\ker\boldsymbol\lambda|_{r=R}=m$ and, hence, $\boldsymbol\lambda_{u^2}|_{r=R}>0$, it follows that
	\begin{equation*}
		\lim\limits_{r\nearrow R}L(r)^*\boldsymbol\Theta_{u^3}=\widehat{L}^*\boldsymbol\Theta_{u^3}|_{r=R}=0.
	\end{equation*}
For $\boldsymbol\Theta_{v^2w}$ \eqref{eqn_dimgen_Theta_vvw} we replace $a=\frac{R^2-1}{R^3}$ and write $\boldsymbol\Theta_{v^2w}=\boldsymbol\Theta_{s^2w}+\boldsymbol\Theta_{u^2w}$ with
	\begin{align}
		\boldsymbol\Theta_{s^2w}&=(3-r^2)^2\left(-r\frac{2R^3-3r(R^2-1)}{R^3}+\frac{(3-r^2)}{R}\right)\langle s,s\rangle w\sqrt{\beta},\label{eqn_dimgen_Theta_ssw}\\
		\boldsymbol\Theta_{u^2w}&=\left(3b^2r^5\frac{2R^3-3r(R^2-1)}{R^3}\langle \widehat{\gamma},u\rangle^2\right.\notag\\
		&\quad\left. + \langle u, (3-r^2)^2\left(-r\frac{2R^3-3r(R^2-1)}{R^3}+(3-r^2)\widehat{q}\right)u\rangle\right) w\sqrt{\beta}.\label{eqn_dimgen_Theta_uuw}
	\end{align}
For the next step we recall the calculations leading to \eqref{eqn_Theta_vvw_pullback_limit_lambda_greater_0} in the previously considered case $\boldsymbol\lambda|_{r=R}>0$ for the pullback of $\boldsymbol\Theta_{u^2w}$, and together with a direct calculation for the $\boldsymbol\Theta_{s^2w}$ verify that
	\begin{align}
		\lim\limits_{r\nearrow R}\left(L(r)\times (\mathbbm{1}_\mathbb{R}/\sqrt{\mu})\right)^*\boldsymbol\Theta_{s^2w}&=\frac{2}{\sqrt{3}}\langle s,s\rangle w,\label{eqn_Theta_ssw_limit_ker_m_0}\\
		\lim\limits_{r\nearrow R}\left(L(r)\times (\mathbbm{1}_\mathbb{R}/\sqrt{\mu})\right)^*\boldsymbol\Theta_{u^2w}&=-\frac{1}{\sqrt{3}}\langle u,u\rangle w.
	\end{align}
Next we deal with the pullback of $\boldsymbol\Theta_{vw^2}$ \eqref{eqn_dimgen_Theta_vww}. Since $\gamma$ is of the form $\gamma=\left(\begin{matrix}0\\ \widehat{\gamma}\end{matrix}\right)$, $\widehat{\gamma}\in\mathbb{R}^{n-1-m}$, we can identify $\boldsymbol\Theta_{vw^2}=\boldsymbol\Theta_{uw^2}$ since $\boldsymbol\Theta_{vw^2}$ does not depend on the $s$-coordinates. Using equation \eqref{eqn_dimgen_theta_vww_by_beta_limit} together with equation \eqref{eqn_dim2_sqrtmu_by_sqrtbeta_limit} and the fact that $\widehat{L}\in\mathrm{GL}(n-1-m)$ does not depend on $r$ shows that
	\begin{equation*}
		\lim\limits_{r\nearrow R}\left(L(r)\times (\mathbbm{1}_\mathbb{R}/\sqrt{\mu})\right)^*\boldsymbol\Theta_{vw^2}=\lim\limits_{r\nearrow R}\left(\widehat{L}\times (\mathbbm{1}_\mathbb{R}/\sqrt{\mu})\right)^*\boldsymbol\Theta_{uw^2}=0.
	\end{equation*}
Next we recall that by \eqref{eqn_dim2_sqrtmu_by_sqrtbeta_limit} and \eqref{eqn_dimgen_theta_www_by_beta_limit}
	\begin{equation*}
		\lim\limits_{r\nearrow R}\left(L(r)\times (\mathbbm{1}_\mathbb{R}/\sqrt{\mu})\right)^*\boldsymbol\Theta_{w^3}=-\frac{2}{3\sqrt{3}}w^3.
	\end{equation*}
Lastly we have to study the pullback of $\boldsymbol\Theta_{su^2}$ \eqref{eqn_dimgen_Theta_suu}. Since $\beta$ and $1-\frac{r}{R}$ both have a simple zero in $r=R$, we find that it exists independently of the initial cubic polynomial $C$ and must be of the form
	\begin{equation*}
		\lim\limits_{r\nearrow R}\left(L(r)\times (\mathbbm{1}_\mathbb{R}/\sqrt{\mu})\right)^*\boldsymbol\Theta_{su^2}=\langle s,F(u)\rangle,
	\end{equation*}
where $F=(F_1,\ldots,F_m)^\mathrm{T}:\mathbb{R}^{n-1-m}\to\mathbb{R}^m$ is a quadratic homogeneous polynomial in each entry. The terms $F_i$ must fulfill the condition in Theorem \ref{thm_limit_behaviour_ccpsr} \eqref{thm_limit_behaviour_ccpsr_iii}, cf. Section \ref{sect_Fi_calcs},
otherwise $\mathcal{H}\subset\mathbb{R}^{n+1}$ cannot be closed by Theorem \ref{thm_Cn}, thereby violating its initial condition. Summarizing, we have shown that in the case $\dim\ker \boldsymbol\lambda|_{r=R}=m$, $q_1=\ldots=q_m=\frac{1}{R}>q_{m+1}\geq\ldots\geq q_{n-1}$ for $1\leq m\leq n-2$, $b\ne 0$, the limit polynomial, that is the pullback of $h_E$ \eqref{eqn_dimgen_hE1_prelimit} with respect to $B(r)\times (\mathbbm{1}_\mathbb{R}/\sqrt{\mu})$ as $r\to R$, is given by
	\begin{align}
		\lim\limits_{r\nearrow R}\left(\mathbbm{1}_\mathbb{R}\times L(r)\times (\mathbbm{1}_\mathbb{R}/\sqrt{\mu})\right)^* h_E&=x^3-x\left(\langle s,s\rangle + \langle u,u\rangle + w^2\right)\notag\\
		&\quad + \left(\frac{2}{\sqrt{3}}\langle s,s\rangle -\frac{1}{\sqrt{3}}\langle u,u\rangle \right)w +\langle s,F(u)\rangle -\frac{2}{3\sqrt{3}}w^3.\label{eqn_limit_poly_ker_m}
	\end{align}
Note that at this point it is not clear whether or not a CCPSR manifold $\mathcal{H}$ fulfilling the initial conditions $b\ne 0$, $\dim\ker \boldsymbol\lambda|_{r=R}=m$, $q_1=\ldots=q_m=\frac{1}{R}$ does exist in every dimension $n=\dim(\mathcal{H})\geq 3$. To clear this up, we will now construct such an example in each dimension.
\begin{Ex}\label{example_Rsqrt2_dimker_m}
	Let $R=\sqrt{2}$, $b$ fixed with $0<|b|<\frac{\sqrt{10}}{12}$, $1\leq m\leq n-2$, and consider $h$ as in equation \eqref{eqn_h_n_geq_3_at_least_m_ker_lambda} with
		\begin{equation}\label{eqn_P3_sqrt2_example_general_dim}
			P_3\left(\left(\begin{smallmatrix}s\\ u\\ w\end{smallmatrix}\right)\right)=-2b \langle s,s\rangle\langle\widehat{\gamma},u\rangle+\frac{1}{\sqrt{2}}\langle s,s\rangle w + b\langle\widehat{\gamma},u\rangle w^2 + \frac{1}{2\sqrt{2}}w^3.
		\end{equation}
	The above defined $P_3$ fulfils equation \eqref{eqn_P_3_n_geq_3_at_least_m_ker_lambda} for $C(v)=-\frac{bR^2}{3-R^2}\langle s,s\rangle \langle \widehat{\gamma},u\rangle$ and $\widehat{q}=0$. We want to show that the connected component $\mathcal{H}\subset\{h=1\}$ that contains the point $(x,s,u,w)=(1,0,0,0)$ is closed in the ambient space independent of the choice for $\widetilde{\gamma}$, which by Theorem \ref{thm_Cn} is equivalent to showing that $\max\limits_{\{\langle s,s\rangle + \langle u,u\rangle + w^2=1\}}P_3\left(\left(\begin{smallmatrix}s\\ u\\ w\end{smallmatrix}\right)\right)\leq \frac{2}{3\sqrt{3}}$. Since $n=\dim(\mathcal{H})\geq 3$ by assumption, it is clear that this is equivalent to showing that for all $3$-dimensional linear subspaces $V\subset\mathbb{R}^{n}$ of the form
		\begin{equation*}
			V=\mathrm{span}\left.\left\{\left(\begin{smallmatrix} s_0\\ 0\\ 0\end{smallmatrix}\right),\ \left(\begin{smallmatrix} 0\\ u_0\\ 0\end{smallmatrix}\right),\ \left(\begin{smallmatrix} 0\\ 0\\ 1\end{smallmatrix}\right)\ \right|\ s_0\ne 0,\ u_0\ne 0\right\}
		\end{equation*}
	it holds that $\max\limits_{V\cap\{\langle s,s\rangle + \langle u,u\rangle + w^2=1\}}P_3\left(\left(\begin{smallmatrix}s\\ u\\ w\end{smallmatrix}\right)\right)\leq \frac{2}{3\sqrt{3}}$. This however means that we can restrict ourselves to the case $\dim(\mathcal{H})=3$, allowing additionally $b=0$, to prove our claim that that $\mathcal{H}$ is closed in general dimension $n\geq 3$.	For $n=3$, $P_3$ is of the form
		\begin{equation}\label{eqn_P3_3dim_reduction_sqrt2_case}
			P_3\left(\left(\begin{smallmatrix}s\\ u\\ w\end{smallmatrix}\right)\right)=-2bs^2u+\frac{1}{\sqrt{2}}s^2 w + buw^2 + \frac{1}{2\sqrt{2}}w^3.
		\end{equation}
	We need to show that for all $|b|<\frac{\sqrt{10}}{12}$, $\max\limits_{\{s^2 + u^2 + w^2=1\}}P_3\left(\left(\begin{smallmatrix}s\\ u\\ w\end{smallmatrix}\right)\right)\leq \frac{2}{3\sqrt{3}}$. This is equivalent to showing that for all $p=(p_s,p_u,p_w)^\mathrm{T}\in\mathbb{R}^3\setminus{0}$ solving
		\begin{equation}\label{eqn_sqrt2_example_crit_points}
			\D P_3|_p=\frac{2}{\sqrt{3}}(p_s \,\D s + p_u \,\D u + p_w \,\D w),
		\end{equation}
	$p$ has Euclidean norm of at least $1$, that is $\|p\|\geq 1$. This follows from the homogeneity of $P_3$. With the help of a computer algebra system, e.g. Maple, we find that for $|b|>0$ the solution set of equation \eqref{eqn_sqrt2_example_crit_points} in $\mathbb{R}^3\setminus\{0\}$ is given by
		\begin{equation}\label{eqn_solution_set_sqrt2_dim2}
			\left\{
				\left(\begin{matrix} \frac{1}{\sqrt{3}}\\ 0\\ \frac{\sqrt{2}}{\sqrt{3}}\end{matrix}\right),\ \left(\begin{matrix} -\frac{1}{\sqrt{3}}\\ 0\\ \frac{\sqrt{2}}{\sqrt{3}}\end{matrix}\right),\ \left(\begin{matrix} 0\\ u_1\\ w_1\end{matrix}\right),\ \left(\begin{matrix} 0\\ u_2\\ w_2\end{matrix}\right)
			\right\},
		\end{equation}
	where
		\begin{equation*}
			\begin{array}{ll}
				u_1=\frac{\left(-3+\sqrt{64b^2+9}\right)^2}{64\sqrt{3}b^3},& w_1=\frac{-3+\sqrt{64b^2+9}}{4\sqrt{6}b^2},\\
				u_2=\frac{\left(3+\sqrt{64b^2+9}\right)^2}{64\sqrt{3}b^3},& w_1=\frac{-3-\sqrt{64b^2+9}}{4\sqrt{6}b^2}.
			\end{array}
		\end{equation*}
	The first two points in the above solution set \eqref{eqn_solution_set_sqrt2_dim2} are of Euclidean norm $1$. Note that this implies by the homogeneity of $P_3$ and equation \eqref{eqn_sqrt2_example_crit_points} that $\max\limits_{\{s^2 + u^2 + w^2=1\}}P_3\left(\left(\begin{smallmatrix}s\\ u\\ w\end{smallmatrix}\right)\right)\geq \frac{2}{3\sqrt{3}}$, so if $\mathcal{H}$ is closed it is also singular at infinity, cf. Lemma \ref{lem_sing_at_infty_iff_max_cond}. For $\mathcal{H}$ being closed it remains to verify that for all $|b|<\frac{\sqrt{10}}{12}$,
		\begin{equation*}
			u_1^2+w_1^2\geq 1,\quad u_2^2+w_2^2\geq 1.
		\end{equation*}
	We find
		\begin{align}
			u_1^2+w_1^2&=\frac{\left(-3+\sqrt{64b^2+9}\right)^2\left(32b^2+3-\sqrt{64b^2+9}\right)}{2048b^6}, \label{eqn_u1_w1_normsquare}\\
			u_2^2+w_2^2&=\frac{\left(3+\sqrt{64b^2+9}\right)^2\left(32b^2+3+\sqrt{64b^2+9}\right)}{2048b^6}. \label{eqn_u2_w2_normsquare}
		\end{align}
	It is now a slightly tedious, but not difficult task to check that the expressions in \eqref{eqn_u1_w1_normsquare} and \eqref{eqn_u2_w2_normsquare} are both strictly monotonously decreasing functions in $b$ for $b\in\left(0,\frac{\sqrt{10}}{12}\right)$, and symmetric in $b$. It thus suffices to check that $\left(u_1^2+w_1^2\right)|_{b=\frac{\sqrt{10}}{12}}=\frac{126}{125}>1$ and $\left(u_2^2+w_2^2\right)|_{b=\frac{\sqrt{10}}{12}}=576>1$. Hence, we have shown that $\max\limits_{\{s^2 + u^2 + w^2=1\}}P_3\left(\left(\begin{smallmatrix}s\\ u\\ w\end{smallmatrix}\right)\right)= \frac{2}{3\sqrt{3}}$ for all $0<|b|<\frac{\sqrt{10}}{12}$. It remains to check that $P_3$ as in \eqref{eqn_P3_3dim_reduction_sqrt2_case} also fulfils $\max\limits_{\{s^2 + u^2 + w^2=1\}}P_3\left(\left(\begin{smallmatrix}s\\ u\\ w\end{smallmatrix}\right)\right)\leq \frac{2}{3\sqrt{3}}$ for $b=0$. We leave that as an easy exercise for the reader.
	
	Summarizing, we have shown that the connected component $\mathcal{H}\subset\{h=1\}$, $h$ defined by $P_3$ as in equation \eqref{eqn_P3_sqrt2_example_general_dim}, that contains the point $(x,s,u,w)^\mathrm{T}=(1,0,0,0)^\mathrm{T}$ is indeed a CCPSR manifold for all initial $b\in\left(-\frac{\sqrt{10}}{12},\frac{\sqrt{10}}{12}\right)$. Note at this point that the range of $b$ has been chosen so that the solution set \eqref{eqn_solution_set_sqrt2_dim2} for critical points of $P_3$ consists of precisely $4$ different points for all $b\ne 0$. It is not excluded here that there exists $|b|\geq\frac{\sqrt{10}}{12}$, so that $\mathcal{H}$ is still closed its ambient space, but the calculations become a lot more complicated.
\end{Ex}

Next we will deal with the case $\dim\ker \boldsymbol\lambda|_{r=R}=m+1$ for $1\leq m\leq n-2$. Again, we use the notation in \eqref{eqn_s_u_def} so that $\boldsymbol\lambda=\boldsymbol\lambda_{s^2}+\boldsymbol\lambda_{u^2}$, cf. equations \eqref{eqn_lambda_ss} and \eqref{eqn_lambda_uu}. Furthermore we have formally the same expressions for $\boldsymbol\Theta_{s^3}$--$\boldsymbol\Theta_{u^3}$ \eqref{eqn_dimgen_Theta_sss}--\eqref{eqn_dimgen_Theta_uuu}, $\boldsymbol\Theta_{s^2w}$ \eqref{eqn_dimgen_Theta_ssw}, $\boldsymbol\Theta_{u^2w}$ \eqref{eqn_dimgen_Theta_uuw}, and $\boldsymbol\Theta_{w^3}$ \eqref{eqn_dimgen_Theta_www}. We will however, in comparison with the case $\dim\ker \boldsymbol\lambda|_{r=R}=m$ and $q_1=\ldots=q_m=\frac{1}{R}$, need to more carefully consider the terms
	\begin{equation}
		\boldsymbol\Theta_{uw^2}= 9b\left(3+ r^2-3ar^3\right)\langle\widehat{\gamma},u\rangle w^2\beta\sqrt{\beta},\label{eqn_dimgen_Theta_uww}
	\end{equation}
which stems from $\boldsymbol\Theta_{vw^2}=\boldsymbol\Theta_{sw^2}+\boldsymbol\Theta_{uw^2}$, cf. equations \eqref{eqn_dimgen_Theta_vww} and \eqref{eqn_gamma_ker_m_and_m_1_cases}, and
	\begin{equation}
		\boldsymbol\chi=-18br\beta\langle\widehat{\gamma},u\rangle w,
	\end{equation}
cf. equations \eqref{eqn_dimgen_chi_def} and \eqref{eqn_gamma_ker_m_and_m_1_cases}. For each equation that contains the $u$-coordinates we will proceed similarly to the case $\frac{1}{R}>q_1\geq\ldots\geq q_{n-1}$ for all $1\leq i\leq n-1$, cf. equation \eqref{eqn_ker_0_1_eta_trafo} onward. In particular we have from equation \eqref{eqn_ker_0_1_condition} the requirement
	\begin{equation}
		\left\|\frac{1}{\sqrt{1-\widehat{q}R}}\widehat{\gamma}\right\|^2=\frac{(3-R^2)^2}{3b^2R^4},\label{eqn_ker_m_1_condition}
	\end{equation}
and we obtain an similar equality as in \eqref{eqn_C_1_by_1_minus_qR_gamma}. Analogously to \eqref{eqn_ker_0_1_B_gamma_condition} we choose a smooth map $\widehat{B}:[0,R]\to\mathrm{SO}(n-1-m)$, such that
	\begin{equation}
		\widehat{B}(r)^\mathrm{T}\frac{1}{\sqrt{1-\widehat{q}r}}\widehat{\gamma}=\left\|\frac{1}{\sqrt{1-\widehat{q}r}}\widehat{\gamma}\right\|\partial_{n-1}.\label{eqn_ker_m_1_B_gamma_condition}
	\end{equation}
As in equation \eqref{eqn_L_trafo} we define
	\begin{equation}
		\widehat{L}:[0,R)\to\mathrm{GL}(n-1-m),\quad \widehat{L}(r)\left(\begin{matrix}\rho\\ \alpha \end{matrix}\right)=\frac{1}{(3-r^2)\sqrt{1-\widehat{q}r}}\widehat{B}(r)\left(\begin{matrix}\rho\\ \frac{\alpha}{\sqrt{1-\left\|\widehat{\nu}\right\|^2}}\end{matrix}\right).\label{eqn_widehat_L_trafo}
	\end{equation}
where
	\begin{equation*}
		\widehat{\nu}=\widehat{\nu}(r):=\frac{\sqrt{3}br^2}{(3-r^2)\sqrt{1-\widehat{q}r}}\widehat{\gamma},\quad r\in[0,R],\quad \|\widehat{\nu}\|^2=\frac{3b^2r^4}{(3-r^2)^2}\left\|\frac{1}{\sqrt{1-\widehat{q}r}}\widehat{\gamma}\right\|^2,
	\end{equation*}
cf. \eqref{eqn_nu_def}. For the following calculations we observe that
	\begin{equation*}
		\widehat{L}(r)^*\langle \widehat{\gamma},u\rangle=\frac{1}{(3-r^2)\sqrt{1-\left\|\widehat{\nu}\right\|^2}}\left\|\frac{1}{\sqrt{1-\widehat{q}r}}\widehat{\gamma}\right\|\alpha,
	\end{equation*}
and we obtain with the notation $w=\frac{1}{\sqrt{\mu}}\overline{w}$
	\begin{align*}
		&\lim\limits_{r\nearrow R}\left(\widehat{L}(r)\times (\mathbbm{1}_\mathbb{R}/\sqrt{\mu})\right)^*\boldsymbol\chi =\boldsymbol\zeta,\\
		&\lim\limits_{r\nearrow R} \left(\widehat{L}(r)\times (\mathbbm{1}_\mathbb{R}/\sqrt{\mu})\right)^*\left(\boldsymbol\Theta_{u^3}+\boldsymbol\Theta_{u^2w}+\boldsymbol\Theta_{uw^2}+\boldsymbol\Theta_{w^3}\right)\\
		&=\boldsymbol\Xi_{\rho^2\alpha}+\boldsymbol\Xi_{\alpha^3}+\boldsymbol\Xi_{\rho^2\overline{w}}+\boldsymbol\Xi_{\alpha^2\overline{w}}+\boldsymbol\Xi_{\alpha\overline{w}^2}+\boldsymbol\Xi_{\overline{w}^3},
	\end{align*}
where $\boldsymbol\zeta$ and $\boldsymbol\Xi_{\rho^2\alpha}$ -- $\boldsymbol\Xi_{\overline{w}^3}$ are symbolically of the same form as their pendants for the case $\frac{1}{R}>q_1\geq\ldots\geq q_{n-1}$ and $\dim\ker \boldsymbol\lambda|_{r=R}=1$, cf. equations \eqref{eqn_Theta_alphaalphaoverlinew_limit}, \eqref{eqn_Theta_www_limit_L_trafo_case}, and \eqref{eqn_zeta_ker_0_1}--\eqref{eqn_Xi_alphaoverlinewoverlinew}. In the following, let
	\begin{equation}
		L:[0,R)\to\mathrm{GL}(n-1),\quad L(r):= \left(\begin{array}{c|c} (3-r^2)^{-1}\left(1-\frac{r}{R}\right)^{-\frac{1}{2}}\mathbbm{1}_m \underset{}{} &  \\ \hline & \overset{}{\widehat{L}(r)} \end{array}\right),\quad \left(\begin{matrix}s\\ u\end{matrix}\right)=L(r)\cdot\left(\begin{matrix}\overline{s}\\ \rho\\ \alpha\end{matrix}\right),\label{eqn_L_r_ker_m_1}
	\end{equation}
$\overline{s}=(\overline{s}_1,\ldots,\overline{s}_m)^\mathrm{T}$, $\rho=(\rho_1,\ldots,\rho_{n-2-m})^\mathrm{T}$,
so that
	\begin{equation*}
		\lim\limits_{r\nearrow R} L(r)^*\left(\boldsymbol\lambda_{s^2}+\boldsymbol\lambda_{u^2}\right)=\langle\overline{s},\overline{s}\rangle + \langle \rho,\rho\rangle + \alpha^2.
	\end{equation*}
Note that formally $L(r)$ as in \eqref{eqn_B_r_ker_m_0} and $L(r)$ as in \eqref{eqn_L_r_ker_m_1} coincide as linear transformation on the $s$-coordinates. Recall that $C(s)\equiv 0$, cf. equation \eqref{eqn_P_3_n_geq_3_at_least_m_ker_lambda}, and, hence,
	\begin{equation}
		\lim\limits_{r\nearrow R} \left(L(r)\times (\mathbbm{1}_\mathbb{R}/\sqrt{\mu})\right)^*\boldsymbol\Theta_{s^3}=0.\label{eqn_Theta_sss_vanishing_ker_m_1}
	\end{equation}
Furthermore we find as for \eqref{eqn_Theta_ssw_limit_ker_m_0}
	\begin{equation}
		\lim\limits_{r\nearrow R}\left(L(r)\times (\mathbbm{1}_\mathbb{R}/\sqrt{\mu})\right)^*\boldsymbol\Theta_{s^2w}=\frac{2}{\sqrt{3}}\langle \overline{s},\overline{s}\rangle \overline{w}.\label{eqn_Theta_ssw_limit_ker_m_1}
	\end{equation}
Next, we will study the limit of the term $\boldsymbol\Theta_{s^2u}$ \eqref{eqn_dimgen_Theta_ssu} transformed under $L$ \eqref{eqn_L_r_ker_m_1} as $r\to R$. We have
	\begin{align}
		L(r)^*\boldsymbol\Theta_{s^2u}&= \left(\frac{br^2}{(3-r^2)(1-r/R)\sqrt{1-\left\|\widehat{\nu}\right\|^2}}\left\|\frac{1}{\sqrt{1-\widehat{q}r}}\widehat{\gamma}\right\|\langle\overline{s},\overline{s}\rangle\alpha\right.\notag \\
		&\quad\left. \vphantom{\frac{br^2}{(3-r^2)(1-r/R)\sqrt{1-\left\|\widehat{\nu}\right\|^2}}\left\|\frac{1}{\sqrt{1-\widehat{q}r}}\widehat{\gamma}\right\|langle\overline{s},\overline{s}\rangle\alpha} + \frac{3}{1-r/R}C\left(\overline{s},\overline{s},\frac{1}{\sqrt{1-\widehat{q}r}}\widehat{B}(r)\left(\begin{matrix}\rho\\ \frac{\alpha}{\sqrt{1-\left\|\widehat{\nu}\right\|^2}}\end{matrix}\right)\right)\right)\sqrt{\beta}.\label{eqn_pre_limit_Theta_ssu_ker_m_1}
	\end{align}
Since both $1-\frac{r}{R}$ and $\beta $ have simple zeros in $r=R$ and $\mathcal{H}$ is assumed to be closed, it follows from equation \eqref{eqn_pre_limit_Theta_ssu_ker_m_1} restricted to $\{\alpha=0\}$ and Theorem \ref{thm_Cn} that
	\begin{equation}
		C\left(\overline{s},\overline{s},\frac{1}{\sqrt{1-\widehat{q}R}}\widehat{B}(R)\left(\begin{matrix}\rho\\ 0\end{matrix}\right)\right)=0\label{eqn_vanishing_C_ssrho_ker_m_1}
	\end{equation}
for all $\overline{s}\in\mathbb{R}^{m}$, $\rho\in\mathbb{R}^{n-2-m}$. With the same argument we obtain from \eqref{eqn_pre_limit_Theta_ssu_ker_m_1} restricted to $\{\rho=0\}$ that the term
	\begin{equation*}
		\lim\limits_{r\nearrow R} \left.\left(\frac{1}{\sqrt{\beta}}\left(1-\frac{r}{R}\right)\sqrt{1-\|\widehat{\nu}\|^2}L(r)^*\boldsymbol\Theta_{s^2u}\right)\right|_{\{\rho=0\}}
	\end{equation*}
must vanish identically. Since $R<\sqrt{3}$, this is equivalent to
	\begin{equation}
		\left.\left(\frac{br^2}{3-r^2}\left\|\frac{1}{\sqrt{1-\widehat{q}r}}\widehat{\gamma}\right\| \langle \overline{s},\overline{s}\rangle + 3 C\left(\overline{s},\overline{s},\frac{1}{\sqrt{1-\widehat{q}r}}\widehat{B}(r) \partial_{n-1}\right)\right)\right|_{r=R}=0,\label{eqn_theta_ssu_vanishing_ker_m_1}
	\end{equation}
and by $\widehat{q}<\frac{1}{R}\langle\cdot,\cdot\rangle$ we see that the left hand side of \eqref{eqn_theta_ssu_vanishing_ker_m_1} has a simple zero in $r=R$. Thus we will need to calculate its $r$-derivative in $r=R$ to obtain a formula for the limit $r\to R$ of \eqref{eqn_pre_limit_Theta_ssu_ker_m_1}. To do so, we will first introduce the following abbreviations which are analogous to \eqref{eqn_K_def} and \eqref{eqn_Y_def}
	\begin{align}
		\widehat{K}&:=\left\|\frac{1}{\sqrt{1-\widehat{q}R}}\widehat{\gamma}\right\|^{-2}\left\|\frac{1}{1-\widehat{q}R}\widehat{\gamma}\right\|^2,\label{eqn_widehat_K_def}\\
		\widehat{Y}&:=9+R^2+(3-R^2)\widehat{K},\label{eqn_widehat_Y_def}
	\end{align}
and as for equation \eqref{eqn_t_Z_split} we check that
	\begin{equation}
		\frac{1}{(1-\widehat{q}R)^2}\widehat{\gamma}=\frac{\widehat{K}}{1-\widehat{q}R} \widehat{\gamma} + \frac{1}{\sqrt{1-\widehat{q}R}}\widehat{B}(R)\left(\begin{smallmatrix}Z\\ 0\end{smallmatrix}\right)\label{eqn_t_Z_split_ker_m_1}
	\end{equation}
with $Z\in\mathbb{R}^{n-2-m}$ uniquely determined. As for equation \eqref{eqn_sqrt_beta_sqrt_1minusnusquare_limit_alternative} we obtain
	\begin{equation}
		\lim\limits_{r\nearrow R}\frac{\sqrt{\beta}}{\sqrt{1-\|\widehat{\nu}\|^2}}=\frac{3-R^2}{\sqrt{\widehat{Y}}}\label{eqn_sqrt_beta_sqrt_1minusnusquare_limit_alternative_ker_m_1}.
	\end{equation}
Now we use equations \eqref{eqn_ker_m_1_B_gamma_condition}, \eqref{eqn_vanishing_C_ssrho_ker_m_1}, and \eqref{eqn_t_Z_split_ker_m_1} and analogous versions of equations \eqref{eqn_r_der_sqrt_stuff_norm} and \eqref{eqn_r_der_stuff_norm} to obtain
	\begin{align}
		&\partial_r  \left.\left(\frac{1}{\sqrt{\beta}}\left(1-\frac{r}{R}\right)\sqrt{1-\|\widehat{\nu}\|^2}L(r)^*\boldsymbol\Theta_{s^2u}\right)\right|_{r=R}\notag\\
		&=\left(\frac{6bR}{(3-R^2)^2}+\frac{bR(-1+\widehat{K})}{2(3-R^2)}\right)\left\|\frac{1}{\sqrt{1-\widehat{q}R}}\widehat{\gamma}\right\|\langle \overline{s},\overline{s}\rangle\notag\\
		&\quad +  \frac{3(-1+\widehat{K})}{2R}\left\|\frac{1}{\sqrt{1-\widehat{q}R}}\widehat{\gamma}\right\|^{-1}C\left(\overline{s},\overline{s},\frac{1}{1-\widehat{q}R}\widehat{\gamma}\right).\label{eqn_partial_r_Cssu_ker_m_1}
	\end{align}
Using \eqref{eqn_vanishing_C_ssrho_ker_m_1} and \eqref{eqn_t_Z_split_ker_m_1}, equation \eqref{eqn_theta_ssu_vanishing_ker_m_1} can be rewritten as
	\begin{equation}
		3\left\|\frac{1}{\sqrt{1-\widehat{q}R}}\widehat{\gamma}\right\|^{-1}C\left(\overline{s},\overline{s},\frac{1}{1-\widehat{q}R}\widehat{\gamma}\right)=-\frac{bR^2}{3-R^2}\left\|\frac{1}{\sqrt{1-\widehat{q}R}}\widehat{\gamma}\right\|\langle \overline{s},\overline{s}\rangle.\label{eqn_Cssu_vanishing_ker_m_1}
	\end{equation}
Note at this point that both \eqref{eqn_vanishing_C_ssrho_ker_m_1} and \eqref{eqn_Cssu_vanishing_ker_m_1} can also be immediately derived by using \eqref{C_ssu_b_ne0_ker_m_mplus1}. We have however decided to keep this way of showing these equivalences in our new coordinates $\overline{s},\rho,\alpha$ since it provides an additional safety check that everything is in order, given the technical nature of this proof. Furthermore the required extra space is small. By inserting \eqref{eqn_Cssu_vanishing_ker_m_1} into \eqref{eqn_partial_r_Cssu_ker_m_1}, using \eqref{eqn_sqrt_beta_sqrt_1minusnusquare_limit_alternative_ker_m_1}, $\partial_r(1-r/R)|_{r=R}=-\frac{1}{R}$, and with the help of L'H\^opital's rule we finally obtain
	\begin{equation}
		\boldsymbol\Xi_{\overline{s}^2\alpha}:=\lim\limits_{r\nearrow R} L(r)^*\boldsymbol\Theta_{s^2u} = \mathrm{sgn}(b)\frac{-2\sqrt{3}}{\sqrt{\widehat{Y}}}\langle \overline{s},\overline{s}\rangle \alpha.\label{eqn_Xi_overlinesoverlinesalpha}
	\end{equation}
Note that the above formula \eqref{eqn_Xi_overlinesoverlinesalpha} is consistent with equation \eqref{C_ssu_b_ne0_ker_m_mplus1} for $u=\frac{1}{1- \widehat{q}R}\widehat{\gamma}$. Next, we will determine the limit of the pullback of $\boldsymbol\Theta_{su^2}$ \eqref{eqn_dimgen_Theta_suu} as $r\to R$. We have
	\begin{equation}
		L(r)^* \boldsymbol\Theta_{su^2} = \frac{3}{\sqrt{1-\frac{r}{R}}} C\left(\overline{s}, \frac{1}{\sqrt{1-\widehat{q}r}}\widehat{B}(r)\left(\begin{matrix} \rho\\ \frac{\alpha}{\sqrt{1-\|\widehat{\nu}\|^2}} \end{matrix}\right), \frac{1}{\sqrt{1-\widehat{q}r}}\widehat{B}(r)\left(\begin{matrix} \rho\\ \frac{\alpha}{\sqrt{1-\|\widehat{\nu}\|^2}} \end{matrix}\right)\right) \sqrt{\beta}, \label{eqn_Theta_suu_pre_limit_ker_m_1}
	\end{equation}
and again with the argument that the initial connected PSR manifold $\mathcal{H}$ is assumed to be closed, Theorem \ref{thm_Cn}, and that $1-r/R$ and $\beta$ both have simple zeros in $r=R$, we deduce that every monomial in \eqref{eqn_Theta_suu_pre_limit_ker_m_1} containing $\alpha$ or $\alpha^2$ must vanish identically near $r=R$ and, hence, will also vanish in the limit $r\to R$. We obtain with the help of L'H\^opital's rule
	\begin{align}
		\boldsymbol\Xi_{\overline{s}\rho^2}:=& \lim\limits_{r\nearrow R} \left(L(r)^* \boldsymbol\Theta_{su^2}\right)\notag\\
		=&\ 3\sqrt{3-R^2} C\left(\overline{s}, \frac{1}{\sqrt{1-\widehat{q}R}}\widehat{B}(R)\left(\begin{matrix} \rho\\ 0 \end{matrix}\right), \frac{1}{\sqrt{1-\widehat{q}R}}\widehat{B}(R)\left(\begin{matrix} \rho\\ 0 \end{matrix}\right)\right).\label{eqn_Xi_overlinesrhorho}
	\end{align}
We further define symmetric bilinear forms $\widetilde{W}_i:\mathbb{R}^{n-2-m}\times\mathbb{R}^{n-2-m}\to \mathbb{R}$, $1\leq i\leq m$, such that with $\widetilde{W}(\rho,\rho)=\left(\widetilde{W}_1(\rho,\rho),\ldots,\widetilde{W}_m(\rho,\rho)\right)^\mathrm{T}$
	\begin{equation}
		\boldsymbol\Xi_{\overline{s}\rho^2}=\left\langle \overline{s},\widetilde{W}(\rho,\rho)\right\rangle.\label{eqn_widetilde_W_def}
	\end{equation}
For the limit of the pullback of $\boldsymbol\Theta_{s^2w}$ \eqref{eqn_dimgen_Theta_ssw} the same calculation as for equation \eqref{eqn_Theta_ssw_limit_ker_m_0} yields
	\begin{equation}
		\boldsymbol\Xi_{\overline{s}^2\overline{w}}:= \lim\limits_{r\nearrow R}\left(L(r)\times (\mathbbm{1}_\mathbb{R}/\sqrt{\mu})\right)^*\boldsymbol\Theta_{s^2w} =\frac{2}{\sqrt{3}}\langle \overline{s},\overline{s}\rangle \overline{w}.\label{eqn_Xi_overlinesoverlinesoverlinew}
	\end{equation}
Summarizing up to this point, we have shown that the limit of the pullback of $h_E$ \eqref{eqn_dimgen_hE1_prelimit} under $\mathbbm{1}_\mathbb{R}\times L(r)\times (\mathbbm{1}_\mathbb{R}/\sqrt{\mu})$ \eqref{eqn_L_r_ker_m_1} is given by
		\begin{align}
			&\lim\limits_{r\nearrow R} \left(\mathbbm{1}_\mathbb{R}\times L(r)\times (\mathbbm{1}_\mathbb{R}/\sqrt{\mu})\right)^* h_{E}\left(\left(\begin{smallmatrix} x\\ \overline{s}\\ \rho\\ \alpha\\ \overline{w}\end{smallmatrix}\right)\right)\notag\\
			&=x^3-x\left(\langle \overline{s},\overline{s}\rangle + \langle \rho,\rho\rangle + \alpha^2 + \boldsymbol\zeta  + \overline{w}^2\right)\notag\\
			&\quad + \boldsymbol\Xi_{\overline{s}^2\alpha} +  \boldsymbol\Xi_{\overline{s}^2\overline{w}} + \boldsymbol\Xi_{\overline{s}\rho^2} + \boldsymbol\Xi_{\rho^2\alpha} + \boldsymbol\Xi_{\alpha^3} + \boldsymbol\Xi_{\rho^2\overline{w}} + \boldsymbol\Xi_{\alpha^2\overline{w}} + \boldsymbol\Xi_{\alpha\overline{w}^2} + \boldsymbol\Xi_{\overline{w}^3},\label{eqn_limit_ker_m_1_case_pre_trafo}
		\end{align}
where (cf. equations \eqref{eqn_Xi_overlinesoverlinesalpha}, \eqref{eqn_Xi_overlinesoverlinesoverlinew}, and \eqref{eqn_widetilde_W_def})
	\begin{equation*}
		\boldsymbol\Xi_{\overline{s}^2\alpha}= \mathrm{sgn}(b)\frac{-2\sqrt{3}}{\sqrt{\widehat{Y}}}\langle \overline{s},\overline{s}\rangle \alpha,\quad
		\boldsymbol\Xi_{\overline{s}^2\overline{w}}= \frac{2}{\sqrt{3}}\langle \overline{s},\overline{s}\rangle \overline{w},\quad
		\boldsymbol\Xi_{\overline{s}\rho^2}= \left\langle \overline{s},\widetilde{W}(\rho,\rho)\right\rangle,
	\end{equation*}
$\boldsymbol\zeta$, $\boldsymbol\Xi_{\rho^2\alpha}$--$\boldsymbol\Xi_{\alpha\overline{w}^2}$ are of the form \eqref{eqn_zeta_ker_0_1}--\eqref{eqn_Xi_alphaoverlinewoverlinew} up to replacing $K$ with $\widehat{K}$ and $Y$ with $\widehat{Y}$, and $\boldsymbol\Xi_{\overline{w}^3}=-\frac{2}{3\sqrt{3}}\overline{w}^3$ which follows from the same calculation as for equation \eqref{eqn_Theta_www_limit_L_trafo_case}. Next, we recall the transformations $M$ \eqref{eqn_M_def} and $N$ \eqref{eqn_N_def} and let $T$, $t$ be as in \eqref{eqn_T_t_dimgen}. We verify that
	\begin{equation*}
		M\cdot N=\left(\begin{matrix} 0 & \frac{1}{2}\left(\frac{\sqrt{T}}{\sqrt{t}} + \frac{\sqrt{t}}{\sqrt{T}}\right)\\
		1 & \frac{1}{2}\left(\frac{\sqrt{T}}{\sqrt{t}} - \frac{\sqrt{t}}{\sqrt{T}}\right) \end{matrix}\right)
	\end{equation*}
and calculate for the final pullback of \eqref{eqn_limit_ker_m_1_case_pre_trafo}, which is analogous to \eqref{eqn_limit_ker_1_case_post_trafo1} and \eqref{eqn_limit_poly_ker_0_1}, but in one step,
	\begin{align}
		&\left(\mathbbm{1}_{\mathbb{R}^{n-1}}\times M\cdot N\right)^*\lim\limits_{r\nearrow R} \left(\mathbbm{1}_\mathbb{R}\times L(r)\times (\mathbbm{1}_\mathbb{R}/\sqrt{\mu})\right)^* h_{E}\left(\left(\begin{smallmatrix} x\\ \overline{s}\\ \rho\\ \widetilde{k}\\ \widetilde{\ell}\end{smallmatrix}\right)\right)\notag\\
		&=x^3-x\left(\langle \overline{s},\overline{s}\rangle + \langle \rho,\rho\rangle + \widetilde{\ell}^2+\widetilde{k}^2\right)\notag\\
		&\quad+ \frac{1}{\sqrt{tT}}\left(\frac{W(\rho,\rho)}{\sqrt{\widehat{Y}}} + \frac{\zeta}{2\sqrt{3}}\langle \rho,\rho\rangle\right)\widetilde{\ell} + \left\langle \overline{s},\widetilde{W}(\rho,\rho)\right\rangle\notag\\
		&\quad+\frac{2}{\sqrt{3}}\langle\overline{s},\overline{s}\rangle\widetilde{k} -\frac{1}{\sqrt{3}}\langle\rho,\rho\rangle\widetilde{k} + \frac{2}{\sqrt{3}}\widetilde{\ell}^2\widetilde{k} - \frac{2}{3\sqrt{3}}\widetilde{k}^3.\label{eqn_limit_poly_m_plus_1}
	\end{align}
By grouping the coordinates $s_1,\ldots,s_m$ and $\widetilde{\ell}$ together, we see that the limit polynomial \eqref{eqn_limit_poly_m_plus_1} is indeed of the form Thm. \ref{thm_limit_behaviour_ccpsr} \eqref{thm_limit_behaviour_ccpsr_iii}.

As for the case $\dim\ker \boldsymbol\lambda|_{r=R}=m$, cf. Example \ref{example_Rsqrt2_dimker_m}, we do not know at this point whether a CCPSR manifold $\mathcal{H}\subset\{h=1\}$ with $h$ of the form \eqref{eqn_h_n_geq_3_at_least_m_ker_lambda} with corresponding $P_3$-term as in \eqref{eqn_P_3_n_geq_3_at_least_m_ker_lambda} exists, so that $\dim\ker \boldsymbol\lambda|_{r=R}=m+1$, $q_1=\ldots=q_m=\frac{1}{R}>q_{m+1}\geq \ldots g_{n-1}$ are satisfied. We will now construct such an example for each dimension $n=\dim(\mathcal{H})\geq 3$. By Theorem \ref{thm_Cn} and equation \eqref{eqn_ker_m_1_condition}, this task is equivalent to finding $P_3\left(\left(\begin{smallmatrix}s\\ u\\ w\end{smallmatrix}\right)\right)$ for any $n\geq 3$ of the form \eqref{eqn_P_3_n_geq_3_at_least_m_ker_lambda} with $b\ne 0$, so that
	\begin{equation}
		\max\limits_{\{\langle s,s\rangle + \langle u,u\rangle + w^2=1\}}P_3\left(\left(\begin{smallmatrix}s\\ u\\ w\end{smallmatrix}\right)\right)\leq \frac{2}{3\sqrt{3}},\quad \left\|\frac{1}{\sqrt{1-\widehat{q}R}}\widehat{\gamma}\right\|^2=\frac{(3-R^2)^2}{3b^2R^4}\label{eqn_condition_kerm1_example}
	\end{equation}
hold.

	\begin{Ex}
		With the assumptions as above, let $P_3\left(\left(\begin{smallmatrix}s\\ u\\ w\end{smallmatrix}\right)\right)$ be as in \eqref{eqn_P_3_n_geq_3_at_least_m_ker_lambda}. Recall that $n\geq 3$ and $1\leq m\leq n-2$. We set 
			\begin{equation}
				R=1,\quad b=\frac{2}{\sqrt{3}},\quad \widehat{q}=0,\quad \widehat{\gamma}=\partial_{n-1},\quad C(s,u,u)=0,\quad C(u)=-\frac{2}{3\sqrt{3}}u_{n-1-m}^3.\label{eqn_values_kerm1_example}
			\end{equation}
		Recall that we have identified $\partial_{n-1}=\partial_{u_{n-1-m}}$ with the $(n-1)$-th Euclidean unit vector $(0,\ldots,1)^\mathrm{T}\in\mathbb{R}^{n-1}$. We claim that $P_3\left(\left(\begin{smallmatrix}s\\ u\\ w\end{smallmatrix}\right)\right)$ fulfils the two conditions in \eqref{eqn_condition_kerm1_example} and is thereby an example for $\dim\ker \boldsymbol\lambda|_{r=R}=m+1$ in any dimension $n\geq 3$. The second condition in \eqref{eqn_condition_kerm1_example} is a very simple calculation. For the first condition it is by restriction of $P_3\left(\left(\begin{smallmatrix}s\\ u\\ w\end{smallmatrix}\right)\right)$ to arbitrary $3$-dimensional subvector spaces of $\mathbb{R}^n\cong\mathbb{R}^m\times\mathbb{R}^{n-1-m}\times\mathbb{R}$ spanned by vectors of the form $(V_s^\mathrm{T},0,0)^\mathrm{T}$, $(0,V_u^\mathrm{T},0)^\mathrm{T}$, and $(0,0,1)^\mathrm{T}$ with $V_s\in\mathbb{R}^m$, $V_u\in\mathbb{R}^{n-1-m}$, sufficient to show that for all $c\in[0,1]$
			\begin{equation}
				\max\limits_{\{s^2+u^2+w^2=1\}}\widetilde{P_3}\left(\left(\begin{smallmatrix}s\\u\\w\end{smallmatrix}\right)\right) = \max\limits_{\{s^2+u^2+w^2=1\}}\left(-\frac{1}{\sqrt{3}}cs^2u -\frac{2}{3\sqrt{3}}c^3u^3 +s^2w+\frac{2}{\sqrt{3}}cuw^2\right)\leq \frac{2}{3\sqrt{3}}.\label{eqn_ker_m1_example_tocheck}
			\end{equation}
		In the above equation, $s,u,w$ denote linear coordinates on $\mathbb{R}^3$. Proceeding as in Example \ref{example_Rsqrt2_dimker_m} we find that the non-zero solutions of $\D \widetilde{P_3}=\frac{2}{\sqrt{3}}(s\,\D s + u\,\D u + w\,\D w)$ are given by
			\begin{equation*}
				\begin{array}{rl}
					\left\{
						\left(
							\begin{matrix}
								\pm\frac{\sqrt{2}}{\sqrt{3}}\\ 0\\ \frac{1}{\sqrt{3}}
							\end{matrix}
						\right)
					\right\} &\text{for }c=0,\\
					\left\{
						\left(
							\begin{matrix}
								0\\ \frac{1}{2}\\ -\frac{\sqrt{3}}{2}
							\end{matrix}
						\right)
					\right\}\cup
					\left\{
						\left.
							\left(
								\begin{matrix}
									\pm\sqrt{2\sqrt{3}t-4t^2}\\ \sqrt{3}t-1\\ t
								\end{matrix}
							\right)
						\ \right|
						t\in\left[0,\frac{\sqrt{3}}{2}\right]
					\right\} &\text{for }c=1,\\
					\left\{
						\left(
							\begin{matrix}
								0\\ \frac{1}{2c}\\ \pm\frac{\sqrt{c^2+2}}{2c}
							\end{matrix}
						\right),
						\left(
							\begin{matrix}
								\pm\frac{\sqrt{2}}{\sqrt{3}}\\ 0\\ \frac{1}{\sqrt{3}}
							\end{matrix}
						\right),
						\left(
							\begin{matrix}
								0\\ -\frac{1}{c^3}\\ 0
							\end{matrix}
						\right)
					\right\} &\text{for }c\in(0,1).
				\end{array}
			\end{equation*}
		In all of the above cases the Euclidean norm of the solutions is either exactly $1$ or bigger than $1$. The only thing that is not immediate when checking this is to show that $\frac{c^2+2}{4c^2}\leq 1$ for all $c\in(0,1)$, which follows by strict monotonicity. Hence, \eqref{eqn_ker_m1_example_tocheck} holds, showing that the maximal connected PSR manifold in standard form defined by $P_3$ fulfilling \eqref{eqn_P_3_n_geq_3_at_least_m_ker_lambda} and \eqref{eqn_values_kerm1_example} is indeed closed.
	\end{Ex}
	
	Lastly we summarise the results of this section in the following table of obtained limit geometries, respectively polynomials, in dependence of properties of the initial CCPSR manifold $\mathcal{H}$:
		\begin{equation*}
			\begin{array}{|c|c|c|c|c|}
				\hline \dim(\mathcal{H}) & b & m &\dim\ker\boldsymbol\lambda|_{r=R} & \text{eqn.}\\
				\hline \hline \geq 3 & \ne 0 & 0 & 0 & \eqref{eqn_dimgen_lambda_nokernel_hE1_limit}\\
				\hline \geq 3 & \ne 0 & 0 & 1 & \eqref{eqn_limit_poly_ker_0_1}\\
				\hline \geq 3 & \ne 0 & \in\{1,\ldots,n-2\} & m & \eqref{eqn_limit_poly_ker_m}\\
				\hline \geq 3 & \ne 0 & \in\{1,\ldots,n-2\} & m+1 & \eqref{eqn_limit_poly_m_plus_1}\\
				\hline
			\end{array}
		\end{equation*}

\subsection{$\dim(\mathcal{H})\geq 3$, $b=0$}\label{subsection_dimHgeq3_b0}
We proceed as in Section \ref{section_dimH_geq3_bnon0} with the difference that we now assume $b=0$. We have the following simpler versions of equations \eqref{eqn_dimgen_lambda_def}--\eqref{eqn_dimgen_Theta_www}:
	\begin{align}
		\boldsymbol\lambda&=\sum\limits_{i=1}^{n-1}(1-q_ir)(3-r^2)^2v_i^2,\label{eqn_dimgen_lambda_def_b0}\\
		\boldsymbol\chi&=0,\\
		\boldsymbol\mu&=3(3-9ar+r^2)\beta w^2,\label{eqn_dimgen_mu_def_b0}\\
		\boldsymbol\Theta_{v^3}&=(3-r^2)^3C(v)\sqrt{\beta},\label{eqn_dimgen_Theta_vvv_b0}\\
		\boldsymbol\Theta_{v^2w}&= \left(-r(2-3ar)(3-r^2)^2\langle v,v\rangle+(3-r^2)^3Q(v)\right)w \sqrt{\beta},\label{eqn_dimgen_Theta_vvw_b0}\\
		\boldsymbol\Theta_{vw^2}&=0,\label{eqn_dimgen_Theta_vww_b0}\\
		\boldsymbol\Theta_{w^3}&=\left((2-27a^2)r^3+27ar^2-18r+27a\right)w^3\beta\sqrt{\beta}.\label{eqn_dimgen_Theta_www_b0}
	\end{align}
As before, we use the convention $\boldsymbol\mu = \mu w^2$. Similarly to the 2-dimensional cases with $b=0$, will differentiate between the different possible values of $\dim\ker\boldsymbol\lambda|_{r=R}$. For this we will start with additional assumption $R<\sqrt{3}$. The the cases with $R=3$ will be treated separately.

For $R<\sqrt{3}$, it is clear that $\dim\ker\boldsymbol\lambda|_{r=R}\in\{0,\ldots,n-1\}$ depending on the values of the $q_i$. Without loss of generality assume that $q_1\geq\ldots\geq q_{n-1}$. Then
	\begin{equation}\label{eqn_dimH_geq3_wlog_q_i_order}
		\dim\ker\boldsymbol\lambda|_{r=R}=m\quad \Leftrightarrow\quad q_1=\ldots=q_m=\frac{1}{R},\ q_{m+1}<\frac{1}{R}.
	\end{equation}
The above equality is meant to include the case $m=n-1$ where $q_1=\ldots=q_m=\frac{1}{\sqrt{3}}$, so that the allowed values for $m$ are $0,\ldots,n-1$.

We start with $\dim\ker\boldsymbol\lambda|_{r=R}=0$. In this case, $R<\sqrt{3}$ is necessarily satisfied. In this case we can proceed exactly as in the analogous case for $b\ne 0$\label{ref to subsubsection or paragraph to be added later here} from equation \eqref{eqn_dimgen_lambda_nokernel_diag} onwards and obtain, just as in equation \eqref{eqn_dimgen_lambda_nokernel_hE1_limit}, that the limit polynomial is given by
	\begin{equation}\label{eqn_limit_poly_dimHgeq3_b0_Rlesssqrt3_m0}
		\lim\limits_{r\nearrow R}\left(B\times(\mathbbm{1}_\mathbb{R}/\sqrt{\mu})\right)^\ast 	h_{E}\left(\left(\begin{smallmatrix} x\\ v\\ w\end{smallmatrix}\right)\right)=x^3-x\left(\langle v,v\rangle +w^2\right) -\frac{1}{\sqrt{3}}\langle v,v\rangle w -\frac{2}{3\sqrt{3}} w^3,
	\end{equation}
cf. \eqref{eqn_dimgen_hE1_prelimit} for the definition of $h_E$. As explained after equation \eqref{eqn_dimgen_lambda_nokernel_hE1_limit} the corresponding CCPSR manifold is homogeneous and, hence, this case is in accordance with Theorem \ref{thm_limit_behaviour_ccpsr}.

Next we will study the cases with $\dim\ker\boldsymbol\lambda|_{r=R}=m$, $1\leq m\leq n-1$, and the additional assumption that $R<\sqrt{3}$. Assume without loss of generality that $q_1=\ldots=q_m=\frac{1}{R}$, cf. \eqref{eqn_dimH_geq3_wlog_q_i_order}. In the following, assume for the calculations that $m<n-1$. It will be clear how the arguments work in the case $m=n-1$. Using the notation \eqref{eqn_s_u_def} to split up our $v$-coordinates, we get
	\begin{equation}\label{eqn_dimH_geq3_C_splitting_b0}
		C(v)=\widetilde{C}(s)+\langle u,\widetilde{\mathcal{F}}(s)\rangle + \langle s,\widehat{\mathcal{F}}(u)\rangle + \widehat{C}(u),
	\end{equation}
where $\widetilde{C}:\mathbb{R}^m\to\mathbb{R}$ and $\widehat{C}:\mathbb{R}^{n-1-m}\to\mathbb{R}$ are homogeneous cubic polynomials, and $\widetilde{F}:\mathbb{R}^m\to\mathbb{R}^{n-1-m}$ and $\widehat{F}:\mathbb{R}^{n-1-m}\to\mathbb{R}^m$ are entry-wise homogeneous quadratic polynomials. We further abbreviate, similar to \eqref{eqn_dimgen_q_def},
	\begin{equation}\label{eqn_q_split}
		q=\left(\begin{array}{ccc|ccc} q_1 & & & & & \\ & \ddots & & & & \\ & & q_m & & & \\ \hline
			& & & q_{m+1} & & \\ & & & & \ddots & \\ & & & & & q_{n-1}\end{array}\right)
		=\left(\begin{array}{c|c} \widetilde{q} & \\ \hline & \overset{}{\widehat{q}}\end{array}\right).
	\end{equation}
Note that in the case we are presently studying, $\widetilde{q}$ is simply the identity matrix $\mathbbm{1}_m$ multiplied with $\frac{1}{\sqrt{3}}$. Let
	\begin{equation}\label{eqn_B_def_b0_cases}
		B(r):=(3-r^2)^{-1}(1-qr)^{-\frac{1}{2}}=(3-r^2)^{-1}\left(\begin{array}{c|c} (1-\frac{r}{R})^{-\frac{1}{2}}\underset{}{\mathbbm{1}_m} & \\ \hline & \overset{}{(1-\widehat{q}r)^{-\frac{1}{2}}}\end{array}\right).
	\end{equation}
We obtain
	\begin{align}
		\left(B(r)\times(\mathbbm{1}_\mathbb{R}/\sqrt{\mu})\right)^\ast\boldsymbol\Theta_{v^3}
		&=\left(\widetilde{C}\left((1-\tfrac{r}{R})^{-\frac{1}{2}}s\right) + \left\langle(1-\widehat{q}r)^{-\frac{1}{2}}u,\widetilde{\mathcal{F}}\left((1-\tfrac{r}{R})^{-\frac{1}{2}}s\right)\right\rangle\right.\notag\\ 
		&\quad + \left.\left\langle(1-\tfrac{r}{R})^{-\frac{1}{2}}s,\widehat{\mathcal{F}}\left((1-\widehat{q}r)^{-\frac{1}{2}}u\right)\right\rangle + \widehat{C}\left((1-\widehat{q}r)^{-\frac{1}{2}}u\right)\right)w\sqrt{\beta}.	\label{eqn_pullback_Theta_vvv_Rlesssqrt3_b0}
	\end{align}
Since $R<\sqrt{3}$, $\beta$ has a simple zero in $r=R$, and it is clear that $1-\frac{r}{R}$ also has a simple zero in $r=R$. Hence, using our assumption that the initial $\mathcal{H}$ is closed in its ambient space and the homogeneity of $\widetilde{C}$, $\widehat{C}$ and $\widetilde{\mathcal{F}}$, $\widehat{\mathcal{F}}$ of degree $3$ and $2$, respectively, this shows with Theorem \ref{thm_Cn} that $\widetilde{C}$ and $\widetilde{\mathcal{F}}$ must identically vanish. Furthermore, $q_i<\frac{1}{\sqrt{3}}$ for all $i>m$ implies that
	\begin{equation*}
		\lim\limits_{r\nearrow R} \widehat{C}\left((1-\widehat{q}r)^{-\frac{1}{2}}u\right)w\sqrt{\beta}=0.
	\end{equation*}
Thus we need to study the limit
	\begin{align}
		\lim\limits_{r\nearrow R} \left\langle(1-\tfrac{r}{R})^{-\frac{1}{2}}s,\widehat{\mathcal{F}}\left((1-\widehat{q}r)^{-\frac{1}{2}}u\right)\right\rangle w\sqrt{\beta}&=\lim\limits_{r\nearrow R}\frac{\sqrt{\beta}}{\sqrt{1-\tfrac{r}{R}}} \left\langle s, \widehat{\mathcal{F}}\left((1-\widehat{q}R)^{-\frac{1}{2}}u\right)\right\rangle \notag \\
		&=\sqrt{3-R^2}\left\langle s, \widehat{\mathcal{F}}\left((1-\widehat{q}R)^{-\frac{1}{2}}u\right)\right\rangle,\label{eqn_lim_F_Rlesssqrt3}
	\end{align}
which can be easily checked by observing that $\frac{\beta}{1-\frac{r}{R}}=1+\frac{r}{R}+r^2\left(\frac{1}{R^2}-1\right)$. The properties and allowed values of $\widehat{\mathcal{F}}$ are controlled by the terms $F_i$ in Theorem \ref{thm_limit_behaviour_ccpsr} \eqref{thm_limit_behaviour_ccpsr_iii}, cf. Section \ref{sect_Fi_calcs} for the proof of these claims. Next we will study the limit of the pullback of $\boldsymbol\Theta_{v^2 w}$. We calculate
	\begin{align}\label{eqn_pullback_Theta_vvw_limit_Rlesssqrt3}
		\lim\limits_{r\nearrow R}\left(B(r)\times(\mathbbm{1}_\mathbb{R}/\sqrt{\mu})\right)^\ast\boldsymbol\Theta_{v^2w} = \left(\frac{2}{\sqrt{3}}\langle s,s\rangle -\frac{1}{\sqrt{3}}\langle u,u\rangle\right)w.
	\end{align}
Together with equation \eqref{eqn_Theta_w3_limit_R_less_sqrt3}, which holds in our present case exactly as stated, we have shown that for $1\leq m<n-1$ the limit polynomial is given by
	\begin{align}
		\lim\limits_{r\nearrow R}\left(B(r)\times(\mathbbm{1}_\mathbb{R}/\sqrt{\mu})\right)^\ast 	h_{E}\left(\left(\begin{smallmatrix} x\\ v\\ w\end{smallmatrix}\right)\right)&=x^3-x\left(\langle v,v\rangle +w^2\right) + \sqrt{3-R^2}\left\langle s, \widehat{\mathcal{F}}\left((1-\widehat{q}R)^{-\frac{1}{2}}u\right)\right\rangle\notag\\
		&\quad + \left(\frac{2}{\sqrt{3}}\langle s,s\rangle -\frac{1}{\sqrt{3}}\langle u,u\rangle\right)w -\frac{2}{3\sqrt{3}} w^3.\label{eqn_limit_poly_dimHgeq3_b0_Rlesssqrt3_mlessn1}
	\end{align}
For $m=n-1$ observe that $C=\widetilde{C}$ \eqref{eqn_dimH_geq3_C_splitting_b0}, and thus it is easy to see with our above calculations that in this case the limit polynomial is of the form
	\begin{equation}
		\lim\limits_{r\nearrow R}\left(B(r)\times(\mathbbm{1}_\mathbb{R}/\sqrt{\mu})\right)^\ast 	h_{E}\left(\left(\begin{smallmatrix} x\\ v\\ w\end{smallmatrix}\right)\right)=x^3-x\left(\langle v,v\rangle +w^2\right) + \frac{2}{\sqrt{3}}\langle v,v\rangle w -\frac{2}{3\sqrt{3}} w^3.\label{eqn_limit_poly_dimHgeq3_b0_Rlesssqrt3_mn1}
	\end{equation}

Next, we need to consider the case $R=\sqrt{3}$. In comparison to the case $R<\sqrt{3}$ we now always have $\dim\ker\boldsymbol\lambda|_{r=R}=n-1$. So instead we will do a case-by-case study of the different possible of
	\begin{equation}\label{eqn_m_def_R3_case}
		\dim\ker\left(\boldsymbol\lambda\cdot(3-r^2)^{-2}\right)|_{r=R}=m,\quad 0\leq m\leq n-1.
	\end{equation}
We adopt the conventions for the used in the case $R<\sqrt{3}$, that is $q_1\geq\ldots\geq q_m$, \eqref{eqn_dimH_geq3_C_splitting_b0}, \eqref{eqn_q_split}, \eqref{eqn_B_def_b0_cases}, and the coordinate relabelling \eqref{eqn_s_u_def}. Similar to equation \eqref{eqn_dimH_geq3_wlog_q_i_order}, we see that $\dim\ker\left(\boldsymbol\lambda\cdot(3-r^2)^{-2}\right)|_{r=R}=m$ if and only if $q_1=\ldots=q_m=\frac{1}{\sqrt{3}},\ q_{m+1}<\frac{1}{\sqrt{3}}$. We start with $m=0$ and calculate
	\begin{align*}
		\lim\limits_{r\nearrow R}\left(B\times(\mathbbm{1}_\mathbb{R}/\sqrt{\mu})\right)^\ast\boldsymbol\Theta_{v^3}
		&=\lim\limits_{r\nearrow R}C\left((1-qr)^{-\frac{1}{2}}v\right)\sqrt{\beta}=0,\\
		\lim\limits_{r\nearrow R}\left(B\times(\mathbbm{1}_\mathbb{R}/\sqrt{\mu})\right)^\ast\boldsymbol\Theta_{v^2w}
		&=-\frac{2}{\sqrt{3}}\langle v,v\rangle w.
	\end{align*}
We further check that
	\begin{equation}\label{eqn_pullback_Theta_www_Rsqrt3_constant}
		\left(B\times(\mathbbm{1}_\mathbb{R}/\sqrt{\mu})\right)^\ast\boldsymbol\Theta_{w^3}=\frac{\boldsymbol\Theta_{w^3}}{\mu\sqrt{\mu}}=\frac{2}{3\sqrt{3}}
	\end{equation}
independent of $r$. Since $\boldsymbol\Theta_{vw^2}=0$, it follows that the limit polynomial is given by
	\begin{equation}
		\lim\limits_{r\nearrow R}\left(B\times(\mathbbm{1}_\mathbb{R}/\sqrt{\mu})\right)^\ast 	h_{E}\left(\left(\begin{smallmatrix} x\\ v\\ w\end{smallmatrix}\right)\right)=x^3-x\left(\langle v,v\rangle +w^2\right) - \frac{2}{\sqrt{3}}\langle v,v\rangle w + \frac{2}{3\sqrt{3}} w^3.\label{eqn_limit_poly_dimHgeq3_b0_Rsqrt3_m0}
	\end{equation}
Recall that the corresponding CCPSR manifold is homogeneous, cf. \cite[Prop.\,6.9]{L1}, and is in accordance with Theorem \ref{thm_limit_behaviour_ccpsr}.

Next consider $1\leq m < n-1$ in \eqref{eqn_m_def_R3_case}. Whit our convention for the $q_i$ this is the case if and only if $q_1=\ldots=q_m=0$, $q_j<\frac{1}{\sqrt{3}}$ for all $m+1\leq j\leq n-1$. Note that equation \eqref{eqn_pullback_Theta_www_Rsqrt3_constant} is independent of $m$. Also, equation \eqref{eqn_pullback_Theta_vvv_Rlesssqrt3_b0} still holds when setting $R=\sqrt{3}$. The difference to the case $R<\sqrt{3}$ is now that $\beta$ has now a double zero in $r=R=\sqrt{3}$. Hence, $\sqrt{\beta}$ has a simple zero in $r=\sqrt{3}$, which implies that $\widetilde{C}$ must identically vanish by Theorem \ref{thm_Cn} since $\mathcal{H}$ was assumed to be closed in its ambient space, and we obtain
	\begin{equation*}
		\lim\limits_{r\nearrow R}\left(\left\langle(1-\tfrac{r}{R})^{-\frac{1}{2}}s,\widehat{\mathcal{F}}\left((1-\widehat{q}r)^{-\frac{1}{2}}u\right)\right\rangle + \widehat{C}\left((1-\widehat{q}r)^{-\frac{1}{2}}u\right)\right)\sqrt{\beta}=0.
	\end{equation*}
Next, analogous to equation \eqref{eqn_lim_F_Rlesssqrt3}, we check that
	\begin{equation*}
		\lim\limits_{r\nearrow R}\left\langle(1-\widehat{q}r)^{-\frac{1}{2}}u,\widetilde{\mathcal{F}}\left((1-\tfrac{r}{R})^{-\frac{1}{2}}s\right)\right\rangle\sqrt{\beta}=\sqrt{3}\left\langle u, (1-\sqrt{3}\widehat{q})^{-\frac{1}{2}}\widetilde{\mathcal{F}}(s)\right\rangle,
	\end{equation*}
and we calculate
	\begin{equation*}
		\lim\limits_{r\nearrow R}\left(B\times(\mathbbm{1}_\mathbb{R}/\sqrt{\mu})\right)^\ast\boldsymbol\Theta_{v^2w} = \left(\frac{1}{\sqrt{3}}\langle s,s\rangle -\frac{2}{\sqrt{3}}\langle u,u\rangle\right)w.
	\end{equation*}
Observe the switched coefficients of $s$ and $u$ when comparing the upper calculation to the limit in \eqref{eqn_pullback_Theta_vvw_limit_Rlesssqrt3}. Summarising, we have shown that in the case $R=\sqrt{3}$, $1\leq m<n-1$ the limit polynomial is given by
	\begin{align}
		\lim\limits_{r\nearrow R}\left(B\times(\mathbbm{1}_\mathbb{R}/\sqrt{\mu})\right)^\ast 	h_{E}\left(\left(\begin{smallmatrix} x\\ v\\ w\end{smallmatrix}\right)\right)&=x^3-x\left(\langle s,s\rangle + \langle u,u\rangle +w^2\right) + \sqrt{3}\left\langle u, (1-\sqrt{3}\widehat{q})^{-\frac{1}{2}}\widetilde{\mathcal{F}}(s)\right\rangle\notag \\
		&\quad + \left(\frac{1}{\sqrt{3}}\langle s,s\rangle -\frac{2}{\sqrt{3}}\langle u,u\rangle\right) w + \frac{2}{3\sqrt{3}} w^3.\label{eqn_limit_poly_dimHgeq3_b0_Rsqrt3_mlessn1}
	\end{align}

Lastly we need to consider the case $R=\sqrt{3}$, $m=n-1$. Using our above calculations with $C=\widetilde{C}$, cf. equation \eqref{eqn_dimH_geq3_C_splitting_b0}, makes it easy to see that the limit polynomial is given by
	\begin{equation}
		\lim\limits_{r\nearrow R}\left(B\times(\mathbbm{1}_\mathbb{R}/\sqrt{\mu})\right)^\ast 	h_{E}\left(\left(\begin{smallmatrix} x\\ v\\ w\end{smallmatrix}\right)\right)=x^3-x\left(\langle v,v\rangle +w^2\right) +\frac{1}{\sqrt{3}}\langle v,v\rangle w +\frac{2}{3\sqrt{3}} w^3.\label{eqn_limit_poly_dimHgeq3_b0_Rsqrt3_mn1}
	\end{equation}

Now, if not already done so, it remains to check that listed possible limit polynomials in this part of the proof, namely
	\begin{equation*}
		\begin{array}{|c|c|c|c|c|}
			\hline \dim(\mathcal{H}) & b & R & m & \text{eqn.}\\
			\hline \hline \geq 3 & 0 & <\sqrt{3} & 0 & \eqref{eqn_limit_poly_dimHgeq3_b0_Rlesssqrt3_m0}\\
			\hline \geq 3 & 0 & <\sqrt{3} & \in\{1,\ldots,n-2\} & \eqref{eqn_limit_poly_dimHgeq3_b0_Rlesssqrt3_mlessn1}\\
			\hline \geq 3 & 0 & <\sqrt{3} & n-1 & \eqref{eqn_limit_poly_dimHgeq3_b0_Rlesssqrt3_mn1}\\
			\hline \geq 3 & 0 & \sqrt{3} & 0 & \eqref{eqn_limit_poly_dimHgeq3_b0_Rsqrt3_m0}\\
			\hline \geq 3 & 0 & \sqrt{3} & \in\{1,\ldots,n-2\} & \eqref{eqn_limit_poly_dimHgeq3_b0_Rsqrt3_mlessn1}\\
			\hline \geq 3 & 0 & \sqrt{3} & n-1 & \eqref{eqn_limit_poly_dimHgeq3_b0_Rsqrt3_mn1}\\
			\hline
		\end{array}
	\end{equation*}
are in accordance with Thm. \ref{thm_limit_behaviour_ccpsr} \eqref{thm_limit_behaviour_ccpsr_iii}, which is easily seen to be true. The only thing that one must be aware of is that Lemma \ref{lem_convergence_limit_in_Cn} implies that the term $\sqrt{3-R^2}\left\langle s, \widehat{\mathcal{F}}\left((1-\widehat{q}R)^{-\frac{1}{2}}u\right)\right\rangle$ in equation \eqref{eqn_limit_poly_dimHgeq3_b0_Rlesssqrt3_mlessn1} when written using the $F_i$-notation in Thm. \ref{thm_limit_behaviour_ccpsr} \eqref{thm_limit_behaviour_ccpsr_iii} must uphold the corresponding eigenvalue bound that we will show in the upcoming section.

\subsection{Proof of the eigenvalue bounds in Thm. \ref{thm_limit_behaviour_ccpsr} \eqref{thm_limit_behaviour_ccpsr_iii}}\label{sect_Fi_calcs}
	We know by Lemma \ref{lem_convergence_limit_in_Cn} That the possible limit polynomials of the form \eqref{eqn_dim_geq_3_limits} obtained in Sections \ref{subsection_dimHgeq3_b0} and \ref{subsection_dimHgeq3_b0} are, in fact, contained in $\mathcal{C}_n$ \eqref{eqn_Cn_gen_set}. It remains to show that for all choices of symmetric matrices $F_i\in \mathrm{Sym}((n-1-m)\times (n-1-m),\mathbb{R})$, $1\leq i\leq m$, such that for all $c\in\mathbb{R}^m$ with $\|c\|=1$ the eigenvalues of $\sum\limits_{i=1}^{m}c_i F_i$ are contained in $[-1,1]$, $\overline{h}$ as in \eqref{eqn_dim_geq_3_limits} is contained in $\mathcal{C}_n$ and that each such polynomial can be realized as a limit polynomial of some given polynomial $h\in\mathcal{C}_n$ with corresponding CCPSR manifold in standard form $\mathcal{H}$. Suppose that we are given symmetric matrices $F_1,\ldots,F_{m}$ that fulfil this eigenvalue condition. In order to show that the corresponding polynomial $\overline{h}$ is contained in $\mathcal{C}_n$, we need to show that
		\begin{equation}
			\max\limits_{\{\langle s,s\rangle + \langle u,u\rangle + w^2=1\}}\left(\sum\limits_{i=1}^{m}  s_i\langle u,F_i u\rangle + \left(\frac{2}{\sqrt{3}}\langle s,s\rangle - \frac{1}{\sqrt{3}}\langle u,u\rangle\right)w - \frac{2}{3\sqrt{3}} w^3\right)\leq \frac{2}{3\sqrt{3}}.\label{eqn_Fi_estimate_Cn}
		\end{equation}
	In order to prove the above estimate \eqref{eqn_Fi_estimate_Cn}, it is sufficient to show that it holds for all unit vectors $(s,u,w)^\mathrm{T}$ restricted to all $3$-dimensional subvector spaces of $\mathbb{R}^n$ of the form
		\begin{equation*}
			\left.\left\{\left(\begin{smallmatrix} t_s V_s\\ t_u V_u\\ t_w V_w\end{smallmatrix}\right)\in\mathbb{R}^n\ \right|\ V_s\in\mathbb{R}^m,\ V_u\in\mathbb{R}^{n-1-m},\ V_w\in\mathbb{R},\ \langle V_s,V_s\rangle + \langle V_u,V_u\rangle + V_w^2=1,\ t_s,t_u,t_w\in\mathbb{R}\right\}.
		\end{equation*}
	This follows from the fact that the weighted average $\sum\limits_{i=1}^{m}c_i F_i$ for $\sum\limits_{i=1}^{m} c_i^2=1$ will always have eigenvalues contained in $[-1,1]$. Hence, by the eigenvalue condition this reduces the problem to the case $n=3$ and showing that for all $F\in[-1,1]$ it holds that
		\begin{equation}
			\max\limits_{\{s^2 + u^2 + w^2 =1\}}\left( Fsu^2 + \left(\frac{2}{\sqrt{3}}s^2 - \frac{1}{\sqrt{3}}u^2\right)w - \frac{2}{3\sqrt{3}}w^3\right)\leq \frac{2}{3\sqrt{3}}.\label{eqn_Fi_reduced}
		\end{equation}
	In fact, we already know that $\max\limits_{\{s^2 + u^2 + w^2 =1\}}\left( Fsu^2 + \left(\frac{2}{\sqrt{3}}s^2 - \frac{1}{\sqrt{3}}u^2\right)w - \frac{2}{3\sqrt{3}}w^3\right)\geq \frac{2}{3\sqrt{3}}$, which follows from the $w^3$-coefficient being of absolute value $\frac{2}{3\sqrt{3}}$. Thus, in order to prove \eqref{eqn_Fi_reduced}, it suffices to solve
		\begin{equation}
			\D \left(Fs^2u + \left(-\frac{1}{\sqrt{3}}s^2 + \frac{2}{\sqrt{3}}u^2\right)w - \frac{2}{3\sqrt{3}}w^3\right)-\frac{2}{\sqrt{3}}(s\,\D s + u\,\D u + w\,\D w)=0\label{eqn_F_reduced_max_points}
		\end{equation}
	for $ s,u,w \in\mathbb{R}$ and showing that every nonzero solution has Euclidean norm at least $1$. Since a sign-flip in $u$ is equivalent to a sign-flip in $F$, we can without loss of generality assume that $F\in[0,1]$. The nonzero solutions of \eqref{eqn_F_reduced_max_points} are given by
		\begin{equation}
			(s,u,w)\in\left\{\left(\pm\frac{\sqrt{3}}{2},0,\frac{1}{2}\right), \left(\frac{1}{\sqrt{3}F},\pm\frac{\sqrt{2}}{\sqrt{3}F},0\right),(0,0,-1)\right\}\label{eqn_sol_set_F}
		\end{equation}
	for all $F\in (0,1)\setminus 0$, for $F=1$ additionally the family of solutions
		\begin{equation*}
			(s,u,w)\in\left\{\left.\left(u,\pm\sqrt{-4u^2+2\sqrt{3}u},-1+\sqrt{3}u\right)\ \right|\ u\in \left(0,\frac{\sqrt{3}}{2}\right) \right\},
		\end{equation*}
	and for $F=0$ we have the solution set
		\begin{equation*}
			(s,u,w)\in\left\{\left(\pm\frac{\sqrt{3}}{2},0,\frac{1}{2}\right),(0,0,-1)\right\}.
		\end{equation*}
	It is now an easy calculation to see that all of the above solutions have Euclidean norm exactly or greater than $1$, showing that \eqref{eqn_Fi_reduced} holds as an equality and, hence, proving \eqref{eqn_Fi_estimate_Cn}.

	On the other hand, suppose that there exists a Euclidean unit vector $c=(c_1,\ldots,c_{m})^\mathrm{T}\in\mathbb{R}^m$, such that $\sum\limits_{i=1}^{m} c_i F_i$ has an eigenvalue $F\in\mathbb{R}$ of absolute value greater than $1$. After restricting $s$ to the corresponding eigenspace and $u$ to $\mathbb{R}\cdot c$, we abuse notation slightly and need to show that \eqref{eqn_Fi_reduced} is false. In fact, one can check that the point $(s,u,w)= \left(\frac{1}{\sqrt{3}F},\pm\frac{\sqrt{2}}{\sqrt{3}F},0\right)$ as in \eqref{eqn_sol_set_F} is still a nonzero solution of \eqref{eqn_F_reduced_max_points}, but it is of Euclidean norm $\frac{1}{|F|}$ is smaller than one. Hence, \eqref{eqn_Fi_reduced} and consequently \eqref{eqn_Fi_estimate_Cn} do not hold, showing that the existence of such a vector $c$ violates the initial assumption that $\mathcal{H}$ is a CCPSR manifold.

	Note that at this point we have shown that every maximal connected PSR manifold in standard form contained in the level set of a polynomial of the form \eqref{eqn_dim_geq_3_limits} is, in fact, closed in its ambient space and thus a CCPSR manifold.

	Lastly we need to prove the claim that every CCPSR manifold of dimension at least $3$ in standard form contained in the level set of a polynomial of the form \eqref{eqn_dim_geq_3_limits} can be realised as a limit geometry of a CCPSR manifold, is singular at infinity, and has continuous symmetry group of dimension at least $1$. For dimensions $1$ and $2$ this holds as the only possible limit geometries are already homogeneous spaces, which are in particular singular at infinity \cite[Prop.\,4.6]{L2}. Hence, we only need to consider dimensions at least $3$. Recall that we have already shown that every possible maximal connected PSR manifold in standard form $\overline{\mathcal{H}}\subset\{\overline{h}=1\}$ with $\overline{h}$ of the form Thm. \ref{thm_limit_behaviour_ccpsr} \eqref{thm_limit_behaviour_ccpsr_iii} is in fact a CCPSR manifold, and the $w^3$-coefficient being of absolute value $\frac{2}{3\sqrt{3}}$ implies with Lemma \ref{lem_sing_at_infty_iff_max_cond} that $\overline{\mathcal{H}}$ is indeed singular at infinity. Thus, it suffices to show  that each such $\overline{\mathcal{H}}$ has continuous symmetry group of dimension at least $1$ which corresponds to moving the reference point for the standard form in $w$-direction, that is along curves of the form \eqref{eqn_p_xvw}. Geometrically this means that the possible limit polynomials in Thm. \eqref{thm_limit_behaviour_ccpsr_iii} and corresponding CCPSR manifolds in standard form are stable under taking their limit geometry along said curves. In order to prove this claim, it follows from Lemma \ref{lem_first_var_P3_sym_grps} that it suffices to show for
		\begin{equation*}
			P_3\left(\left(\begin{smallmatrix}s\\u\\w\end{smallmatrix}\right)\right)= \sum\limits_{i=1}^{m}  s_i\langle u,F_i u\rangle + \left(\frac{2}{\sqrt{3}}\langle s,s\rangle - \frac{1}{\sqrt{3}}\langle u,u\rangle\right)w - \frac{2}{3\sqrt{3}} w^3
		\end{equation*}
	corresponding to Thm. \ref{thm_limit_behaviour_ccpsr} \eqref{thm_limit_behaviour_ccpsr_iii} that $\delta P_3\left(\left(\begin{smallmatrix}s\\u\\w\end{smallmatrix}\right)\right)(\partial_w)=0$ for some choice of $L:\mathbb{R}^n\to\mathfrak{so}(n)$. After setting $L=0$, this is a slightly tedious but not difficult calculation and in fact turns out to be true.

	This finishes the proof of Theorem \ref{thm_limit_behaviour_ccpsr}.\hfill\qed

\section{Proof of Theorem \ref{thm_generic_limit_behaviour}}
	First note that the CCPSR manifold in standard form $\overline{\mathcal{H}}\subset\{\overline{h}=1\}$, $\overline{h}$ as in \eqref{eqn_h_lie_grp_generic_limit}, is in fact isomorphic to $\mathbb{R}_{>0}\ltimes\mathbb{R}^{n-1}$, which is slightly easier to check in different linear coordinates, cf. the proof of \cite[Prop.\,6.9]{L1}. For the $2$-dimensional case see \cite[Thm.\,1\,b)]{CDL} and \cite[Ex.\,3.2]{L2}.
	
	As for the proof of our main theorem we will proceed by first considering CCPSR surfaces and then CCPSR manifolds of dimension at least $3$. For $\dim(\mathcal{H})=2$, a CCPSR surface $\mathcal{H}$ can have as limit geometry either Thm. \ref{thm_limit_behaviour_ccpsr} \eqref{thm_limit_behaviour_ccpsr_ii_a} or Thm. \ref{thm_limit_behaviour_ccpsr} \eqref{thm_limit_behaviour_ccpsr_ii_b}. We need to show that $\mathcal{H}$ having regular boundary behaviour automatically excludes the second case. Assume that $\mathcal{H}\subset\{h=x^3-x(v^2+w^2)+cv^3+qv^2w+bvw^2+aw^3\}$ is in standard form. We have in the proof of the main theorem, cf. Sections \ref{subsect_dim2_b_ne_0} and \ref{sect_b_0}, seen that $\mathcal{H}$ can have Thm. \ref{thm_limit_behaviour_ccpsr} \eqref{thm_limit_behaviour_ccpsr_ii_b} as limit geometry if either $b\ne 0$ and \eqref{eqn_dim2_lambda_0_limit} holds, cf. equation \eqref{eqn_dim2_bne0_lambda0}, or $b=0$, $R=\sqrt{3}$, and $q<\frac{1}{\sqrt{3}}$, cf. equation \eqref{eqn_dim2_b0_Rsqrt3_qlessmax}. In the latter case $\mathcal{H}$ is singular at infinity by Lemma \ref{lem_beta_zeros} and, hence, in particular does not have regular boundary behaviour. In the former case recall that we have proven that $q=\frac{-3b^2R^4+(3-R^2)^2}{R(3-R^2)^2}$ and $c=\frac{bR^2 (b^2R^4-(3-R^2)^2)}{(3-R^2)^3}$ must hold, cf. equations \eqref{eqn_dim2_q_lamda_zero} and \eqref{eqn_dim2_c_bR_dep}. Also recall that $a=\frac{R^2-1}{R^3}$. We will now show that for all $b\ne 0$, $R\in\left(\frac{\sqrt{3}}{2},\sqrt{3}\right)$, such that $P_3\left(\left(\begin{smallmatrix}v\\w\end{smallmatrix}\right)\right)=cv^3+qv^2w+bvw^2+aw^3$ fulfils $\max\limits_{\{v^2+w^2=1\}}P_3\left(\left(\begin{smallmatrix}v\\ w\end{smallmatrix}\right)\right)\leq \frac{2}{3\sqrt{3}}$, equality of these two expression must in fact hold, showing that the corresponding CCPSR surface in standard form is singular at infinity and, hence, does not have regular boundary behaviour. The maybe easiest way to see that is, instead of trying to solve $\D P_3=\frac{2}{\sqrt{3}}(v\,\D v+ w\,\D w)$ and checking that one of the solutions has Euclidean norm $1$, to rotate the coordinates $\left(\begin{smallmatrix}v\\w\end{smallmatrix}\right)$ via
		\begin{equation}
			M(f):=\frac{1}{\sqrt{1+f^2}}\left(\begin{array}{rr} f & -1\\ 1 & f\end{array}\right)\in \mathrm{O}(2),\quad f\in\mathbb{R},\label{eqn_Mf_rotation}
		\end{equation}
	such that the $vw^2$-term in $P_3$ vanishes. Recall that we have in this case assumed without loss of generality $b>0$ which can always be achieved by a sign-flip in $v$. One finds with the help of a computer algebra system like Maple that with $f=\frac{3-R^2}{bR^3}$ that
		\begin{align*}
			\widetilde{P_3}\left(\left(\begin{smallmatrix}v\\w\end{smallmatrix}\right)\right):=\left(M(f)^*P_3\right)\left(\left(\begin{smallmatrix}v\\w\end{smallmatrix}\right)\right)&= \frac{\sqrt{9-6R^2+R^4+b^2R^6}}{(3-R^2)R}v^2w\\
			&\quad -\frac{(9-15R^2+7R^4+(b^2-1)R^6)\sqrt{9-6R^2+R^4+b^2R^6}}{(3-R^2)^3R^3}w^3.
		\end{align*}
	Now solving $\D \widetilde{P_3}-\frac{2}{\sqrt{3}}(v\,\D v+ w\,\D w)$ works well in the sense that a computer algebra system yields a solution more readable to human eyes, and we obtain the solution set consisting of three non-zero solutions
		\begin{align*}
			(v,w)&\in\left\{ \left(0,-\frac{2(3-R^2)^3R^3}{3\sqrt{3}(9-15R^2+7R^4+(b^2-1)R^6)\sqrt{9-6R^2+R^4+b^2R^6}}\right),\right.\\
				&\quad\left.\left(\pm\frac{\sqrt{27-27R^2+9R^4+(3b^2-1)R^6}}{\sqrt{3}\sqrt{9-6R^2+R^4+b^2R^6}},\frac{(3-R^2)R}{\sqrt{3}\sqrt{9-6R^2+R^4+b^2R^6}}\right)\right\}
		\end{align*}
	and, trivially, the origin.
	One quickly checks that the second pair of solutions above each have Euclidean norm $1$ and, hence, $\mathcal{H}$ is singular at infinity as claimed. We conveniently obtain a sharp upper limit on $b$, respectively $|b|$ if we had not assumed $b>0$, depending on the choice of $R\in\left(\frac{\sqrt{3}}{2},\sqrt{3}\right)$ by analysing when the first solution has Euclidean norm at least $1$. Any lower number would violate the maximality condition on $\widetilde{P_3}$ by its homogeneity of degree $3$. We find that the first solution has Euclidean norm at least one if and only if
		\begin{equation*}
			0<b\leq \frac{(3-R^2)\sqrt{12R^2-9}}{3R^3},
		\end{equation*}
	this can easily be verified with the help of a computer algebra system like Maple. A similar result can be obtained if one does not restrict $b$ to be positive. One also verifies that in fact for all $b\in\left(0,\frac{(3-R^2)\sqrt{12R^2-9}}{3R^3}\right]$, each of the three solutions is real and the solutions are always distinct. For the extremal case $b=\frac{(3-R^2)\sqrt{12R^2-9}}{3R^3}$ we see that the first solution is exactly of Euclidean norm one. Hence, $\mathcal{H}$ is singular at infinity along two distinct rays in $\partial(\mathbb{R}_{>0}\cdot\mathcal{H})$ if $b<\frac{(3-R^2)\sqrt{12R^2-9}}{3R^3}$ and along three distinct rays if $b= \frac{(3-R^2)\sqrt{12R^2-9}}{3R^3}$. In the latter case we have
		\begin{equation*}
			\widetilde{P_3}\left(\left(\begin{smallmatrix}v\\w\end{smallmatrix}\right)\right)=\frac{2}{\sqrt{3}}v^2w - \frac{2}{3\sqrt{3}}w^3.
		\end{equation*}
	Hence, $\mathcal{H}$ is in this case already a homogeneous space, cf. \cite[Thm.\,1\,a)]{CDL} with \cite[Ex.\,3.2]{L2}, and the same holds for its considered limit geometry.

	Next we consider $\dim(\mathcal{H})\geq 3$. We first need to check which kind of initial CCPSR manifolds in standard form allow for a limit geometry that is not isomorphic to $\mathbb{R}_{>0}\ltimes \mathbb{R}^{n-1}$. For $b=0$ studied in Section \ref{subsection_dimHgeq3_b0}, we can ignore the case $R=\sqrt{3}$ entirely since then $\mathcal{H}$ is singular at infinity by Lemma \ref{lem_beta_zeros}. For $R<\sqrt{3}$ the limit geometry in $w$-direction might only not be isomorphic to $\mathbb{R}_{>0}\ltimes\mathbb{R}^{n-1}$ or, equivalently, the limit polynomial can only not be equivalent to
		\begin{equation}
			\overline{h}=x^3-x(\langle v,v\rangle +w^2) -\frac{1}{\sqrt{3}}\langle v,v\rangle w -\frac{2}{3\sqrt{3}} w^3\label{eqn_overline_h_generic_proof}
		\end{equation}
	if $\dim\ker\boldsymbol\lambda|_{r=R}>0$. We have seen in equation \eqref{eqn_dimH_geq3_wlog_q_i_order} that this holds if at least one $q_i$ as in \eqref{eqn_dimgen_q_def} is equal to $\frac{1}{R}$, where we recall our convention $q_1\geq\ldots\geq g_{n-1}$ so that $q_1=\frac{1}{R}$ must hold. It now suffices to show that $R<\sqrt{3}$ and $q_1=\frac{1}{R}$ already imply that the initial CCPSR manifold in standard form $\mathcal{H}$ corresponding to $h=x^3-x\left(\langle v,v\rangle +w^2\right) + C(v) +Q(v)w + aw^3$, cf. equation \eqref{eqn_h_n_geq_3_first_form} for $b=0$, is singular at infinity. From equation \eqref{eqn_pullback_Theta_vvv_Rlesssqrt3_b0} we have concluded that with $C(v)$ and $\widetilde{C}(v)$ as in \eqref{eqn_dimH_geq3_C_splitting_b0}, $\widetilde{C}(v)$ must identically vanish. Thus, in order to show that $\mathcal{H}$ is singular at infinity, it suffices by restricting $h$ to the $s_1$-$w$-plane and using Lemma \ref{lem_sing_at_infty_iff_max_cond} to show that
		\begin{equation}
			\max\limits_{\{s_1^2+w^2=1\}} \left(\frac{1}{R}s_1^2w+aw^3\right)=\frac{2}{3\sqrt{3}},\label{eqn_2dim_2spikes}
		\end{equation}
	where $a=\frac{R^2-1}{R^3}$. One can easily verify that this is true by inserting one of the unit vectors
		\begin{equation}
			(s_1,w)^\mathrm{T}=\left(\pm\frac{\sqrt{3-R^2}}{\sqrt{3}},\frac{R}{\sqrt{3}}\right)^\mathrm{T}.\label{eqn_2dim_2spikes_solution}
		\end{equation}
	Hence, $\mathcal{H}$ is singular at infinity as claimed.

	Next we will have to deal with the cases with $b\ne 0$ studied in Section \ref{section_dimH_geq3_bnon0}. If $\boldsymbol\lambda$ as in \eqref{eqn_dimgen_lambda_def} fulfils $\boldsymbol\lambda|_{r=R}>0$ we have shown that the corresponding limit geometry is isomorphic to $\mathbb{R}_{>0}\ltimes\mathbb{R}^{n-1}$, cf. equation \eqref{eqn_dimgen_lambda_nokernel_hE1_limit}. Hence we do not need to study these cases any further. Next we need to consider the possible cases with $\dim\ker\boldsymbol\lambda|_{r=R}>0$. We will as before assume without loss of generality $q_1\geq\ldots\geq q_{n-1}$. Suppose first that $\frac{1}{R}>q_1$ and $\dim\ker\boldsymbol\lambda|_{r=R}=1$. We have seen that these cases have precisely the requirement \eqref{eqn_ker_0_1_condition} and the possible limit geometries arising are given by \eqref{eqn_limit_poly_ker_0_1}. For our current question it suffices again to show that each of the initial CCPSR manifolds in these cases are singular at infinity. To do so we will use equation \eqref{eqn_C_1_by_1_minus_qR_gamma}. We now restrict $P_3$ as in \eqref{eqn_P3_term_bnon0} to the plane spanned by the Euclidean unit vectors
		\begin{equation*}
			\left\|\frac{1}{1-qR}\gamma\right\|^{-1}\frac{1}{1-qR}\gamma,\quad \partial_w,
		\end{equation*}
	and define
		\begin{equation*}
			\widehat{P_3}\left(\left(\begin{smallmatrix}t\\s\end{smallmatrix}\right)\right):= P_3\left(t\left\|\frac{1}{1-qR}\gamma\right\|^{-1}\frac{1}{1-qR}\gamma +s\partial_w\right).
		\end{equation*}
	Before further studying $\widehat{P_3}\left(\left(\begin{smallmatrix}t\\s\end{smallmatrix}\right)\right)$, recall equation \eqref{eqn_K_def} and observe that $K>0$ by $q<\frac{1}{R}\langle \cdot,\cdot\rangle$. Setting $J:=\frac{1}{\sqrt{K}}$ and using $a=\frac{R^2-1}{R^3}$, we find
		\begin{equation*}
			\widehat{P_3}\left(\left(\begin{smallmatrix}t\\s\end{smallmatrix}\right)\right)= t^3\frac{J(J^2-3)}{3\sqrt{3}} + t^2s\frac{-J^2+1}{R} + ts^2\frac{J(3-R^2)}{\sqrt{3}R^2} + s^3\frac{R^2-1}{R^3}.
		\end{equation*}
	Next, similar to the $2$-dimensional case, we rotate the coordinates $\left(\begin{smallmatrix}t\\s\end{smallmatrix}\right)$ with $M(f)$ as in equation \eqref{eqn_Mf_rotation} and obtain for $f=\frac{\sqrt{3}}{JR}$
		\begin{equation*}
			\widetilde{P_3}\left(\left(\begin{smallmatrix}t\\s\end{smallmatrix}\right)\right):=\left(M(f)^*\widehat{P_3}\right)\left(\left(\begin{smallmatrix}t\\s\end{smallmatrix}\right)\right)=t^2s\frac{\sqrt{3+J^2R^2}}{\sqrt{3}R} + s^3\frac{(-3+3R^2-J^2R^2)\sqrt{3+J^2R^2}}{3\sqrt{3}R^3}.
		\end{equation*}
	The next step is solving $\D \widetilde{P_3}-\frac{2}{\sqrt{3}}(t\,\D t+ s\,\D s)$. We obtain the set of nonzero solutions
		\begin{align}
			(t,s)&\in\left\{\left(0,\frac{2R^3}{(-3+3R^2-J^2R^2)\sqrt{3+J^2R^2}}\right),\right.\notag\\
			&\quad\left.\left(\pm \frac{\sqrt{3+R^2(-1+J^2)}}{\sqrt{3+J^2R^2}},\frac{R}{\sqrt{3+J^2R^2}}\right)\right\}.\label{eqn_JR_solution}
		\end{align}
	It is now straightforward to check that the second two solutions have in fact Euclidean norm $1$ and are furthermore well defined for all $R\in\left(\frac{\sqrt{3}}{2},\sqrt{3}\right)$ and all $J>0$. Hence, the initial CCPSR manifold in standard form is, as claimed, singular at infinity and, hence, of non-regular boundary behaviour. Note that the above result yields more, namely boundaries for $J$ and, hence, $K$. The initial CCPSR manifold being closed implies that the first solution in \eqref{eqn_JR_solution} must have Euclidean norm at least $1$. With $J>0$ one can show that this holds if and only if
		\begin{equation*}
			J\in \left(0,\frac{\sqrt{4R^2-3}}{R}\right]\quad \Leftrightarrow\quad K\in \left[\frac{R^2}{4R^2-3},\infty\right),
		\end{equation*}
	allowing us to obtain a lower bound for $K$ imposed by the initial connected PSR manifold being closed. The extremal case $K=\frac{R^2}{4R^2-3}$ corresponds to an additional ray in the boundary of the cone of the initial CCPSR manifold along which it is singular at infinity.

	It remains to deal with the cases $\dim\ker \boldsymbol\lambda|_{r=R}>0$ and $q_1=\frac{1}{R}$, where we again assume $q_1\geq\ldots\geq q_{n-1}$ since the possible limit polynomials \eqref{eqn_limit_poly_ker_m}, respectively \eqref{eqn_limit_poly_m_plus_1}, are possibly not equivalent to \eqref{eqn_overline_h_generic_proof}. In all of these cases we have shown that the $P_3$-term of the defining polynomial of the initial CCPSR manifold in standard form is of the form \eqref{eqn_P_3_n_geq_3_at_least_m_ker_lambda}. Assume that $\frac{1}{R}=q_1=\ldots=q_m>q_{m+1}$ for some $1\leq m\leq n-1$ and $\dim\ker \boldsymbol\lambda|_{r=R}=m$ or $\dim\ker \boldsymbol\lambda|_{r=R}=m+1$. There is no need to differentiate between the possible two options for the dimensions of $\ker \boldsymbol\lambda|_{r=R}$ as the proofs are completely analogous. After restricting $P_3\left(\left(\begin{smallmatrix}s\\ u\\ w\end{smallmatrix}\right)\right)$ as in \eqref{eqn_P_3_n_geq_3_at_least_m_ker_lambda} to the $(m+1)$-dimensional linear subspace $\{u=0\}\subset\mathbb{R}^m\times\mathbb{R}^{n-1-m}\times\mathbb{R}$, it follows with equation \eqref{eqn_gamma_ker_m_and_m_1_cases} as for the $2$-dimensional case that
		\begin{equation*}
			\max\limits_{\{\langle s,s\rangle + w^2=1\}}P_3\left(\left(\begin{smallmatrix}s\\ 0\\ w\end{smallmatrix}\right)\right)=\frac{2}{3\sqrt{3}},
		\end{equation*}
	cf. \eqref{eqn_2dim_2spikes} and \eqref{eqn_2dim_2spikes_solution}. Alternatively and more explicitly see equation \eqref{eqn_ker_m_m_plus_1_explicit_singatinfpoints} and the comment afterwards. Hence, the initial CCPSR manifold in standard form is singular at infinity and thus of non-regular boundary behaviour. This finishes the proof of Theorem \ref{thm_generic_limit_behaviour}.\hfill\qed

\section{Further applications and outlook}
	In the following we will discuss some applications of our results and give an outlook of what can and should be subject of future studies.
	
	\begin{rem}
		One thing we have not yet answered is what the number $m\in\{0,\ldots,n-1\}$ in Thm. \ref{thm_limit_behaviour_ccpsr} \eqref{thm_limit_behaviour_ccpsr_iii} corresponds to. The answer depends on the value of $R$, that is if $R=\sqrt{3}$ or $R<\sqrt{3}$, where $R$ is the smallest positive zero of $\beta$ \eqref{eqn_beta_def_general}. If $R<\sqrt{3}$, we find that $m=\dim\ker\boldsymbol\lambda|_{r=R}$, cf. equation \eqref{eqn_dimgen_lambda_def}. We might interpret this \textbf{heuristically} as follows: The more $-\partial^2 h$ degenerates when we move to the boundary of the cone spanned by the initial CCPSR manifold, the more the limit geometry is different from the only possible limit geometry for CCPSR manifold with regular boundary behaviour, cf. Theorem \ref{thm_generic_limit_behaviour}, if said infinity does not correspond to a direction with $R=\sqrt{3}$. In the case $R=\sqrt{3}$, that is when we move towards infinity that is non-regular in the boundary behaviour sense, the situation is exactly the opposite. We find that by analysing equations \eqref{eqn_limit_poly_dimHgeq3_b0_Rsqrt3_m0}, \eqref{eqn_limit_poly_dimHgeq3_b0_Rsqrt3_mlessn1}, and \eqref{eqn_limit_poly_dimHgeq3_b0_Rsqrt3_mn1} that the more $-\partial^2 h$ degenerates when we move to the boundary of the cone spanned by the initial CCPSR manifold in the sense of increasing $\dim\ker\boldsymbol\lambda|_{r=R}$, the \textit{closer} our limit geometry is to $\mathbb{R}_{>0}\ltimes\mathbb{R}^{n-1}$. In that case, the integer $m$ in Thm. \ref{thm_limit_behaviour_ccpsr} \eqref{thm_limit_behaviour_ccpsr_iii} fulfils $m=n-1-\dim\ker\boldsymbol\lambda|_{r=R}$. Note that the two formulas for $m$ are consistent with the $2$-dimensional case.
	\end{rem}

	One application of our results lies in quaternionic K\"ahler geometry via the so-called r- and c-map constructions from supergravity theory \cite{ACD,F,DV}. The r-map is a construction that takes a CCPSR manifold $\mathcal{H}$ of dimension $n$ and yields a projective special K\"ahler manifold $M$ of real dimension $2n+2$. The supergravity c-map takes a connected projective special K\"ahler manifold $M$ of real dimension $2n+2$ and yields a quaternionic K\"ahler manifold of negative scalar curvature of real dimension $4n+8$. It has been shown in \cite{CHM} that these construction preserve geodesic completeness, and they were used in \cite{CDJL} to construct an example of a complete non-compact locally inhomogeneous quaternionic K\"ahler manifold of negative scalar curvature of dimension $4n+8$ for all $n\geq 3$, corresponding to the $\mathrm{q}:=\mathrm{c}\circ\mathrm{r}$-map image of the CCPSR manifold in standard form $\mathcal{H}\subset\{h=x^3-x\langle y,y\rangle=1\}$. In \cite{CHM} it was shown that each initial CCPSR manifold $\mathcal{H}$ admits a totally geodesic embedding in its r- and also q-map image. Since the q-map image is on the level of pseudo-Riemannian fibre bundles a principle fibre bundle over the CCPSR manifold, it makes sense to lift our construction of limit geometries to r- and q-map images. For the example constructed in \cite{CDJL}, this means by Theorem \ref{thm_generic_limit_behaviour} the quaternionic K\"ahler manifold $\mathrm{q}(\mathcal{H})$ behaves in the limit of curves leaving every compact subset of the totally geodesically embedded CCPSR manifold with regular boundary behaviour $\mathcal{H}$ like a symmetric quaternionic K\"ahler manifold, cf. \cite{C,DV}, in the sense that its metric tensor converges in the same vein as in Proposition \ref{prop_metric_convergence}. In general it is very difficult to determine if a quaternionic K\"ahler manifold in the image of the q-map is locally inhomogeneous. Our results might be used as follows in future studies. If one can find a closed formula of the Kretschmann scalar of q-map images of $\overline{\mathcal{H}}\subset\{\overline{h}=1\}$, $\overline{h}$ as in Thm. \ref{thm_limit_behaviour_ccpsr} \eqref{thm_limit_behaviour_ccpsr_iii}, at the point corresponding to $\left(\begin{smallmatrix}x\\y\end{smallmatrix}\right)=\left(\begin{smallmatrix}1\\ 0\end{smallmatrix}\right)\in\overline{\mathcal{H}}$ with the help of the recent results in \cite{CST}, and show that these necessarily differ for different values of e.g. $m$ in said formula, one only had to show that the initial CCPSR manifold has different limit geometries yielding different values for the Kretschmann scalar with the help of the determined formula. This endeavour might prove to be very technical, but once completed should easily yield many new examples of non-compact locally inhomogeneous quaternionic K\"ahler manifolds of negative scalar curvature.

	Another open problem related purely to limit geometry of CCPSR manifolds is the following.

	\begin{open}
		Classify all polynomials $\overline{h}$ of the form Thm. \ref{thm_limit_behaviour_ccpsr} \eqref{thm_limit_behaviour_ccpsr_iii} up to equivalence.
	\end{open}

	Solving the above is probably very difficult. By \cite[Prop.\,4.6]{L2} these polynomials include standard forms for all homogeneous spaces. It might be a more realistic ansatz to try solving the following problem.
	
	\begin{open}
		Determine all polynomials $\overline{h}$ of the form Thm. \ref{thm_limit_behaviour_ccpsr} \eqref{thm_limit_behaviour_ccpsr_iii} with corresponding CCPSR manifold in standard form being a homogeneous space.
	\end{open}

	Note that homogeneous CCPSR manifolds have been classified, although not in our standard form, in \cite{DV}, so that reference will most likely prove to be very useful to solve the above problem.

	In Corollary \ref{cor_non_hd} we have seen that the quotient topology of the moduli space of CCPSR manifolds, respectively their defining polynomials, is not Hausdorff. We have shown that a certain class corresponding to the homogeneous space $\mathbb{R}_{>0}\ltimes \mathbb{R}^{n-1}$ can not be separated via disjoint open sets from any class corresponding to CCPSR manifolds with regular boundary behaviour. A natural question one should ask is as follows.

	\begin{open}
		How many distinct points in the moduli space of CCPSR manifolds equipped with the quotient topology that can be separated by disjoint open sets are there, depending on the dimension of the manifolds?
	\end{open}

	Another question that we have not dealt with yet is possible limit geometries of \textbf{incomplete} maximal connected PSR manifolds. For example the maximal connected PSR curve in standard form $\mathcal{H}\subset\{h=x^3-xy^2+y^3=1\}$ has in our sense one direction, that is the negative $y$-direction, in which it makes sense to say that its limit geometry is well defined and given by the CCPSR curve in standard form $\overline{\mathcal{H}}\subset\left\{\overline{h}=x^3-xy^2+\frac{2}{3\sqrt{3}}y^3=1\right\}$. This can be easily checked. However, in the positive $y$-direction, it does not make sense to define a limit geometry since one can show that $\mathcal{H}$ is equivalent to the maximal connected PSR curve in standard form contained in the level set $\{\widetilde{h}=x^3-xy^2+cy^3=1\}$ for all $c>\frac{2}{3\sqrt{3}}$. A reasonable problem to approach would be the following.

	\begin{open}\label{open_incomplete_limits}
		Prove that a maximal incomplete connected PSR manifold in standard form $\mathcal{H}\subset\{h=1\}$ has a CCPSR manifold in standard form $\overline{\mathcal{H}}\subset\{\overline{h}=1\}$ as a limit geometry in the sense that there exists a curve $\gamma:[0,1)\to \mathcal{H}$ along which $h|_t$ analogous to equation \eqref{eqn_h_t} converges to $\overline{h}$ as $t\to 1$, if and only if
			\begin{equation*}
				\overline{\{\mathbb{R}_{>0}\cdot\mathcal{H}\}}\not\subset\{h>0\}.
			\end{equation*}
		The above equation means that we can find a convergent sequence with non-zero limit of points in the cone spanned by $\mathcal{H}$, such that its limit is contained in $\{h=0\}$.
	\end{open}
	
	A counterexample to the existence to a sequence as in Open problem \ref{open_incomplete_limits} in any dimension $n\geq 2$ is the maximal connected incomplete CCPSR manifold in standard form corresponding to
		\begin{equation*}
			h=x^3-x\langle y,y\rangle + r\left(-\frac{2}{\sqrt{3}}(y_1^2+\ldots y_{n-1}^2)y_n^2 + \frac{2}{3\sqrt{3}}y_n^3\right)
		\end{equation*}
	for any $r>1$. For $r=1$, the corresponding maximal connected PSR manifold in standard form is in fact a homogeneous space. We leave it as an exercise for the reader to verify the latter two claims.

	The concept of limit geometry can be readily generalised to so-called \textbf{CCGPSR} (short for \textbf{c}losed \textbf{c}onnected \textbf{g}eneralised \textbf{p}rojective \textbf{s}pecial \textbf{r}eal) manifolds introduced in \cite{L1}. These manifolds correspond to higher homogeneity degrees of the defining polynomial. An open and probably extremely difficult problem is the generalisation of our results to quartic CCGPSR manifolds $\mathcal{H}\subset\{h=1\}$, that is the connected components of level sets in quartic homogeneous polynomials of the form
		\begin{equation*}
			h=x^4-x^2\langle y,y\rangle + xP_3(y) + P_4(y),\quad \left(\begin{smallmatrix}x\\y\end{smallmatrix}\right)=\left(\begin{smallmatrix}1\\ 0\end{smallmatrix}\right)\in\mathcal{H},
		\end{equation*}
	such that $\mathcal{H}$ consists only of hyperbolic points of $h$. Even for quartic CCGPSR curves, which have been classified in \cite[Thm.\,7.2]{L1}, this is a highly non-trivial problem. In \cite[Prop.\,4.4]{L3} we found a complete solution that also includes maximal incomplete connected quartic GPSR curves, which have been classified in that work as well, cf. \cite[Thm.\,3.2]{L3}. However, once successful, it might be useful in proving the still open question whether or not quartic CCGPSR manifolds are automatically geodesically complete.

	Lastly, we ask the following question which relates our results to the theory of the K\"ahler-Ricci flow of compact K\"ahler 3-folds and is probably a good place to start studying a likely and useful relation.
	
	\begin{open}
		Find an explicit example of a compact K\"ahler 3-fold $X$ with $c_1(X)\ne 0$, such that the PSR manifold $\mathcal{H}$ defined by $\{h=1\}\cap\mathcal{K}$, where $h$ is defined as in equation \eqref{eqn_h_KR_flow} and $\mathcal{K}$ is the K\"ahler cone of $X$, is closed and connected and has regular boundary behaviour. Describe the limit of the volume preserving K\"ahler-Ricci flow, or the inverse flow, starting with an initial value such that the solution curve on the level of classes leaves every compact subset of $\mathcal{H}$.
	\end{open}

	The above does not define what \textit{limit} should mean in this context. A possible answer to this might be a certain type of degeneration of the K\"ahler form, the complex structure, or a topological degeneration of $X$ in a meaningful sense. This is a difficult problem but we expect useful and interesting results, in particular for the theory of time-incomplete K\"ahler-Ricci flow.

\end{document}